\pdfoutput=1
\documentclass[a4paper,11pt]{article}
\usepackage[T1]{fontenc}
\usepackage[utf8]{inputenc}
\usepackage{lmodern}
\usepackage[margin=25mm]{geometry}
\usepackage{amsmath,amssymb,amsthm,mathtools}
\usepackage{microtype}
\usepackage{graphicx}
\usepackage{booktabs}
\usepackage{enumitem}
\usepackage[section]{placeins}
\usepackage{tikz}
\usetikzlibrary{arrows.meta,calc}
\usepackage[numbers,sort&compress]{natbib}
\usepackage[hidelinks]{hyperref}
\usepackage{bookmark}

\allowdisplaybreaks[2]
\numberwithin{equation}{section}
\setlist{itemsep=2pt,topsep=4pt}

\newtheorem{theorem}{Theorem}[section]
\newtheorem{proposition}[theorem]{Proposition}
\newtheorem{corollary}[theorem]{Corollary}
\newtheorem{lemma}[theorem]{Lemma}
\theoremstyle{definition}
\newtheorem{definition}[theorem]{Definition}
\newtheorem{example}[theorem]{Example}
\theoremstyle{remark}
\newtheorem{remark}[theorem]{Remark}

\newcommand{\dd}{\mathrm{d}}
\newcommand{\B}{\mathcal{B}}
\newcommand{\minor}{\mathcal{D}}
\newcommand{\Z}{\mathbf{Z}}
\newcommand{\Y}{\mathbf{Y}}
\newcommand{\K}{\boldsymbol{\mathcal{K}}}
\newcommand{\al}{\boldsymbol{\alpha}}
\newcommand{\be}{\boldsymbol{\beta}}
\newcommand{\om}{\boldsymbol{\omega}}
\newcommand{\Pib}{\mathbf{\Pi}}
\newcommand{\Lop}{\boldsymbol{\mathcal{L}}}

\DeclareMathOperator{\spanop}{span}
\DeclareMathOperator{\proj}{proj}
\DeclareMathOperator{\atantwo}{atan2}
\DeclareMathOperator*{\argmin}{arg\,min}

\title{A Poisson-bracket operator for Hamiltonian separability:\\
intrinsic geometry and exact defects}
\author{Vladimir Martinusi\thanks{Corresponding author. E-mail:
\href{mailto:vmartinusi@technion.ac.il}{vmartinusi@technion.ac.il}.
ORCID: \href{https://orcid.org/0000-0001-9689-7657}{0000-0001-9689-7657}.}\\[1ex]
\small Stephen B.\ Klein Faculty of Aerospace Engineering,\\
\small Technion -- Israel Institute of Technology, Technion City, Haifa 3200003, Israel}
\date{}

\begin{document}
\maketitle

\begin{abstract}
	Additive separation of the Hamilton--Jacobi equation depends on the
	Hamiltonian and on the canonical variables, and this paper introduces a
	five-input Poisson-bracket operator whose values on canonical pairs
	reproduce the Levi-Civita separability expressions. On a
	symplectic product structure its off-diagonal values are intrinsic, since
	they equal a symplectic pairing of the sector Hamiltonian fields and measure
	the Frobenius obstruction of a Lagrangian distribution, and their
	vanishing yields a local separation theorem. Because the
	operator is cubic in the Hamiltonian, its perturbation expansion is exact
	and defines a linearised separability operator and a defect functional. For
	Hamiltonians with a cyclic angle, the linearised operator divided by the
	squared rotor frequency is an angle derivative of the Melnikov integrand,
	so that the defect integrated along a homoclinic orbit decides the
	first-order splitting of the separatrix. Closed-form defects are obtained
	for the polar-nodal normal forms of satellite theory and, through a
	covariance of the operator on level sets, for Kerr geodesics. Used
	for coordinate selection, the defect singles out the separating focal
	distance of the Euler two-centre problem as its unique nondegenerate
	minimum, and for axisymmetric potentials it shows that the prescriptions of
	Sanders and Vasiliev are two weightings of the same St\"ackel residual. The
	residual selector and Binney's shell-orbit fit improve both actions over
	the median prescription for 48 and 50 of 50 orbits in a flattened
	logarithmic potential, whereas the operator selector lowers the
	90th-percentile action spread of 50 orbits in a Miyamoto--Nagai potential
	from $2.90\%$ to $0.91\%$.
\end{abstract}


\section{Introduction}
\label{sec:intro}

The existence of an additively separated complete integral of the
stationary Hamilton--Jacobi equation
$H(\mathbf q,\partial_{\mathbf q}\mathcal S)=E$ depends both on the Hamiltonian and on the canonical
variables in which it is written. Levi-Civita \cite{LeviCivita} reduced this
question to a finite set of differential identities, one for every pair of
distinct degrees of freedom, each of which requires a cubic expression in the
first and second derivatives of $H$ to vanish. Llibre and Ram\'{\i}rez
\cite{LlibreRamirez2020} have rewritten these identities through the
partial Hamiltonian fields of the degrees of freedom and related them to the
Frobenius integrability of the distribution that these fields span. For
natural Hamiltonians,
St\"ackel's construction \cite{Stackel} and its geometric development through
Killing tensors and separable webs
\cite{Eisenhart,KalninsMiller,Benenti1997,Benenti2016} describe all separable
coordinate systems, and Waksj\"o and Rauch-Wojciechowski have turned this
description into an algorithm that finds them \cite{Waksjo}. In a different
direction, bi-Hamiltonian geometry characterises separable foliations through
invariant Lagrangian distributions \cite{FalquiPedroni,Blaszak}, while
Nijenhuis and Haantjes geometry clarifies the role played by compatible
endomorphisms \cite{BolsinovKonyaevMatveev,TondoTempesta,TempestaTondo2022}.

All of these theories answer an exact, yes-or-no question. Two related
questions are less well served: (i) how far a given Hamiltonian is from
separating in a given coordinate system, measured in a way that does not
depend on inessential choices of coordinates, and (ii) which member of a
family of trial coordinates should be used when no member separates exactly.
The second question arises directly in practice. The St\"ackel fudge
\cite{Binney2012,SandersBinney2016}, for instance, estimates actions in a
non-separable axisymmetric galactic potential by applying the quadratures of a
locally separated model in prolate spheroidal coordinates, and its accuracy
depends on the choice of a single focal distance
\cite{Sanders2012,Binney2014,VasilievAgama}.

This paper develops a scalar calculus that addresses both questions. Its
basic object is a five-input operator $\B_{a,b;c,d}[H]$, built from Poisson
brackets of four observables with the Hamiltonian, which reproduces the
Levi-Civita expressions on canonical pairs. On general observables, the
operator is represented by a decomposable two-form and a decomposable
bivector, so that its algebra, its covariance and its kernels can all be
stated invariantly.

\paragraph{Contributions.} On a symplectic product structure, the values of
the operator on canonical pairs are intrinsic, and since the operator is
cubic in the Hamiltonian, its perturbation expansion is exact and defines a
separability defect in closed form. Along a homoclinic orbit of a
Hamiltonian with a cyclic angle, the integrated defect is the derivative of
the Melnikov function, and over a family of trial coordinates the defect is
an objective for coordinate selection. The results are grouped below under
these three headings, and each item states which of its parts are new and
which are known.

\noindent\emph{Operator and intrinsic geometry.}
\begin{enumerate}[label=(\roman*)]
		\item \emph{New.} The five-input operator itself, which extends to
	arbitrary observables the form of the Levi-Civita conditions
	written by Llibre and Ram\'{\i}rez through the sector fields
	\cite[eqs.~(13)--(15)]{LlibreRamirez2020}, together with its structure
	theorem (Proposition~\ref{prop:structure}), its exact law under energy
	reparametrisation (Proposition~\ref{prop:algebra}), whose restriction to
	canonical pairs is the identity proved in
	\cite[Proposition~5]{LlibreRamirez2020}, its interpretation as a
	projected turning rate (Proposition~\ref{prop:turning}), and its
	identification, on the angles of an integrable Hamiltonian, with the
	isoenergetic nondegeneracy determinant of Arnold
	\cite{ArnoldSmallDen1963} (Remark~\ref{rem:arnold}).
	\item \emph{New formulas for a known condition.} For a symplectic product
	structure, the off-diagonal values are shown to equal a mixed Hessian and a
	symplectic pairing of the sector Hamiltonian fields
	(Proposition~\ref{prop:intrinsic}), which explains both their symmetry and
	their invariance under canonical changes within the factors.
	\item \emph{Sharpening of known criteria.} Llibre and Ram\'{\i}rez
	\cite[Theorem~10 and Remark~11]{LlibreRamirez2020} showed that, where every $H_{p_i}$ or every
	$H_{q_i}$ is nonzero, the Levi-Civita conditions are equivalent to the
	involutivity of the distribution spanned by the sector fields, and Falqui
	and Pedroni \cite{FalquiPedroni} characterised separation through the
	involutivity of an invariant Lagrangian distribution. Here the Frobenius
	obstruction of the distribution spanned by the sector fields, which is
	Lagrangian, is expressed by a single
	invariant scalar for each pair of degrees of freedom
	(Proposition~\ref{prop:frobenius}), and the separation theorem is proved
	under the sole requirement that the sector fields do not vanish, which
	includes the turning points of the separated motions
	(Theorem~\ref{thm:sector}). A symplectic--Nijenhuis form of the criterion
	is also given (Proposition~\ref{prop:nijenhuis}).
\end{enumerate}

\noindent\emph{Exact defects and homoclinic splitting.}
\begin{enumerate}[label=(\roman*),resume]
	\item \emph{New.} An exact cubic perturbation expansion, a linearised
	separability operator with explicit kernel directions, and a weighted
	defect functional built on them (Sections~\ref{sec:linearised}
	and~\ref{sec:defect-functional}).
	\item \emph{New.} For Hamiltonians with a cyclic angle, an exact identity
	between the linearised operator and the Melnikov integrand of Holmes and
	Marsden \cite{HolmesMarsden1982}, by which the integrated defect along a
	homoclinic orbit is the derivative of the Melnikov function and decides the
	first-order splitting of the separatrix (Proposition~\ref{prop:rotor},
	Theorem~\ref{thm:melnikov} and Corollary~\ref{cor:dichotomy}).
	\item \emph{New applications.} A closed-form defect for the polar-nodal
	normal forms of satellite theory, in which the defect and the Melnikov
	integrand share one angular factor (Proposition~\ref{prop:polar-nodal}),
	and a covariance law of the operator on level sets that gives the defect of
	Kerr geodesics and of their deformations in closed form
	(Lemma~\ref{lem:levelset} and Proposition~\ref{prop:kerr}).
\end{enumerate}

\noindent\emph{Coordinate selection.}
\begin{enumerate}[label=(\roman*),resume]
	\item \emph{New.} For the Euler two-centre problem, an exact factorisation
	of the integrated defect, which identifies the separating focal distance on
	any admissible fixed physical domain (Theorem~\ref{thm:euler}).
	\item \emph{New derivation of known prescriptions.} For axisymmetric
	potentials, the operator factorises into the St\"ackel residual and the
	product of the two prolate coordinate velocities. Its pointwise zero is the local prescription of
	Sanders \cite{Sanders2012,SandersBinney2016}, and removing the weight yields
	a closed-form implementation of the residual criterion of Vasiliev
	\cite{VasilievDoc}, while the shell-orbit fit of Binney \cite{Binney2014}
	corresponds to the operator objective evaluated on the shell orbit
	(Section~\ref{sec:stackel}).
	\item \emph{New comparison.} A study of four focal selectors, among them
	Binney's shell-orbit fit, on 100 orbits in two potentials, with an accuracy
	measure that requires no reference actions, a consistency check against
	independent reference actions, a relative residual defect that ranks the
	orbits a priori, and a turning-point diagnostic that quantifies how
	differently the two least-squares objectives weight an orbit
	(Sections~\ref{sec:numerics}--\ref{sec:weighting}).
\end{enumerate}

The operator and the theorems of
Sections~\ref{sec:operator}--\ref{sec:defects} apply to arbitrary smooth
Hamiltonians, and Propositions~\ref{prop:polar-nodal} and~\ref{prop:kerr}
compute the defect in closed form for the polar-nodal normal forms of
satellite theory, which are not natural, and for the geodesics of the Kerr
metric. The coordinate-selection applications of
Sections~\ref{sec:euler}--\ref{sec:weighting} concern natural Hamiltonians,
for which algebraic criteria \cite{Waksjo,Benenti2016} already identify
exactly separating coordinates without any optimisation, so the variational
formulation of these sections is meant for quantitative comparison and,
above all, for the approximate case, in which no member of the trial family
separates.

\paragraph{Outline.} Section~\ref{sec:operator} fixes the notation, defines
the operator and establishes its structure, and Section~\ref{sec:product}
develops the geometry of product structures together with the separation
theorem. Section~\ref{sec:defects} contains the perturbation theory, the defect
functional, the relation between the defect and homoclinic splitting, and
the covariance on level sets with its application to Kerr geodesics. The
defect functional is then used in Section~\ref{sec:euler} for exact
coordinate identification. Sections~\ref{sec:stackel}--\ref{sec:weighting}
turn to the St\"ackel application and its numerical study, and Sections~\ref{sec:limitations} and~\ref{sec:conclusions} contain the
limitations and the conclusions, respectively. The reference-action protocol used for the consistency check is documented
in Appendix~\ref{app:reference}.

\section{The five-input operator}
\label{sec:operator}

This section introduces the operator and the algebraic facts about it on
which the rest of the paper relies. Since the operator mixes Poisson brackets,
exterior forms and bivectors, the section begins by fixing the notation and the sign
and normalisation conventions.

\subsection{Notation and conventions}
\label{sec:notation}

Throughout, $M$ is a smooth symplectic manifold of dimension $2n$ with
symplectic form $\om$ and Poisson bracket $\{\cdot,\cdot\}$, all functions are assumed to be
smooth, and $C^\infty(M)$ denotes the smooth functions on $M$.

\paragraph{Partial derivatives.} In a canonical chart
$(q_1,\ldots,q_n,p_1,\ldots,p_n)$, subscripts $q_i$ and $p_i$ on a function
denote partial derivatives, so that $H_{q_i}=\partial H/\partial q_i$ and
$H_{p_i}=\partial H/\partial p_i$, and repeated subscripts denote higher
derivatives, for example $H_{q_ip_j}=\partial^2H/\partial q_i\partial p_j$.
The same convention is used for other coordinates, such as
$\Phi_R=\partial\Phi/\partial R$ and $\Phi_{Rz}=\partial^2\Phi/\partial R\,\partial z$
in Section~\ref{sec:stackel}. Only coordinate names act as derivative
subscripts, and only on functions such as $H$, $\Phi$ or $D$. Conjugate
momenta such as $p_x$, $p_R$ or $p_u$ are coordinates in their own right,
and all other subscripts, as in $q_i$, $b_{ij}$, $\Y_i$, $R_0$ or
$\Phi_{\log}$, are labels.

\paragraph{Wedge products.} For one-forms $\theta_1,\ldots,\theta_k$, the
wedge product $\theta_1\wedge\cdots\wedge\theta_k$ is the alternating
$k$-form whose value on vectors $\mathbf X_1,\ldots,\mathbf X_k$ is the
determinant $\det\bigl(\theta_i(\mathbf X_j)\bigr)$; in particular, for two
one-forms,
\begin{equation}
 (\theta\wedge\phi)(\mathbf X,\mathbf V)
 =\theta(\mathbf X)\,\phi(\mathbf V)-\theta(\mathbf V)\,\phi(\mathbf X).
 \label{eq:wedge}
\end{equation}
Dually, the wedge product of two vector fields is the bivector that acts on
pairs of one-forms by
$(\mathbf X\wedge\mathbf V)(\theta,\phi)=\theta(\mathbf X)\,\phi(\mathbf V)-\theta(\mathbf V)\,\phi(\mathbf X)$.
This is the determinant normalisation used, for instance, in \cite{Arnold}.
Part of the differential-geometry literature instead defines the wedge
product as the antisymmetrisation of the tensor product, which inserts a
factor $\tfrac12$ in \eqref{eq:wedge} and, more generally, a factor $1/k!$
for $k$ one-forms. The two conventions differ only by these constants, but
only the determinant normalisation gives, for $\om=\sum_i\dd q_i\wedge\dd p_i$,
the relations $\om(\partial_{q_i},\partial_{p_i})=1$ and
$\iota_{\Z_F}\om=\dd F$, as well as the invariant formula
$\dd\theta(\mathbf X,\mathbf V)=\mathbf X\bigl(\theta(\mathbf V)\bigr)-\mathbf V\bigl(\theta(\mathbf X)\bigr)-\theta([\mathbf X,\mathbf V])$
used in Section~\ref{sec:product}, without extra factors. With the factor
$\tfrac12$, the representation $\B_{a,b;c,d}[H]=\be(\Z_a,\Z_b)$ of
Proposition~\ref{prop:structure}, for example, would read
$\B_{a,b;c,d}[H]=2\,\be(\Z_a,\Z_b)$. A two-form or bivector is called \emph{decomposable} if it is the wedge
product of two factors; its exterior square then vanishes, and the algebraic
normal form of skew forms \cite[Chapter~I, Theorem~1.7]{BryantEtAl1991},
which holds for bivectors with the vector space replaced by its dual, shows
that it has rank zero or two.

\paragraph{Other notations.} The symbol $\dd$ denotes the exterior
derivative, $\iota_{\mathbf X}$ the contraction with a vector field
$\mathbf X$, $[\mathbf X,\mathbf V]$ the Lie bracket of vector fields,
$\Gamma(E)$ the smooth sections of a vector bundle $E$, and $\Phi^*$ the
pullback by a map $\Phi$; a symplectomorphism is a diffeomorphism with $\Phi^*\om=\om$.
Finally, the Poisson bracket and the Hamiltonian vector field $\Z_F$ of a
function $F$ are fixed by
\begin{equation}
 \{F,G\}=\sum_{i=1}^{n}\bigl(F_{q_i}G_{p_i}-F_{p_i}G_{q_i}\bigr),
 \qquad \Z_F(G)=\{G,F\},
 \qquad \om=\sum_{i=1}^n\dd q_i\wedge\dd p_i,
 \label{eq:convention}
\end{equation}
so that $\iota_{\Z_F}\om=\dd F$ and $\Z_H$ generates the Hamiltonian flow:
along a trajectory, $\dot G=\{G,H\}=\Z_H(G)$ for every function $G$.

\subsection{Definition and structure}

\begin{definition}[The five-input operator]
	\label{def:B}
	For $a,b,c,d,H\in C^\infty(M)$ put
	\begin{equation}
		\minor_{a,b}(F,G)=\{a,F\}\{b,G\}-\{a,G\}\{b,F\},\qquad
		\Y^H_{a,b}(F)=\minor_{a,b}(F,H),
		\label{eq:minor}
	\end{equation}
	and define the response rates $u=\{H,c\}$ and $v=\{H,d\}$. The five-input
	operator is
	\begin{equation}
		\B_{a,b;c,d}[H]=u\,\Y^H_{a,b}(v)-v\,\Y^H_{a,b}(u).
		\label{eq:B}
	\end{equation}
\end{definition}

The functions $a$ and $b$, referred to as the \emph{drivers}, enter only
through the derivation $\Y^H_{a,b}$, while the \emph{responses} $c$ and $d$
enter through their rates of change $u=-\dot c$ and $v=-\dot d$ along the
flow. Any four smooth functions can serve as probes -- canonical coordinates
are the natural choice, but nothing in the definition depends on them.

\begin{proposition}[Structure theorem]
	\label{prop:structure}
	For given $H$, $c$ and $d$, let
	\begin{equation}
		\al=u\,\dd v-v\,\dd u,\qquad \be=\al\wedge\dd H,\qquad
		\mathbf W=u\Z_v-v\Z_u,\qquad \K=\mathbf W\wedge\Z_H.
		\label{eq:forms}
	\end{equation}
	Then the following hold for all $a,b\in C^\infty(M)$.
	\begin{enumerate}[label=(\alph*)]
		\item \emph{Representations.}
		\begin{equation}
			\B_{a,b;c,d}[H]=\al\bigl(\Y^H_{a,b}\bigr)=\be(\Z_a,\Z_b)=\K(\dd a,\dd b).
			\label{eq:representations}
		\end{equation}
		\item \emph{Rank and exterior derivatives.} The two-form $\be$ and the
		bivector $\K$ are decomposable and therefore have rank zero or two, and
		\[
		\dd\al=2\,\dd u\wedge\dd v,\qquad \dd\be=2\,\dd u\wedge\dd v\wedge\dd H.
		\]
		\item \emph{Ratio form.} Where $u\ne0$, the ratio $r=u^{-1}v$ of the response
		rates satisfies
		\begin{equation}
			\B_{a,b;c,d}[H]=u^2\,\Y^H_{a,b}(r),\qquad
			\be=u^2\,\dd r\wedge\dd H,\qquad \K=u^2\,\Z_r\wedge\Z_H.
			\label{eq:ratio}
		\end{equation}
		\item \emph{Angle form.} Where $u^2+v^2\ne0$, the polar angle
		$\vartheta=\atantwo(v,u)$ of the vector $(u,v)$, taken on any continuous
		local branch, satisfies
		\begin{equation}
			\B_{a,b;c,d}[H]=(u^2+v^2)\,\Y^H_{a,b}(\vartheta).
			\label{eq:angle}
		\end{equation}
		\item \emph{Bordered Hessian.} With the derivations $D_x=\{x,\cdot\}$,
		\begin{equation}
			\B_{a,b;c,d}[H]=\det\begin{pmatrix}
				0 & D_cH & D_dH\\
				D_aH & D_aD_cH & D_aD_dH\\
				D_bH & D_bD_cH & D_bD_dH
			\end{pmatrix}.
			\label{eq:bordered}
		\end{equation}
	\end{enumerate}
\end{proposition}

\begin{proof}
	\emph{(a)} Since $\Y^H_{a,b}$ is a derivation, evaluating $\al$ on it gives
	\[
	\al\bigl(\Y^H_{a,b}\bigr)=u\,\Y^H_{a,b}(v)-v\,\Y^H_{a,b}(u),
	\]
	which is the definition \eqref{eq:B}. For the other two equalities, the
	conventions \eqref{eq:convention} give
	\[
	\al(\Z_a)=u\{v,a\}-v\{u,a\}=-\mathbf W(a),\qquad
	\dd H(\Z_a)=\{H,a\}=-\Z_H(a),
	\]
	and the wedge rule \eqref{eq:wedge}, together with its dual for bivectors,
	reduces both exterior evaluations to the same expression,
	\[
	\be(\Z_a,\Z_b)=\K(\dd a,\dd b)
	=\mathbf W(a)\,\Z_H(b)-\mathbf W(b)\,\Z_H(a).
	\]
	Expanding $\mathbf W(a)=u\{a,v\}-v\{a,u\}$ and $\Z_H(b)=\{b,H\}$ and
	grouping the terms with \eqref{eq:minor} gives
	\[
	\mathbf W(a)\,\Z_H(b)-\mathbf W(b)\,\Z_H(a)
	=u\,\Y^H_{a,b}(v)-v\,\Y^H_{a,b}(u),
	\]
	which is again \eqref{eq:B}.
	
	\emph{(b)} By \eqref{eq:forms}, $\be$ and $\K$ are wedge products of two
	factors, so they are decomposable and have rank zero or two. Exterior
	differentiation of $\al$ gives
	$\dd\al=\dd u\wedge\dd v-\dd v\wedge\dd u=2\,\dd u\wedge\dd v$, and since
	$\dd(\dd H)=0$, $\dd\be=\dd\al\wedge\dd H=2\,\dd u\wedge\dd v\wedge\dd H$.
	
	\emph{(c)} Substituting $v=ur$ and using the product rules
	$\dd(ur)=r\,\dd u+u\,\dd r$ and $\Z_{ur}=u\Z_r+r\Z_u$ gives
	\[
	\al=u^2\,\dd r,\qquad \mathbf W=u^2\,\Z_r,
	\]
	and inserting these into \eqref{eq:forms} and \eqref{eq:representations}
	yields \eqref{eq:ratio}.
	
	\emph{(d)} On every continuous branch, the differential of the
	two-argument arctangent is
	\[
	\dd\vartheta=\frac{u\,\dd v-v\,\dd u}{u^2+v^2}=\frac{\al}{u^2+v^2},
	\]
	which does not depend on the branch, since two branches differ by a
	constant multiple of $2\pi$. With (a), this gives
	\[
	\B_{a,b;c,d}[H]=\al\bigl(\Y^H_{a,b}\bigr)
	=(u^2+v^2)\,\dd\vartheta\bigl(\Y^H_{a,b}\bigr)
	=(u^2+v^2)\,\Y^H_{a,b}(\vartheta),
	\]
	which is \eqref{eq:angle}.
	
	\emph{(e)} By \eqref{eq:convention}, $D_cH=-u$, $D_dH=-v$,
	$D_aD_cH=-\{a,u\}$ and $D_aD_dH=-\{a,v\}$, and similarly for $b$. A change
	of sign of the last two columns leaves the determinant unchanged, and the
	expansion along the first row gives
	\[
	\det\begin{pmatrix}
		0 & u & v\\
		\{a,H\} & \{a,u\} & \{a,v\}\\
		\{b,H\} & \{b,u\} & \{b,v\}
	\end{pmatrix}
	=u\bigl(\{a,v\}\{b,H\}-\{a,H\}\{b,v\}\bigr)
	-v\bigl(\{a,u\}\{b,H\}-\{a,H\}\{b,u\}\bigr),
	\]
	which is $u\,\Y^H_{a,b}(v)-v\,\Y^H_{a,b}(u)$ by \eqref{eq:minor}.
\end{proof}

Equation~\eqref{eq:angle} gives the operator a simple kinematic meaning. Up
to sign, $(u,v)$ is the velocity of the response pair $(c,d)$ along the flow,
and $\B$ is its squared speed $u^2+v^2$ multiplied by the rate at which the
driver derivation $\Y^H_{a,b}$ turns its direction. This turning picture
recurs throughout the paper, most concretely in Section~\ref{sec:product},
where it describes how motion in one degree of freedom turns the phase
velocity of another.

\begin{proposition}[Algebraic properties]
	\label{prop:algebra}
	The five-input operator introduced in Definition~\ref{def:B} has the
	following properties, for all $a,a_1,a_2,b,c,d,H\in C^\infty(M)$.
	\begin{enumerate}[label=(\alph*)]
		\item \emph{Multilinearity.} $\B_{a,b;c,d}[H]$ is $\mathbb R$-linear in
		each of $a$, $b$, $c$ and $d$.
		\item \emph{Alternation in the drivers.}
		$\B_{b,a;c,d}[H]=-\B_{a,b;c,d}[H]$.
		\item \emph{Alternation in the responses.}
		$\B_{a,b;d,c}[H]=-\B_{a,b;c,d}[H]$.
		\item \emph{Leibniz rule in the drivers.}
		$\B_{a_1a_2,b;c,d}[H]=a_1\,\B_{a_2,b;c,d}[H]+a_2\,\B_{a_1,b;c,d}[H]$,
		and likewise in $b$.
		\item \emph{Cubic dependence on the Hamiltonian.} $\B_{a,b;c,d}[H]$ is a
		homogeneous cubic differential polynomial in the first and second
		derivatives of $H$; in particular,
		$\B_{a,b;c,d}[tH]=t^3\,\B_{a,b;c,d}[H]$ for every $t\in\mathbb R$.
		\item \emph{Energy reparametrisation.} For every smooth function $\chi$ of
		one variable,
		\begin{equation}
			\B_{a,b;c,d}[\chi(H)]=\chi'(H)^3\,\B_{a,b;c,d}[H].
			\label{eq:energychange}
		\end{equation}

		\item \emph{Quadratic form in the first brackets.} Let
		$\mathbf g=\bigl(\{a,H\},\{b,H\},\{c,H\},\{d,H\}\bigr)^{\mathsf T}$ be the
		tuple of first brackets of $H$ with the four probes, and let $\mathbf h$
		be the $2\times2$ matrix of double brackets $h_{xy}=\{x,\{y,H\}\}$, with
		rows indexed by the drivers $x\in\{a,b\}$ and columns indexed by the
		responses $y\in\{c,d\}$. Then
		\begin{equation}
			\B_{a,b;c,d}[H]=\mathbf g^{\mathsf T}\mathbf Q\,\mathbf g,\qquad
			\mathbf Q=\frac12\begin{pmatrix}
				0 & 0 & -h_{bd} & h_{bc}\\
				0 & 0 & h_{ad} & -h_{ac}\\
				-h_{bd} & h_{ad} & 0 & 0\\
				h_{bc} & -h_{ac} & 0 & 0
			\end{pmatrix}
			=-\frac12\begin{pmatrix}
				\mathbf 0 & \operatorname{cof}\mathbf h\\
				(\operatorname{cof}\mathbf h)^{\mathsf T} & \mathbf 0
			\end{pmatrix},
			\label{eq:quadratic}
		\end{equation}
		where $\operatorname{cof}\mathbf h$ is the cofactor matrix of
		$\mathbf h$, so that the driver--driver and response--response blocks of
		$\mathbf Q$ vanish and only mixed double brackets enter.
	\end{enumerate}
\end{proposition}

\begin{proof}
	The seven properties are verified in turn.
	\begin{enumerate}[label=(\alph*)]
		\item The Poisson bracket is bilinear, so the minor \eqref{eq:minor} and
		hence the derivation $\Y^H_{a,b}$ are linear in $a$ and in $b$. The
		response rates $u=\{H,c\}$ and $v=\{H,d\}$ are linear in $c$ and in $d$,
		and \eqref{eq:B} is bilinear in $(u,v)$.
		\item By \eqref{eq:minor}, $\minor_{b,a}(F,G)=-\minor_{a,b}(F,G)$, so
		$\Y^H_{b,a}=-\Y^H_{a,b}$.
		\item Exchanging $c$ and $d$ exchanges $u$ and $v$, which reverses the sign
		of \eqref{eq:B}.
		\item Since $\{a_1a_2,F\}=a_1\{a_2,F\}+a_2\{a_1,F\}$ for every $F$, the
		derivation satisfies
		\[
		\Y^H_{a_1a_2,b}=a_1\,\Y^H_{a_2,b}+a_2\,\Y^H_{a_1,b},
		\]
		and \eqref{eq:B} depends linearly on the derivation. The same argument
		applies to $b$.
		\item After the derivation is expanded, every term of \eqref{eq:B} is a
		product of three factors: a response rate $u$ or $v$, a bracket $\{a,H\}$
		or $\{b,H\}$, and a bracket of a driver with $u$ or $v$. Each factor is
		linear in the first and second derivatives of $H$, which gives the cubic
		homogeneity and, in particular, the scaling by $t^3$.
		\item Writing $s=\chi'(H)$, the chain rule gives $\{\chi(H),F\}=s\{H,F\}$,
		so that
		\[
		u\mapsto su,\qquad v\mapsto sv,\qquad \Y^{\chi(H)}_{a,b}=s\,\Y^H_{a,b}.
		\]
		Substituting into \eqref{eq:B} and expanding with the Leibniz rule for the
		derivation gives
		\[
		\B_{a,b;c,d}[\chi(H)]
		=su\,s\Y^H_{a,b}(sv)-sv\,s\Y^H_{a,b}(su)
		=s^3\,\B_{a,b;c,d}[H]+s^2uv\,\Y^H_{a,b}(s)-s^2vu\,\Y^H_{a,b}(s),
		\]
		and the last two terms cancel.
		\item With $g_x=\{x,H\}$, the identities $D_xH=g_x$ and
		$D_xD_yH=h_{xy}$ turn the bordered form \eqref{eq:bordered} into the
		matrix $\mathbf h$ bordered by $(g_c,g_d)$ and $(g_a,g_b)$. The expansion
		along the first row, followed by the expansion of each $2\times2$ minor
		along its first column, gives
		\[
		\B_{a,b;c,d}[H]
		=-g_ag_c\,h_{bd}+g_ag_d\,h_{bc}+g_bg_c\,h_{ad}-g_bg_d\,h_{ac},
		\]
		in which the coefficient of $g_xg_y$, for a driver $x$ and a response
		$y$, is minus the $(x,y)$ entry of $\operatorname{cof}\mathbf h$. This is
		the expansion of $\mathbf g^{\mathsf T}\mathbf Q\,\mathbf g$, and each
		term is quadratic in the first brackets and linear in the double
		brackets, which makes the cubic structure of (e) explicit. \qedhere
	\end{enumerate}
\end{proof}

The next result shows that the operator respects the symplectic structure
and describes exactly when it vanishes for every choice of drivers.

\begin{proposition}[Covariance and kernel]
\label{prop:covariance}
The operator is covariant under symplectomorphisms: if $\Phi$ is a
symplectomorphism, then
\begin{equation}
 \B_{\Phi^*a,\Phi^*b;\Phi^*c,\Phi^*d}[\Phi^*H]=\Phi^*\B_{a,b;c,d}[H].
 \label{eq:covariance}
\end{equation}
Moreover, let $U\subset M$ be an open set on which $u\ne0$ and $\dd H\ne0$.
Then the following two statements are equivalent:
\begin{enumerate}[label=(\alph*)]
 \item $\B_{a,b;c,d}[H]=0$ on $U$ for all drivers $a,b\in C^\infty(M)$;
 \item the ratio $r=u^{-1}v$ is locally a function of the energy, that is, about
 each point of $U$ there is a smooth function $F$ of one variable such that
 $r=F(H)$.
\end{enumerate}
\end{proposition}

\begin{proof}
	\emph{Covariance.} A symplectomorphism preserves Poisson brackets,
	$\{\Phi^*F,\Phi^*G\}=\Phi^*\{F,G\}$ \cite{Arnold}, and every ingredient of
	\eqref{eq:B}, namely the minor, the derivation and the response rates, is
	built from Poisson brackets. Applying this identity to each of them gives
	\eqref{eq:covariance}.
	
	\emph{Equivalence.} By \eqref{eq:representations} and \eqref{eq:ratio}, on
	$U$ the operator is represented by
	\[
	\B_{a,b;c,d}[H]=\K(\dd a,\dd b),\qquad \K=u^2\,\Z_r\wedge\Z_H .
	\]
	
	\emph{(a)$\Rightarrow$(b).} Suppose that (a) holds. Since the differentials
	of functions span every cotangent space, $\K$ vanishes on all pairs of
	one-forms, so $\K=0$ on $U$, and since $u\ne0$ this gives
	$\Z_r\wedge\Z_H=0$. The map $\dd F\mapsto\Z_F$ is a linear isomorphism by
	the nondegeneracy of $\om$, so it follows that $\dd r\wedge\dd H=0$. Because
	$\dd H\ne0$, every point of $U$ has a neighbourhood with box-shaped local
	coordinates $(y_1,\ldots,y_{2n})$ in which $y_1=H$. In these coordinates
	\[
	\dd r\wedge\dd H=\sum_{k=2}^{2n}\frac{\partial r}{\partial y_k}\,\dd y_k\wedge\dd y_1=0
	\]
	forces $\partial r/\partial y_k=0$ for $k\ge2$, so $r$ depends on $y_1$
	alone, that is, $r=F(H)$ for a smooth function $F$ of one variable.
	
	\emph{(b)$\Rightarrow$(a).} Suppose that $r=F(H)$ near a point of $U$. Then
	\[
	\dd r=F'(H)\,\dd H,\qquad \Z_r=F'(H)\,\Z_H,\qquad \Z_r\wedge\Z_H=0,
	\]
	so $\K=0$ and $\B_{a,b;c,d}[H]=\K(\dd a,\dd b)=0$ near that point for all
	drivers $a,b$. Since the point was arbitrary, (a) holds on $U$.
\end{proof}

For given $H,c,d$, the operator depends on the drivers only through their
differentials, since $\B_{a,b;c,d}[H]=\K(\dd a,\dd b)$ by
\eqref{eq:representations}. Its values on pairs of coordinate functions are
therefore the components of the bivector $\K$, and they transform
tensorially: under a change of coordinates $x\mapsto x'$,
\[
\B_{x'^\beta,x'^\gamma;c,d}[H]
=\sum_{\alpha,\delta}\frac{\partial x'^\beta}{\partial x^\alpha}\,
\frac{\partial x'^\gamma}{\partial x^\delta}\,\B_{x^\alpha,x^\delta;c,d}[H],
\]
so only the Jacobian of the change enters. Dependence on the responses
is different: the functions $c$ and $d$ enter through $u=\{H,c\}$ and
$v=\{H,d\}$, which are differentiated once more by $\Y^H_{a,b}$, so second
derivatives of the responses appear. For a product response, the product
rule gives
\begin{equation}
	\B_{a,b;fc,d}[H]=f\B_{a,b;c,d}[H]+c\B_{a,b;f,d}[H]
	-v\bigl(\{H,c\}\Y^H_{a,b}(f)+\{H,f\}\Y^H_{a,b}(c)\bigr).
	\label{eq:lastslot}
\end{equation}
For a tensorial dependence on $c$ the last term would be absent, and it is
this term that makes the full four-probe array depend on the chosen response
coordinates. For the Hamiltonian $H=p$ on $T^*\mathbb R$, every evaluation of the operator on the probes $q,p$ vanishes, whereas in the
canonical chart
\[
Q=e^q,\qquad P=e^{-q}p,
\]
in which the same Hamiltonian reads $H=QP$, one finds
\[
\B_{Q,P;Q,P}[QP]=2QP.
\]
A dependence of this kind is to be expected, since separability is a
property of the Hamiltonian together with the chosen coordinates.
Section~\ref{sec:product} identifies the part of the array that is invariant
under the transformations preserving separation.

\begin{corollary}[Induced Poisson bracket]
	\label{cor:jacobi}
	Let $H,c,d$ be given, let $U$ be an open set on which $u\ne0$ and
	$\dd r\wedge\dd H\ne0$, and put
	\[
	[a,b]_{\K}=\B_{a,b;c,d}[H],\qquad s=\{r,H\}=\dot r .
	\]
	The bracket $[\cdot,\cdot]_{\K}$ is skew and satisfies the Leibniz rule in
	each argument, and its Jacobiator, the cyclic sum whose vanishing is the
	Jacobi identity \cite[Section~2.1]{Balseiro2014}, is
	\begin{equation}
		[[a,b]_{\K},e]_{\K}+[[b,e]_{\K},a]_{\K}+[[e,a]_{\K},b]_{\K}
		=-u^4\,(\dd r\wedge\dd H\wedge\dd s)(\Z_a,\Z_b,\Z_e).
		\label{eq:jacobiator}
	\end{equation}
	Consequently, the following three statements are equivalent on $U$:
	\begin{enumerate}[label=(\alph*)]
		\item $[\cdot,\cdot]_{\K}$ is a Poisson bracket;
		\item $\dd s\wedge\dd r\wedge\dd H=0$;
		\item the ratio $r$ obeys a closed first-order equation on each energy
		level, that is, about each point of $U$ there is a smooth function $G$ of
		two variables such that $\dot r=G(r,H)$.
	\end{enumerate}
	When they hold, the symplectic leaves of $[\cdot,\cdot]_{\K}$ are the
	two-dimensional integral surfaces of $\spanop\{\Z_r,\Z_H\}$, and its Casimir
	functions are the common first integrals of $r$ and $H$, namely the
	functions $f$ with $\{f,r\}=\{f,H\}=0$.
\end{corollary}

\begin{proof}
	By \eqref{eq:representations} and \eqref{eq:ratio},
	$[a,b]_{\K}=\K(\dd a,\dd b)$ with $\K=g\,\mathbf X\wedge\mathbf V$, where
	$g=u^2$, $\mathbf X=\Z_r$ and $\mathbf V=\Z_H$. A bracket defined by a
	bivector is always skew and satisfies the Leibniz rule. For a bivector of
	this form, a direct computation gives
	\[
	[[a,b]_{\K},e]_{\K}+[[b,e]_{\K},a]_{\K}+[[e,a]_{\K},b]_{\K}
	=-g^2\,\bigl(\mathbf X\wedge\mathbf V\wedge[\mathbf X,\mathbf V]\bigr)(\dd a,\dd b,\dd e),
	\]
	in which the terms containing derivatives of $g$ cancel because
	$\K(\dd g,\cdot)$ lies in $\spanop\{\mathbf X,\mathbf V\}$. With the
	conventions \eqref{eq:convention}, $[\Z_r,\Z_H]=\Z_{\{H,r\}}=-\Z_s$, and
	$\dd w(\Z_f)=-\dd f(\Z_w)$ for all functions $f,w$; together these turn the
	right-hand side into \eqref{eq:jacobiator}.
	
	\emph{(a)$\Leftrightarrow$(b).} The Hamiltonian fields $\Z_a$ span every
	tangent space and $u\ne0$, so by \eqref{eq:jacobiator} the Jacobiator
	vanishes identically if and only if the three-form
	$\dd r\wedge\dd H\wedge\dd s$ does.
	
	\emph{(b)$\Leftrightarrow$(c).} Since $\dd r\wedge\dd H\ne0$, every point of
	$U$ has a neighbourhood with box-shaped local coordinates
	$(y_1,\ldots,y_{2n})$ in which $y_1=r$ and $y_2=H$. In these coordinates
	\[
	\dd s\wedge\dd r\wedge\dd H
	=\sum_{k=3}^{2n}\frac{\partial s}{\partial y_k}\,\dd y_k\wedge\dd y_1\wedge\dd y_2,
	\]
	which vanishes if and only if $s$ depends on $(y_1,y_2)$ alone, that is,
	$s=G(r,H)$. Finally, $s=\{r,H\}=\dot r$ by \eqref{eq:convention}.
	
	\emph{Leaves and Casimirs.} The map $\dd f\mapsto\K(\dd f,\cdot)$ sends
	$\dd f$ to $u^2\bigl(\{f,r\}\,\Z_H-\{f,H\}\,\Z_r\bigr)$. Its image is
	$\spanop\{\Z_r,\Z_H\}$, whose integral surfaces exist by (b) and are the
	symplectic leaves, and its kernel consists of the differentials of the
	functions with $\{f,r\}=\{f,H\}=0$.
\end{proof}

Corollary \ref{cor:jacobi} is the rank-two case of a general fact. A bivector
$\boldsymbol\pi$ defines a Poisson bracket exactly when its Schouten bracket
with itself vanishes, $[\boldsymbol\pi,\boldsymbol\pi]=0$ \cite{Lichnerowicz1977},
\cite[Chapter~1]{Vaisman}. The Schouten bracket \cite{Schouten1940,Nijenhuis1955}
extends the Lie bracket of vector fields to multivector
fields, and for a bivector, $[\boldsymbol\pi,\boldsymbol\pi]$ is, up to a
convention-dependent constant, the three-vector whose values on
$(\dd a,\dd b,\dd e)$ form the Jacobiator of the bracket. For a bivector of
rank two, this condition reduces to the Frobenius integrability of its image
distribution, with any nonvanishing conformal factor allowed
\cite[Section~2]{Ramos}. The five-input
operator identifies that distribution as $\spanop\{\Z_r,\Z_H\}$ and gives its
integrability a dynamical reading: the induced bracket is Poisson precisely
when the response ratio $r=u^{-1}v$ evolves
autonomously on each energy level. The condition holds automatically when $r$ is a first integral of $H$, since
then $s=0$ and the flows of $r$ and $H$ commute, so that the leaves are
their joint orbits. For $n=1$ it is empty, because every three-form on a
surface vanishes, in agreement with the fact that every bivector on a
surface is Poisson. The case $r=F(H)$, on the other hand, is excluded by the
hypothesis $\dd r\wedge\dd H\ne0$, since by Proposition~\ref{prop:covariance}
the operator then vanishes for all drivers and the bivector is identically zero.

Summing the operator over canonical driver pairs produces a chart-independent
quantity with a transparent geometric meaning: it is the signed curvature,
in the classical sense of plane curves \cite{doCarmo}, of the curve traced
by the responses $(c,d)$ along a trajectory.

\begin{proposition}[Projected turning]
	\label{prop:turning}
	Let $H,c,d\in C^\infty(M)$ be given functions, and define the canonical
	trace
	\[
	S_{H;c,d}=\sum_{i=1}^n\B_{q_i,p_i;c,d}[H]
	\]
	in a canonical chart $(q_1,\ldots,q_n,p_1,\ldots,p_n)$. Then the following
	hold.
	\begin{enumerate}[label=(\alph*)]
		\item \emph{Trace formula.}
		\begin{equation}
			S_{H;c,d}=u\,\Z_H(v)-v\,\Z_H(u).
			\label{eq:trace}
		\end{equation}
		\item \emph{Chart independence.} $S_{H;c,d}$ does not depend on the
		canonical chart used to define it.
		\item \emph{Ratio form.} Where $u\ne0$, $S_{H;c,d}=u^2\{r,H\}$.
		\item \emph{Curvature.} Along a trajectory $z(t)$ of $\Z_H$, the curve
		$\mathbf C(t)=\bigl(c(z(t)),d(z(t))\bigr)$ in the Euclidean $(c,d)$-plane
		has signed curvature
		\[
		\kappa_{\mathbf C}=\frac{S_{H;c,d}}{(u^2+v^2)^{3/2}}
		\]
		wherever $u^2+v^2\ne0$.
	\end{enumerate}
\end{proposition}

\begin{proof}
	\emph{(a)} In the chart, the sector derivations are
	\[
	\Y^H_{q_i,p_i}=H_{p_i}\partial_{q_i}-H_{q_i}\partial_{p_i},
	\qquad \sum_{i=1}^n\Y^H_{q_i,p_i}=\Z_H,
	\]
	and \eqref{eq:B} is linear in the derivation, so summing over $i$ gives
	\eqref{eq:trace}.
	
	\emph{(b)} The right-hand side of \eqref{eq:trace} involves only $H$, $c$,
	$d$ and the Hamiltonian vector field $\Z_H$, none of which refers to a
	chart.
	
	\emph{(c)} Substituting $v=ur$ into \eqref{eq:trace} gives
	$u\,\Z_H(ur)-ur\,\Z_H(u)=u^2\,\Z_H(r)$, and $\Z_H(r)=\{r,H\}$ by
	\eqref{eq:convention}.
	
	\emph{(d)} Along the flow, $\dot G=\Z_H(G)$ for every function $G$, so
	\[
	\dot{\mathbf C}=-(u,v),\qquad \ddot{\mathbf C}=-\bigl(\Z_H(u),\Z_H(v)\bigr).
	\]
	The signed curvature of a plane curve is the determinant formed by its
	velocity and acceleration, divided by the cube of its speed
	\cite{doCarmo}. Here the determinant is $u\,\Z_H(v)-v\,\Z_H(u)=S_{H;c,d}$
	by (a), and the speed is $(u^2+v^2)^{1/2}$.
\end{proof}

\begin{remark}[Isoenergetic nondegeneracy]
\label{rem:arnold}
Let $H_0=H_0(I_1,I_2)$ be an integrable Hamiltonian written in action--angle
variables $(\varphi_1,\varphi_2,I_1,I_2)$, with frequencies
$\omega_i=\partial_{I_i}H_0$, which form the frequency vector
$\partial H_0/\partial I=(\omega_1,\omega_2)^{\mathsf T}$. With both angles taken as drivers and
as responses, $(a,b;c,d)=(\varphi_1,\varphi_2;\varphi_1,\varphi_2)$, the
derivations in \eqref{eq:bordered} are $D_{\varphi_i}=\partial_{I_i}$, which
gives the first equality in
\begin{equation}
 \B_{\varphi_1,\varphi_2;\varphi_1,\varphi_2}[H_0]
 =\det\begin{pmatrix}
  0&\omega_1&\omega_2\\[2pt]
  \omega_1&\dfrac{\partial^2H_0}{\partial I_1^2}
  &\dfrac{\partial^2H_0}{\partial I_1\partial I_2}\\[10pt]
  \omega_2&\dfrac{\partial^2H_0}{\partial I_1\partial I_2}
  &\dfrac{\partial^2H_0}{\partial I_2^2}
 \end{pmatrix}
 =\det\begin{pmatrix}
  \dfrac{\partial^2H_0}{\partial I^2}&\dfrac{\partial H_0}{\partial I}\\[10pt]
  \Bigl(\dfrac{\partial H_0}{\partial I}\Bigr)^{\mathsf T}&0
 \end{pmatrix},
 \label{eq:arnold}
\end{equation}
and moving the first row and the first column to the last position gives
the second. The last determinant is the one introduced by Arnold, and its
nonvanishing is the isoenergetic nondegeneracy condition, under which the
invariant tori of the KAM theorem persist within each energy level
\cite{ArnoldSmallDen1963}, \cite[Appendix~8]{Arnold}. The identity
\eqref{eq:arnold} follows directly from \eqref{eq:bordered}, and the ratio
and angle forms of Proposition~\ref{prop:structure} give the condition its
classical geometric reading. Since $u=-\omega_1$ and $v=-\omega_2$, the
sector derivation
\[
 \Y^{H_0}_{\varphi_1,\varphi_2}=\omega_2\,\partial_{I_1}-\omega_1\,\partial_{I_2}
\]
is tangent to the energy curves $\{H_0=E\}$ in the action plane, and the
ratio form \eqref{eq:ratio} reads
\[
 \B_{\varphi_1,\varphi_2;\varphi_1,\varphi_2}[H_0]
 =\omega_1^2\,\Y^{H_0}_{\varphi_1,\varphi_2}\Bigl(\frac{\omega_2}{\omega_1}\Bigr)
\]
wherever $\omega_1\ne0$. The condition therefore holds at an action $I$
exactly when the frequency ratio has nonzero derivative at $I$ along the
energy curve. By the angle form \eqref{eq:angle}, the same quantity is
$\omega_1^2+\omega_2^2$ times the rate at which the direction of the
frequency vector, normal to the energy curve, turns along the flow of $\Y^{H_0}_{\varphi_1,\varphi_2}$. This
flow traverses the energy curve with speed $(\omega_1^2+\omega_2^2)^{1/2}$,
so the energy curve,
oriented by $\Y^{H_0}_{\varphi_1,\varphi_2}$, has signed curvature
$(\omega_1^2+\omega_2^2)^{-3/2}\,\B_{\varphi_1,\varphi_2;\varphi_1,\varphi_2}[H_0]$
\cite{doCarmo}, which is the analogue in the action plane of
Proposition~\ref{prop:turning}(d).

Proposition~\ref{prop:algebra} adds two properties. By item (f), a
reparametrisation $\chi(H_0)$ with $\chi'\ne0$ multiplies \eqref{eq:arnold}
by $\chi'(H_0)^3$, in agreement with the geometric reading, since such a
reparametrisation leaves the energy curves unchanged. By item (g), the matrix
of double brackets is the Hessian $\mathbf h=\partial^2H_0/\partial I^2$,
whose determinant is Kolmogorov's nondegeneracy condition
\cite{Kolmogorov1954}, \cite[Appendix~8]{Arnold}, and the same
reparametrisation gives
\[
 \det\mathbf h\bigl[\chi(H_0)\bigr]
 =\chi'(H_0)^2\det\mathbf h[H_0]
 -\chi'(H_0)\,\chi''(H_0)\,\B_{\varphi_1,\varphi_2;\varphi_1,\varphi_2}[H_0].
\]
Kolmogorov's condition thus depends on the parametrisation of the energy,
and the two conditions are independent. The Hamiltonian
$H_0=I_1+\frac12I_2^2$ has $\det\mathbf h=0$ and \eqref{eq:arnold} equal to
$-1$, and the Hamiltonian $H_0=\ln I_1-\ln I_2$ on $I_1,I_2>0$ has
$\det\mathbf h=-(I_1I_2)^{-2}$ and \eqref{eq:arnold} identically zero,
because its energy curves are rays through the origin.

For $n\ge3$ degrees of freedom, Arnold's determinant is the
$(n+1)\times(n+1)$ bordered Hessian built in the same way. Sylvester's
determinant identity \cite[Chapter~II, \S3, eq.~(28)]{Gantmacher}, applied to the first
matrix form with the border and $\varphi_1$ as pivot block, expresses it
through five-input operators,
\[
 \det\begin{pmatrix}
  \dfrac{\partial^2H_0}{\partial I^2}&\dfrac{\partial H_0}{\partial I}\\[10pt]
  \Bigl(\dfrac{\partial H_0}{\partial I}\Bigr)^{\mathsf T}&0
 \end{pmatrix}
 =\bigl(-\omega_1^2\bigr)^{2-n}
 \det\bigl(\B_{\varphi_1,\varphi_i;\varphi_1,\varphi_j}[H_0]\bigr)_{i,j=2}^{n},
\]
wherever $\omega_1\ne0$. The probes in \eqref{eq:arnold} are the angles
alone. With the canonical pairs $(\varphi_1,I_1)$ and $(\varphi_2,I_2)$ as
probes, the last column of \eqref{eq:bordered} vanishes, because $H_0$ does
not depend on the angles, so that
$\B_{\varphi_1,I_1;\varphi_2,I_2}[H_0]=0$, in agreement with the
separability of every function of the actions (Section~\ref{sec:product}).
\end{remark}

\section{Product structures and separation}
\label{sec:product}

Additive separation of the Hamilton--Jacobi equation is a property of a
Hamiltonian together with a splitting of phase space into degrees of
freedom, and the classical theory is accordingly formulated in a fixed
system of canonical coordinates \cite{LeviCivita,Eisenhart,Benenti2016}. In
this section, such a splitting is fixed, in the form of a local product of
symplectic surfaces, and the five-input operator is evaluated on the
canonical pairs of the factors. These evaluations reproduce the Levi-Civita
expressions, and their off-diagonal part turns out to be intrinsic to the
splitting, since it can be written both as a mixed Hessian and as a
symplectic pairing of the sector Hamiltonian fields
(Proposition~\ref{prop:intrinsic}). Its vanishing is equivalent to the
Frobenius integrability of the Lagrangian distribution spanned by the sector
fields (Proposition~\ref{prop:frobenius}). This extends the theorem of Llibre
and Ram\'{\i}rez \cite[Theorem~10 and Remark~11]{LlibreRamirez2020}, who
established the equivalence where every $H_{p_i}$, or every $H_{q_i}$, is
nonzero, and leads to a local separation theorem in which the classical
momentum regularity is replaced by the nonvanishing of the sector fields
(Theorem~\ref{thm:sector}). The same
condition is finally expressed through a Nijenhuis operator adapted to the
splitting (Proposition~\ref{prop:nijenhuis}), in line with the
bi-Hamiltonian approach of Falqui and Pedroni \cite{FalquiPedroni} and with
Nijenhuis geometry \cite{BolsinovKonyaevMatveev}.

\subsection{Sector fields and the intrinsic off-diagonal values}

The starting point is the observation that, on canonical pairs, the
five-input operator reproduces the classical Levi-Civita expressions.
Levi-Civita derived them in a letter in Italian to St\"ackel, dated Padua,
3~February 1904, and published the same year in \emph{Mathematische Annalen}
\cite[eqs.~(3)]{LeviCivita},\footnote{In the first term of his equations (3), the second derivative is
	printed with $p_i$ in place of $x_i$, as $H_{p_ix_j}$ instead of $H_{x_ix_j}$.
	With $x_i=q_i$, the corrected condition is $L_{ij}=0$ in \eqref{eq:LC}. The correct form appears on p.~386, where he restates (3) with the kinetic
	energy $K$ in place of $H$ as his group (I).} in which he showed that
the Hamilton--Jacobi equation separates in the given coordinates exactly when
$H$ satisfies one such second-order equation for each pair of degrees of
freedom.

\begin{proposition}[Levi-Civita expressions]
	\label{prop:canonical}
	Let $(q_1,\ldots,q_n,p_1,\ldots,p_n)$ be a canonical chart. For all
	$i,j\in\{1,\ldots,n\}$,
	\begin{equation}
		\B_{q_i,p_i;q_j,p_j}[H]
		=-H_{p_i}H_{p_j}H_{q_iq_j}+H_{q_i}H_{p_j}H_{p_iq_j}
		+H_{p_i}H_{q_j}H_{q_ip_j}-H_{q_i}H_{q_j}H_{p_ip_j}
		=:-L_{ij}(H),
		\label{eq:LC}
	\end{equation}
	where, for $i\ne j$, $L_{ij}$ is the classical Levi-Civita expression
	\cite{LeviCivita,Benenti2016}.
\end{proposition}

\begin{proof}
	With $c=q_j$ and $d=p_j$, the response rates are $u=-H_{p_j}$ and
	$v=H_{q_j}$, and the derivation is
	$\Y^H_{q_i,p_i}=H_{p_i}\partial_{q_i}-H_{q_i}\partial_{p_i}$. Substituting
	into \eqref{eq:B} gives
	\[
	\B_{q_i,p_i;q_j,p_j}[H]
	=-H_{p_j}\,\Y^H_{q_i,p_i}(H_{q_j})+H_{q_j}\,\Y^H_{q_i,p_i}(H_{p_j}),
	\]
	which expands to \eqref{eq:LC}.
\end{proof}

The overall sign in \eqref{eq:LC} is fixed by Definition~\ref{def:B}. The
vanishing of the expressions $L_{ij}$ for $i\ne j$ is Levi-Civita's
separability condition; its equivalence with separation is recalled in
Proposition~\ref{prop:LC} and extended in Theorem~\ref{thm:sector}.

The expression obtained in the proof before the expansion,
$-H_{p_j}\,\Y^H_{q_i,p_i}(H_{q_j})+H_{q_j}\,\Y^H_{q_i,p_i}(H_{p_j})$, is the
form in which Llibre and Ram\'{\i}rez write the Levi-Civita conditions
\cite[eq.~(13)]{LlibreRamirez2020}. Where $H_{p_j}\ne0$, the ratio form
\eqref{eq:ratio} with $u=-H_{p_j}$ and $v=H_{q_j}$ gives
\[
 \B_{q_i,p_i;q_j,p_j}[H]=H_{p_j}^2\,\Y^H_{q_i,p_i}(R_j),\qquad
 R_j=-\frac{H_{q_j}}{H_{p_j}},
\]
which is their condition $\Y^H_{q_i,p_i}(R_j)=0$
\cite[eqs.~(14) and (15)]{LlibreRamirez2020}, and
Proposition~\ref{prop:algebra}(f), restricted to canonical pairs, is the
identity proved in \cite[Proposition~5]{LlibreRamirez2020}. The five-input
operator extends these canonical-pair forms to arbitrary drivers and
responses.

Suppose now that $M$ is locally a product of symplectic surfaces $M_i$, with
$\om=\sum_i\om_i$. Let $E_i\subset TM$ denote the tangent plane of the $i$th
factor, $\Pib_i:TM\to E_i$ the corresponding projection, and $(q_i,p_i)$ any
canonical chart on $M_i$. The product splitting is part of the data, and
with respect to it the \emph{sector fields} and the \emph{off-diagonal
	values} are defined by
\begin{equation}
	\Y_i=\Pib_i\Z_H=H_{p_i}\partial_{q_i}-H_{q_i}\partial_{p_i},\qquad
	\Z_H=\sum_{i=1}^n\Y_i,\qquad
	b_{ij}=\B_{q_i,p_i;q_j,p_j}[H]\quad(i\ne j).
	\label{eq:sectorfields}
\end{equation}
For $i\ne j$, the product structure also defines the mixed second
differential $\mathbf H_{ij}\in\Gamma(E_i^*\otimes E_j^*)$ of $H$: vectors
in $E_i$ and $E_j$ are extended to fields $\mathbf X_i,\mathbf V_j$ whose
coefficients depend only on their own factors, and
$\mathbf H_{ij}(\mathbf X_i,\mathbf V_j)=\mathbf X_i(\mathbf V_jH)$. In
product coordinates this is the mixed block of the Hessian, and, unlike the
full Hessian, it transforms as a tensor. Under a general change of
coordinates, the Hessian of a function acquires terms in which the second
derivatives of the change multiply the first derivatives of the function, so
that it is a tensor only at critical points. Independent changes of chart
within the factors have vanishing second derivatives with respect to
variables of two different factors, and no such terms appear in the mixed
block.

\begin{proposition}[Intrinsic off-diagonal values]
	\label{prop:intrinsic}
	For $i\ne j$, the off-diagonal values have the following properties.
	\begin{enumerate}[label=(\alph*)]
		\item \emph{Intrinsic representations.}
		\begin{equation}
			b_{ij}=-\mathbf H_{ij}(\Y_i,\Y_j)=\om\bigl(\Y_j,[\Y_i,\Y_j]\bigr).
			\label{eq:intrinsic}
		\end{equation}
		\item \emph{Symmetry.} $b_{ij}=b_{ji}$.
		\item \emph{Invariance.} $b_{ij}$ is unchanged by independent canonical
		changes of chart within the factors, and a permutation of the factors
		permutes its indices accordingly.
	\end{enumerate}
\end{proposition}

\begin{proof}
	\emph{(a)} In product coordinates, $\mathbf H_{ij}(\Y_i,\Y_j)$ is the mixed
	block of the Hessian evaluated on the components of $\Y_i$ and $\Y_j$, and
	this reproduces $L_{ij}$ in \eqref{eq:LC}, which gives the first equality.
	For the second, the $E_j$ component of $[\Y_i,\Y_j]$ is
	\begin{equation}
		\Pib_j[\Y_i,\Y_j]=\Y_i(H_{p_j})\,\partial_{q_j}-\Y_i(H_{q_j})\,\partial_{p_j},
		\label{eq:sectorbracket}
	\end{equation}
	and its pairing with $\Y_j$ is
	\[
	\om\bigl(\Y_j,\Pib_j[\Y_i,\Y_j]\bigr)
	=-H_{p_j}\,\Y_i(H_{q_j})+H_{q_j}\,\Y_i(H_{p_j})=b_{ij};
	\]
	the $E_i$ component of the bracket does not contribute, because $E_i$ and
	$E_j$ are $\om$-orthogonal.
	
	\emph{(b)} Mixed partial derivatives commute, so
	$\mathbf H_{ij}(\mathbf X,\mathbf V)=\mathbf H_{ji}(\mathbf V,\mathbf X)$,
	and the first equality in \eqref{eq:intrinsic} gives $b_{ij}=b_{ji}$.
	
	\emph{(c)} Both representations in \eqref{eq:intrinsic} involve only $H$,
	$\om$ and the splitting. A canonical change of chart within the factors
	preserves all three, and in particular the sector fields $\Y_i=\Pib_i\Z_H$,
	so it leaves $b_{ij}$ unchanged. Relabelling the factors relabels the sector
	fields, and hence the indices of $b_{ij}$.
\end{proof}

\begin{figure}[t]
	\centering
	\begin{tikzpicture}[scale=1.15,>=Stealth]
		\draw[->,gray] (-0.3,0) -- (4.6,0) node[right,black] {$q_j$};
		\draw[->,gray] (0,-0.3) -- (0,3.0) node[above,black] {$p_j$};
		\draw[thick,gray!70] (0.400,0.450) .. controls (1.500,0.750) and (2.600,1.450) .. (4.200,2.800);
		\fill[blue!8] (1.823,1.058) -- (3.342,1.928) -- (3.124,3.013) -- (1.605,2.143) -- cycle;
		\draw[->,very thick,blue!70!black] (1.823,1.058) -- (3.342,1.928) node[right] {$\Y_j$};
		\draw[->,very thick,orange!80!black] (1.823,1.058) -- (1.605,2.143) node[above] {$\Pib_j[\Y_i,\Y_j]$};
		\node[blue!50!black] at (2.474,2.035) {$b_{ij}$};
		\fill (1.823,1.058) circle (1.4pt);
		\node[gray!80!black] at (1.45,0.28) {\small projected orbit};
		\begin{scope}[xshift=6.4cm]
			\draw[->,gray] (-0.3,0) -- (4.6,0) node[right,black] {$q_j$};
			\draw[->,gray] (0,-0.3) -- (0,3.0) node[above,black] {$p_j$};
			\draw[thick,gray!70] (0.400,0.450) .. controls (1.500,0.750) and (2.600,1.450) .. (4.200,2.800);
			\draw[->,very thick,orange!80!black] (1.823,1.058) -- (3.819,2.201) node[right] {$\Pib_j[\Y_i,\Y_j]$};
			\draw[->,very thick,blue!70!black] (1.823,1.058) -- (2.908,1.679);
			\node[blue!70!black,below right] at (2.908,1.679) {$\Y_j$};
			\fill (1.823,1.058) circle (1.4pt);
			\node at (3.4,0.55) {$b_{ij}=0$};
		\end{scope}
	\end{tikzpicture}
	\caption{Geometry of the off-diagonal values in the $(q_j,p_j)$-plane of a
		factor chart. In both panels $\Y_j$ is tangent to the projected orbit, since
		it is the velocity of the $j$th projected motion. Left:
		$b_{ij}=\om(\Y_j,\Pib_j[\Y_i,\Y_j])$ is the oriented area spanned by $\Y_j$
		and its variation along sector $i$. Right: when $b_{ij}=0$ this variation is
		parallel to $\Y_j$, hence also tangent to the orbit, so motion in sector $i$
		changes the speed of the $j$th projected motion but not its direction.}
	\label{fig:turning}
\end{figure}
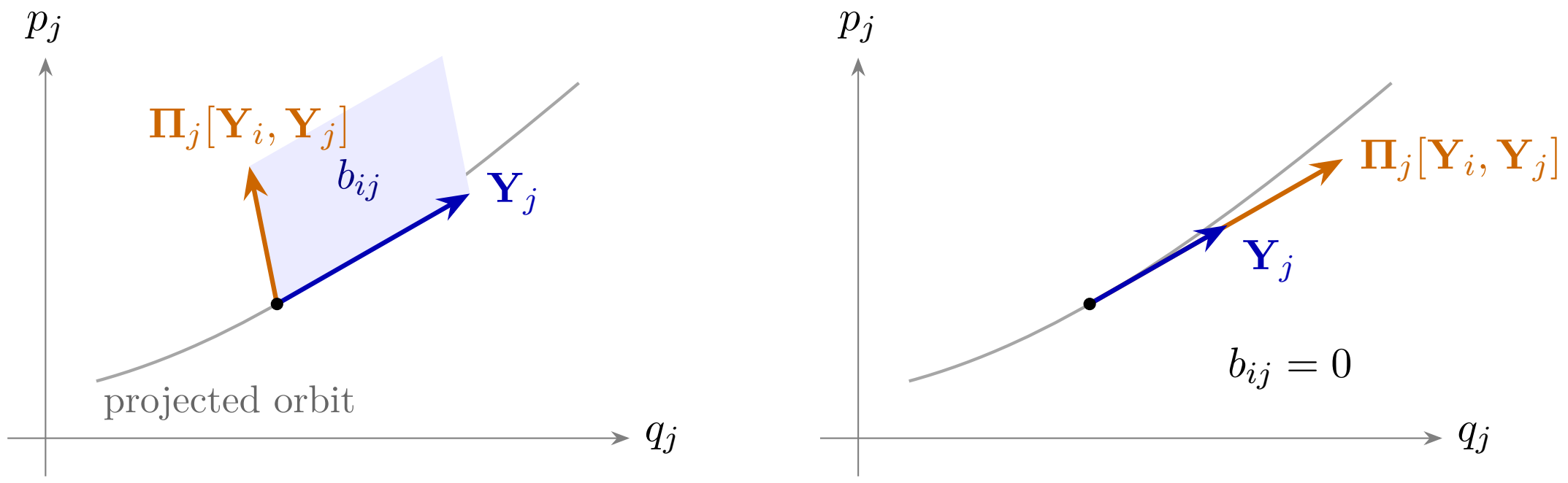

The diagonal entries $b_{jj}$, by contrast, do not vanish in general, because
the operator is alternating within the drivers and within the responses but
not under the exchange of the two pairs. They play no role in separation, and
they are not invariant: for $H=\frac12(q^2+p^2)$, $b_{11}=-2H$ in the chart
$(q,p)$ and $b_{11}=0$ in action--angle variables, and the example after
\eqref{eq:lastslot} is another diagonal entry. The off-diagonal values also inherit the turning interpretation
of Section~\ref{sec:operator}. Let $\vartheta_j=\atantwo(-H_{q_j},H_{p_j})$
be the direction angle of $\Y_j$ in a chart of $M_j$; it differs from the
angle of $(u,v)$ in \eqref{eq:angle} by the constant $\pi$, so
\eqref{eq:angle} gives
\begin{equation}
	b_{ij}=|\Y_j|^2\,\Y_i(\vartheta_j),\qquad
	S_{H;q_j,p_j}=\sum_{i=1}^nb_{ij}.
	\label{eq:turning-sector}
\end{equation}
Here $\Y_j$ is the velocity of the $j$th projected phase portrait, so
$b_{ij}$ measures how fast motion in sector $i$ turns the direction of the
phase velocity in sector $j$, and the $b_{ij}$ are the cross-sector
contributions to the curvature of that portrait. Neither the factor
$|\Y_j|^2$ nor the angle rate is chart independent on its own, but their
product is. The invariant form of this statement, drawn in Figure~\ref{fig:turning}, is
that when $\Y_j\ne0$, $b_{ij}=0$ if and only if $\Pib_j[\Y_i,\Y_j]$ is
parallel to $\Y_j$, in which case cross-sector variation rescales the $j$th
phase velocity without turning it.

\subsection{The Frobenius obstruction}

For a separable Hamiltonian, the leaves of a separated complete integral are
locally products of curves, one in each factor $M_i$, and on each leaf the
sector field $\Y_i$ is tangent to the corresponding curve. The sector fields
of a separable system are therefore tangent to a foliation, and a first test of separability is the integrability of the distribution they
span.
By the Frobenius theorem \cite[Chapter~19]{Lee}, integrability is equivalent to involutivity, so the
obstruction lies in the components of the brackets $[\Y_i,\Y_j]$ that leave
the distribution. For a product structure these components are controlled
by the off-diagonal values alone. The next proposition shows that, where
the sector fields do not vanish, the vanishing of all off-diagonal values
$b_{ij}=\B_{q_i,p_i;q_j,p_j}[H]$, $i\ne j$, defined in
\eqref{eq:sectorfields}, is exactly
the integrability condition of the distribution they span. Where every
$H_{p_i}$, or every $H_{q_i}$, is nonzero, the equivalence between the
Levi-Civita conditions and the involutivity of this distribution was
established by Llibre and Ram\'{\i}rez \cite[Theorem~10 and
Remark~11]{LlibreRamirez2020}, and the same condition is the geometric
content behind the invariant Lagrangian distributions of the bi-Hamiltonian
approach \cite{FalquiPedroni}. In comparison with both, the obstruction is expressed here through a single
invariant scalar for each pair of degrees of freedom.

\begin{proposition}[Frobenius obstruction]
	\label{prop:frobenius}
	Let $U$ be an open set on which every sector field $\Y_i$ is nonzero, and
	put $\mathcal P_H=\spanop\{\Y_1,\ldots,\Y_n\}$. Then the following hold on
	$U$.
	\begin{enumerate}[label=(\alph*)]
		\item \emph{Lagrangian distribution.} $\mathcal P_H$ is a Lagrangian
		distribution of rank $n$.
		\item \emph{Tangency to the energy levels.} $\mathcal P_H$ is tangent to
		the level sets of $H$.
		\item \emph{Frobenius obstruction.} For vector fields $\mathbf X$ and
$\mathbf V$ tangent to $\mathcal P_H$, let $\mathbf F(\mathbf X,\mathbf V)$
be the class of $[\mathbf X,\mathbf V]$ in the quotient bundle
$TM/\mathcal P_H$, that is, the part of the bracket that leaves the
distribution. This Frobenius tensor of $\mathcal P_H$ satisfies
\begin{equation}
	\mathbf F=0\quad\Longleftrightarrow\quad b_{ij}=0\ \text{for all }i\ne j.
	\label{eq:frobenius}
\end{equation}
	\end{enumerate}
\end{proposition}

\begin{proof}
	\emph{(a)} Each $\Y_i$ is a nonzero vector in $E_i$, and the planes $E_i$
	form a direct sum, so $\Y_1,\ldots,\Y_n$ are linearly independent and
	$\mathcal P_H$ has rank $n$. For $i\ne j$, $\om(\Y_i,\Y_j)=0$ because $E_i$
	and $E_j$ are $\om$-orthogonal, and $\om(\Y_i,\Y_i)=0$ by skew-symmetry.
	Hence $\mathcal P_H$ is isotropic of half the dimension of $M$, that is,
	Lagrangian in the usual sense \cite{Weinstein1971,McDuffSalamon,Arnold}.
	
	\emph{(b)} By \eqref{eq:sectorfields},
	$\dd H(\Y_i)=H_{p_i}H_{q_i}-H_{q_i}H_{p_i}=0$ for every $i$.
	
	\emph{(c)} Fix $i\ne j$. The coefficients of $\Y_j$ multiply only
	$\partial_{q_j},\partial_{p_j}$ and those of $\Y_i$ only
	$\partial_{q_i},\partial_{p_i}$, so $[\Y_i,\Y_j]$ lies in $E_i\oplus E_j$.
	In the symplectic plane $E_j$, the $\om$-orthogonal complement of the
	nonzero vector $\Y_j$ is its own span, so by
	Proposition~\ref{prop:intrinsic}(a)
	\[
	\Pib_j[\Y_i,\Y_j]\in\spanop\{\Y_j\}
	\quad\Longleftrightarrow\quad
	\om\bigl(\Y_j,[\Y_i,\Y_j]\bigr)=b_{ij}=0 .
	\]
	The same argument in $E_i$, applied to $[\Y_j,\Y_i]=-[\Y_i,\Y_j]$, shows
	that $\Pib_i[\Y_i,\Y_j]\in\spanop\{\Y_i\}$ if and only if $b_{ji}=0$, which
	by Proposition~\ref{prop:intrinsic}(b) is again $b_{ij}=0$. Together,
	\begin{equation*}
		[\Y_i,\Y_j]\in\spanop\{\Y_i,\Y_j\}
		\quad\Longleftrightarrow\quad b_{ij}=0 .
	\end{equation*}
	($\Rightarrow$) If $\mathbf F=0$, then $[\Y_i,\Y_j]\in\mathcal P_H$. Since
	$\Y_k\in E_k$ for every $k$ and the planes form a direct sum,
	$\mathcal P_H\cap(E_i\oplus E_j)=\spanop\{\Y_i,\Y_j\}$, so
	$[\Y_i,\Y_j]\in\spanop\{\Y_i,\Y_j\}$ and therefore $b_{ij}=0$.
	
	\noindent ($\Leftarrow$) If all $b_{ij}$ vanish, then
	$[\Y_i,\Y_j]\in\spanop\{\Y_i,\Y_j\}\subset\mathcal P_H$ for all $i,j$, so
	$\mathbf F$ vanishes on the local frame $\Y_1,\ldots,\Y_n$. Since
	\[
	[f\mathbf X,g\mathbf V]=fg\,[\mathbf X,\mathbf V]
	+f\,\mathbf X(g)\,\mathbf V-g\,\mathbf V(f)\,\mathbf X
	\]
	and the last two terms lie in $\mathcal P_H$ for
	$\mathbf X,\mathbf V\in\mathcal P_H$, the tensor $\mathbf F$ is
	$C^\infty$-bilinear, so that its vanishing on a local frame implies that it
	vanishes identically, which is the local-frame criterion for involutivity
	\cite[Chapter~19]{Lee}.
\end{proof}

\subsection{Separation under sector regularity}

Proposition~\ref{prop:frobenius} describes separability through the
distribution $\mathcal P_H$, whereas separation of variables in the classical
sense refers to a complete integral of the Hamilton--Jacobi equation
\cite{Stackel,LeviCivita,Eisenhart}. The link between the two descriptions is
provided by the leaves of $\mathcal P_H$, since a Lagrangian submanifold that
is a graph over configuration space is locally the graph of the differential
of a generating function \cite{Arnold,Weinstein1971}, and this generating
function is additively separated exactly when the leaf splits into curves in
the individual factors. The classical criterion of Levi-Civita assumes that every
$H_{p_i}$ is nonzero, which guarantees that the leaves are momentum graphs
$\mathbf p=\mathbf P(\mathbf q)$. This criterion is recalled first, with a
proof adapted to the sign conventions of Section~\ref{sec:operator}, and the
assumption is then replaced by the weaker requirement that the sector fields
do not vanish. The condition $H_{p_i}\ne0$ depends on the canonical chart
chosen in the $i$th factor, while $\Y_i\ne0$ is intrinsic to the product
structure. The resulting theorem therefore links the off-diagonal values to
the invariant-Lagrangian approach to separation of \cite{FalquiPedroni}.

\begin{proposition}[Levi-Civita criterion for momentum graphs]
\label{prop:LC}
Suppose that $H_{p_i}\ne0$ for all $i$ on an open canonical domain. Then
$b_{ij}=0$ there for all $i\ne j$ if and only if, locally about each point,
the equation $H(\mathbf q,\partial_{\mathbf q}\mathcal S)=\kappa_1$ has a
complete additively separated integral
$\mathcal S(\mathbf q;\boldsymbol\kappa)=\sum_i\mathcal S_i(q_i;\boldsymbol\kappa)$,
with $\det(\partial^2\mathcal S/\partial q_i\partial\kappa_j)\ne0$.
\end{proposition}

\begin{proof}
Put $\gamma_i=H_{q_i}/H_{p_i}$ and
$\mathbf D_i=\partial_{q_i}-\gamma_i\partial_{p_i}=\Y_i/H_{p_i}$. By
\eqref{eq:ratio}, $\mathbf D_i(\gamma_j)=-b_{ij}/(H_{p_i}H_{p_j}^2)$, and the
symmetry of $b_{ij}$ then gives
\begin{equation}
 [\mathbf D_i,\mathbf D_j]
 =\frac{b_{ij}}{H_{p_i}^2H_{p_j}^2}\bigl(H_{p_i}\partial_{p_j}-H_{p_j}\partial_{p_i}\bigr).
 \label{eq:commutator}
\end{equation}
If all $b_{ij}$ vanish, the independent fields $\mathbf D_i$ commute, and
their integral leaves are momentum graphs
$\mathbf p=\mathbf P(\mathbf q;\boldsymbol\zeta)$ normalised by
$\mathbf P(\mathbf q^0;\boldsymbol\zeta)=\boldsymbol\zeta$. Tangency gives
$\partial P_i/\partial q_j=0$ for $i\ne j$, and since $\mathbf D_iH=0$, the
Hamiltonian is constant on the leaves. Because
$\partial H/\partial\zeta_k=H_{p_k}\ne0$ at $\mathbf q^0$, one leaf parameter
can be replaced by the energy $\kappa_1$; integrating each $P_i$ in $q_i$
then gives $\mathcal S$, and completeness follows because
$\partial\mathbf P/\partial\boldsymbol\zeta$ is the identity at
$\mathbf q^0$. Conversely, a complete separated integral defines a local
foliation by the graphs $P_i=\partial\mathcal S_i/\partial q_i$.
Differentiating the Hamilton--Jacobi equation shows that the tangent fields
of these graphs are the $\mathbf D_j$, whose commutators are tangent to the
graphs and have no configuration component; hence the commutators vanish,
and \eqref{eq:commutator} gives $b_{ij}=0$.
\end{proof}

The hypothesis $H_{p_i}\ne0$ serves only to write the leaves as momentum
graphs $\mathbf p=\mathbf P(\mathbf q)$ and imposes no condition on the
distribution $\mathcal P_H$ itself. Theorem~\ref{thm:sector} replaces it by the
weaker requirement that the sector fields do not vanish.

\begin{theorem}[Separation under sector regularity]
	\label{thm:sector}
	Suppose that $\Y_i(z_*)\ne0$ for all $i$. Then there is a neighbourhood of
	$z_*$ on which the following statements are equivalent.
	\begin{enumerate}[label=(\alph*)]
		\item \emph{Vanishing off-diagonal values.} $b_{ij}=0$ identically for all
		$i\ne j$.
		\item \emph{Involutivity.} The distribution $\mathcal P_H$ is involutive.
		\item \emph{Product foliation.} $\mathcal P_H$ is tangent to a foliation
		whose leaves are local products of curves in the factors $M_i$, each leaf
		lying in a level set of $H$.
		\item \emph{Separated complete integral.} In suitable canonical charts
		$(Q_i,P_i)$ on the factors, the Hamilton--Jacobi equation
		$H(\mathbf Q,\partial_{\mathbf Q}\mathcal S)=\kappa_1$ has a complete
		additively separated integral
		$\mathcal S(\mathbf Q;\boldsymbol\kappa)=\sum_i\mathcal S_i(Q_i;\boldsymbol\kappa)$
		in the sense of Proposition~\ref{prop:LC}.
	\end{enumerate}
\end{theorem}

\begin{proof}
		Throughout the proof, $U$ denotes an open neighbourhood of $z_*$ on which
	every $\Y_i$ is nonzero. The Frobenius theorem, the constant-rank theorem,
	the inverse function theorem and the continuity argument for the
	derivatives $H_{P_i}$ below each hold on some smaller open neighbourhood of
	$z_*$, and at each of these steps $U$ is replaced by that neighbourhood.
	Since there are finitely many steps, the final $U$ is still an open
	neighbourhood of $z_*$, which is all that the local statement of the
	theorem requires.
	
	\emph{(a)$\Leftrightarrow$(b).} Since every $\Y_i$ is nonzero on $U$, this
	is Proposition~\ref{prop:frobenius}(c).
	
\emph{(b)$\Rightarrow$(c).} By the Frobenius theorem \cite{Lee}, the
involutive rank-$n$ distribution $\mathcal P_H$ is tangent to a foliation of
$U$ by $n$-dimensional leaves. Let $L$ be such a leaf, and let
$\mathrm{pr}_i:M\to M_i$ be the projection onto the $i$th factor. Its
derivative maps $E_i$ isomorphically onto $TM_i$ and annihilates $E_j$ for
$j\ne i$, so that every tangent space of $L$ is mapped onto the line
\[
\mathrm{pr}_{i*}(T_zL)=\spanop\{\mathrm{pr}_{i*}\Y_i(z)\},\qquad z\in L.
\]
The restriction of $\mathrm{pr}_i$ to $L$ therefore has constant rank one,
and the constant-rank theorem \cite{Lee} shows that its image is locally an
embedded curve $\gamma_i\subset M_i$. Taking the $n$ projections together
gives the map
\[
\mathrm{pr}=(\mathrm{pr}_1,\ldots,\mathrm{pr}_n):\;
L\longrightarrow\gamma_1\times\cdots\times\gamma_n,
\]
whose derivative acts on tangent vectors of $L$ as
\[
\mathrm{pr}_*\Bigl(\sum_{j=1}^nc_j\Y_j\Bigr)
=\bigl(c_1\,\mathrm{pr}_{1*}\Y_1,\ldots,c_n\,\mathrm{pr}_{n*}\Y_n\bigr).
\]
Because the vectors $\mathrm{pr}_{i*}\Y_i$ are nonzero, this derivative is
injective, and since $L$ and the product of the curves both have dimension
$n$, the inverse function theorem shows that $L$ is locally an open subset of
$\gamma_1\times\cdots\times\gamma_n$. Finally, the leaf lies in an energy
level, because $\dd H(\Y_i)=0$ for every $i$ and $H$ is therefore constant on
each connected leaf.
	
	\emph{(c)$\Rightarrow$(b).} The Lie bracket of two vector fields tangent to
	the leaves of a foliation is again tangent to the leaves, so $\mathcal P_H$
	is involutive.
	
	\emph{(a)$\Rightarrow$(d).} In each factor, put
	\[
	(Q_i,P_i)=
	\begin{cases}
		(q_i,p_i), & H_{p_i}(z_*)\ne0,\\
		(p_i,-q_i), & H_{p_i}(z_*)=0.
	\end{cases}
	\]
	Both choices are canonical, because $\dd p_i\wedge\dd(-q_i)=\dd q_i\wedge\dd p_i$.
	In the second case $\Y_i=H_{p_i}\partial_{q_i}-H_{q_i}\partial_{p_i}$ is
	nonzero at $z_*$ only if $H_{q_i}(z_*)\ne0$, whence
	$H_{P_i}(z_*)=-H_{q_i}(z_*)\ne0$. In either case $H_{P_i}(z_*)\ne0$, and
	after shrinking $U$ all $H_{P_i}$ are nonzero on $U$. The off-diagonal
	values do not depend on the canonical chart chosen in each factor
	(Proposition~\ref{prop:intrinsic}(c)), so (a) holds in the charts
	$(Q_i,P_i)$, and Proposition~\ref{prop:LC} yields the separated complete
	integral of (d).
	
	\emph{(d)$\Rightarrow$(a).} Completeness means that the graphs
	$P_i=\partial\mathcal S_i/\partial Q_i$ foliate a neighbourhood of $z_*$, and
	on each graph $P_i$ depends on $Q_i$ alone. Differentiating
	$H(\mathbf Q,\partial_{\mathbf Q}\mathcal S)=\kappa_1$ with respect to $Q_i$
	therefore gives
	\[
	H_{Q_i}+H_{P_i}\,\frac{\partial^2\mathcal S_i}{\partial Q_i^2}=0 .
	\]
	At a point where $H_{P_i}$ vanished, $H_{Q_i}$ would vanish as well, and so
	would $\Y_i=H_{P_i}\partial_{Q_i}-H_{Q_i}\partial_{P_i}$, contrary to the
	choice of $U$. Hence all $H_{P_i}$ are nonzero, Proposition~\ref{prop:LC}
	applies in the charts $(Q_i,P_i)$ and gives $b_{ij}=0$ there, and
	Proposition~\ref{prop:intrinsic}(c) transfers this to every factorwise
	canonical chart, in particular to the original one.
\end{proof}

Theorem~\ref{thm:sector} goes beyond Proposition~\ref{prop:LC} at the turning
points of the separated motions, which every orbit that librates in a sector
reaches. A typical case is a Liouville Hamiltonian
\cite{Eisenhart,KalninsMiller}
\[
H=\frac{\tfrac12(p_1^2+p_2^2)+V_1(q_1)+V_2(q_2)}{W_1(q_1)+W_2(q_2)},
\qquad W_1+W_2>0,
\]
a class that contains the two-centre problem in its separating elliptic chart
(Section~\ref{sec:euler}). Its off-diagonal value vanishes identically, and
on the level $H=h$ the Hamilton--Jacobi equation splits into
\[
\tfrac12p_1^2+V_1(q_1)-hW_1(q_1)=\kappa_2,\qquad
\tfrac12p_2^2+V_2(q_2)-hW_2(q_2)=-\kappa_2 .
\]
At a turning point of the first motion, where $p_1=0$, the hypothesis of
Proposition~\ref{prop:LC} fails because $H_{p_1}=p_1/(W_1+W_2)$ vanishes, and
the leaf is indeed not a graph over $q_1$ there, since the two branches
$p_1=\pm\sqrt{2(\kappa_2-V_1+hW_1)}$ meet with a vertical tangent. The sector
field
\[
\Y_1=\frac{p_1\,\partial_{q_1}-(V_1'-HW_1')\,\partial_{p_1}}{W_1+W_2}
\]
nevertheless remains nonzero as long as $V_1'(q_1)\ne hW_1'(q_1)$, and under
the same condition the implicit function theorem gives $q_1$ as a smooth
function of $p_1$ on the leaf, which is precisely the graph described by the
chart $(Q_1,P_1)=(p_1,-q_1)$ in the proof. The theorem therefore applies at
every turning point except those at which the first separated motion is at
rest, that is, on the thin orbits along which $q_1$ stays constant, such as the shell orbits of a St\"ackel potential, which lie on a surface of
constant confocal radial coordinate and whose ellipse fits determine the
focal distance in \cite{Binney2014} (see also \cite{SandersBinney2016}). For
a prescribed product structure, the theorem can also be read as a concrete
realisation of the invariant-Lagrangian criterion of
\cite[Theorem~3.1 and Remark~3.8]{FalquiPedroni}.

The Frobenius formulation of Llibre and Ram\'{\i}rez
\cite[Section~4]{LlibreRamirez2020} uses the partial Hamiltonian fields
$H_{p_i}\partial_{q_i}-H_{q_i}\partial_{p_i}$, which are the sector fields
$\Y_i$ of \eqref{eq:sectorfields}, and shows, in the given canonical
coordinates, that the Levi-Civita conditions hold exactly when the brackets
$[\Y_i,\Y_j]$ close on the span of the sector fields
\cite[Theorem~10]{LlibreRamirez2020}. Their bracket identity
\cite[eq.~(20)]{LlibreRamirez2020} exhibits the components of $[\Y_i,\Y_j]$
outside this span as multiples of the Levi-Civita expressions, and the
structure coefficients of the closure are the derivatives
$\Y_j(\log|H_{p_i}|)$, or, in \cite[Remark~11]{LlibreRamirez2020},
$\Y_j(\log|H_{q_i}|)$, so that the formulation requires every $H_{p_i}$, or
every $H_{q_i}$, to be nonzero. Proposition~\ref{prop:frobenius} and
Theorem~\ref{thm:sector} differ from it in the following respects.
\begin{itemize}
	\item The obstruction is the single scalar
	$b_{ij}=\om(\Y_j,[\Y_i,\Y_j])$ for each pair of degrees of freedom, which
	is obtained without division by a derivative of $H$ and is invariant under
	canonical changes within the factors (Proposition~\ref{prop:intrinsic}).
	\item The distribution spanned by the sector fields is identified as
	Lagrangian (Proposition~\ref{prop:frobenius}(a)).
	\item The only regularity assumption is $\Y_i\ne0$, under which $H_{p_i}$
	may vanish in some factors and $H_{q_i}$ in others.
	\item The equivalence is extended to the product foliation and to the
	separated complete integral, items (c) and (d) of
	Theorem~\ref{thm:sector}, whereas in \cite{LlibreRamirez2020} the link with
	a complete integral rests on the classical theorem of Levi-Civita.
\end{itemize}

\subsection{A symplectic--Nijenhuis formulation}

Separation of variables also has a more algebraic description in the
bi-Hamiltonian framework. There, a Nijenhuis operator compatible with the
symplectic form is the primary object, and its eigenvalues, together with
suitable conjugate momenta, provide the separation coordinates
\cite{FalquiPedroni,Blaszak,TondoTempesta}. On a prescribed product
structure, a torsion-free operator of this kind is available without any
reference to the Hamiltonian, and the off-diagonal values reappear in the
exterior derivative of the one-form that it induces from $\dd H$. To make
this precise, let $\lambda_i$ be functions on the individual factors,
pairwise distinct on the domain, and set
\begin{equation}
 \mathbf L=\sum_{i=1}^n\lambda_i\Pib_i,\qquad
 \al_{\mathbf L}=\mathbf L^*\dd H,\qquad
 \{F,G\}_{\mathbf L}=\sum_{i=1}^n\lambda_i\bigl(F_{q_i}G_{p_i}-F_{p_i}G_{q_i}\bigr).
 \label{eq:nijenhuis-data}
\end{equation}

\begin{proposition}[Nijenhuis characterisation]
\label{prop:nijenhuis}
Let $\mathbf L$, $\al_{\mathbf L}$ and $\{\cdot,\cdot\}_{\mathbf L}$ be
given by \eqref{eq:nijenhuis-data}. Then the following hold.
\begin{enumerate}[label=(\alph*)]
 \item \emph{Vanishing torsion.} The Nijenhuis torsion of $\mathbf L$
 vanishes.
 \item \emph{Compatible Poisson bracket.} The bracket
 $\{\cdot,\cdot\}_{\mathbf L}$ is Poisson and compatible with
 $\{\cdot,\cdot\}$, that is, $\{\cdot,\cdot\}+t\{\cdot,\cdot\}_{\mathbf L}$
 is a Poisson bracket for every $t\in\mathbb R$.
 \item \emph{Off-diagonal values.} For $i\ne j$,
 \begin{equation}
  \dd\al_{\mathbf L}(\Y_i,\Y_j)=(\lambda_i-\lambda_j)\,b_{ij}.
  \label{eq:nijenhuis-defect}
 \end{equation}
 \item \emph{Cyclic description.} Where every $\Y_i$ is nonzero,
 \[
  \mathcal P_H=\spanop\bigl\{\Z_H,\mathbf L\Z_H,\ldots,\mathbf L^{n-1}\Z_H\bigr\}.
 \]
 \item \emph{Characterisation.} Where every $\Y_i$ is nonzero,
 \begin{equation}
  b_{ij}=0\ \text{for all }i\ne j\quad\Longleftrightarrow\quad
  \dd\al_{\mathbf L}\big|_{\mathcal P_H\times\mathcal P_H}=0
  \quad\Longleftrightarrow\quad \mathcal P_H\ \text{is involutive.}
  \label{eq:nijenhuis-equiv}
 \end{equation}
\end{enumerate}
\end{proposition}

\begin{proof}
\emph{(a)} The Nijenhuis torsion
\[
 \mathbf N_{\mathbf L}(\mathbf X,\mathbf V)
 =[\mathbf L\mathbf X,\mathbf L\mathbf V]-\mathbf L[\mathbf L\mathbf X,\mathbf V]
 -\mathbf L[\mathbf X,\mathbf L\mathbf V]+\mathbf L^2[\mathbf X,\mathbf V]
\]
is a tensor \cite{BolsinovKonyaevMatveev}, so it suffices to evaluate it on
the coordinate fields of a product chart. Let $\mathbf X\in E_i$ and
$\mathbf V\in E_j$ be two such fields, so that
$\mathbf L\mathbf X=\lambda_i\mathbf X$, $\mathbf L\mathbf V=\lambda_j\mathbf V$
and $[\mathbf X,\mathbf V]=0$. If $i\ne j$, then $\mathbf X(\lambda_j)=0$ and
$\mathbf V(\lambda_i)=0$, because each eigenvalue depends only on its own
factor, and all four terms vanish. If $i=j$, put $f=\lambda_i$; then
\[
 [f\mathbf X,f\mathbf V]=f\,\mathbf X(f)\,\mathbf V-f\,\mathbf V(f)\,\mathbf X,\qquad
 \mathbf L[f\mathbf X,\mathbf V]=-f\,\mathbf V(f)\,\mathbf X,\qquad
 \mathbf L[\mathbf X,f\mathbf V]=f\,\mathbf X(f)\,\mathbf V,
\]
and the four terms of the torsion cancel.

\emph{(b)} The bracket $\{\cdot,\cdot\}+t\{\cdot,\cdot\}_{\mathbf L}$ is
\[
\{F,G\}_t=\sum_{i=1}^n(1+t\lambda_i)\bigl(F_{q_i}G_{p_i}-F_{p_i}G_{q_i}\bigr).
\]
Where every $1+t\lambda_i$ is nonzero, this is the Poisson bracket of the
two-form
\[
\om_t=\sum_{i=1}^n\frac{\dd q_i\wedge\dd p_i}{1+t\lambda_i},
\]
which is nondegenerate, and closed because each coefficient depends only on
$(q_i,p_i)$, so that the Jacobi identity holds there \cite{Arnold}. For given
functions, the Jacobiator of $\{\cdot,\cdot\}_t$ at a point is a polynomial in
$t$ that vanishes for all sufficiently small $|t|$, and hence for every
$t\in\mathbb R$.

\emph{(c)} Put $\mathbf C_{ij}=[\Y_i,\Y_j]$, which lies in $E_i\oplus E_j$
by the proof of Proposition~\ref{prop:frobenius}(c). Since
$\dd H(\mathbf X)=\om(\Z_H,\mathbf X)$ and the planes $E_k$ are mutually
$\om$-orthogonal, \eqref{eq:intrinsic} and
Proposition~\ref{prop:intrinsic}(b) give
\[
 \dd H(\Pib_j\mathbf C_{ij})=\om(\Y_j,[\Y_i,\Y_j])=b_{ij},\qquad
 \dd H(\Pib_i\mathbf C_{ij})=-\om(\Y_i,[\Y_j,\Y_i])=-b_{ij}.
\]
Moreover, $\al_{\mathbf L}(\mathbf X)=\sum_k\lambda_k\,\dd H(\Pib_k\mathbf X)$,
so that $\al_{\mathbf L}(\Y_k)=\lambda_k\,\dd H(\Y_k)=0$ by
Proposition~\ref{prop:frobenius}(b). The formula for the exterior derivative
of a one-form then leaves only the bracket term,
\[
 \dd\al_{\mathbf L}(\Y_i,\Y_j)=-\al_{\mathbf L}(\mathbf C_{ij})
 =-\lambda_i\,\dd H(\Pib_i\mathbf C_{ij})-\lambda_j\,\dd H(\Pib_j\mathbf C_{ij})
 =(\lambda_i-\lambda_j)\,b_{ij}.
\]

\emph{(d)} Since $\mathbf L$ acts on $E_i$ as multiplication by $\lambda_i$,
\[
 \mathbf L^m\Z_H=\sum_{i=1}^n\lambda_i^m\,\Y_i,\qquad m=0,\ldots,n-1.
\]
The coefficient matrix $(\lambda_i^m)$ is a Vandermonde matrix with
determinant $\prod_{i<j}(\lambda_j-\lambda_i)\ne0$, so the fields
$\mathbf L^m\Z_H$ span the same space as $\Y_1,\ldots,\Y_n$, which is
$\mathcal P_H$.

\emph{(e)} For the first equivalence, suppose that all $b_{ij}$ vanish. Then
\eqref{eq:nijenhuis-defect} gives $\dd\al_{\mathbf L}(\Y_i,\Y_j)=0$ for
$i\ne j$, skew-symmetry gives the same for $i=j$, and bilinearity extends
this to all of $\mathcal P_H\times\mathcal P_H$. Conversely, if
$\dd\al_{\mathbf L}$ vanishes on $\mathcal P_H\times\mathcal P_H$, then
$(\lambda_i-\lambda_j)\,b_{ij}=0$ for $i\ne j$, and $b_{ij}=0$ because the
$\lambda_i$ are pairwise distinct. The second equivalence is
Proposition~\ref{prop:frobenius}(c), since the vanishing of the Frobenius
tensor is the definition of involutivity.
\end{proof}

In \eqref{eq:nijenhuis-equiv}, the two-form $\dd\al_{\mathbf L}$ is evaluated
on $\mathcal P_H$ only, and the same condition imposed on the whole tangent
bundle would exclude separated Hamiltonians. For $H=p_1+p_2$ and
$\lambda_i=q_i$ on $q_1\ne q_2$, for instance, the Hamiltonian is separated and
$b_{12}=0$, whereas
\[
\al_{\mathbf L}=q_1\,\dd p_1+q_2\,\dd p_2,\qquad
\dd\al_{\mathbf L}=\om,\qquad
\dd\al_{\mathbf L}\wedge\dd H
=\dd q_1\wedge\dd p_1\wedge\dd p_2+\dd q_2\wedge\dd p_2\wedge\dd p_1\ne0,
\]
so that the condition $\dd\al_{\mathbf L}\wedge\dd H=0$ fails on the whole
tangent bundle.
The torsion-free operator $\mathbf L$ itself is determined by the product
structure and the choice of the $\lambda_i$ alone, in agreement with the
splitting theorem for Nijenhuis operators
\cite[Theorem~3.1]{BolsinovKonyaevMatveev}, and the Hamiltonian enters the
characterisation only through \eqref{eq:nijenhuis-defect}.

\section{Exact perturbation theory and defects}
\label{sec:defects}

In practical applications, exact separability is exceptional, since perturbed
oscillators, galactic potentials and most other models of interest are at
best close to a separable system, and the relevant questions are then
how far from separation a given Hamiltonian lies and how this distance
responds to a perturbation \cite{Arnold,SandersBinney2016}. This section
therefore turns from the zero set of the operator to its size. Because the
operator is cubic in the Hamiltonian, its behaviour under perturbation can be
described without any approximation. With the probes held fixed, the
operator is the diagonal $\B[H]=\mathfrak T(H,H,H)$ of the trilinear form
\begin{equation}
 \mathfrak T(F,G,K)=\{F,c\}\,\minor_{a,b}(\{G,d\},K)-\{F,d\}\,\minor_{a,b}(\{G,c\},K),
 \label{eq:trilinear}
\end{equation}
and for $H_\varepsilon=H_0+\varepsilon V$ trilinearity gives the expansion
\begin{equation}
 \B[H_\varepsilon]=\B[H_0]+\varepsilon\,b_1+\varepsilon^2b_2+\varepsilon^3\B[V],
 \label{eq:cubic}
\end{equation}
which is exact, with no remainder term, and whose coefficients are
\begin{align*}
 b_1&=\mathfrak T(V,H_0,H_0)+\mathfrak T(H_0,V,H_0)+\mathfrak T(H_0,H_0,V),\\
 b_2&=\mathfrak T(H_0,V,V)+\mathfrak T(V,H_0,V)+\mathfrak T(V,V,H_0).
\end{align*}

For natural Hamiltonians, the three nonzero coefficients of this expansion
already appear in Levi-Civita's letter. With
$H_0=K=\frac12\sum_{k,l}a^{(kl)}(q)\,p_kp_l$, $V=U(q)$ and $\varepsilon=-1$,
the expansion describes $H=K-U$. The cubic coefficient $b_{ij}[U]$ vanishes
because $U$ does not depend on the momenta, and the coefficient of
$\varepsilon^k$ is homogeneous of degree $4-2k$ in the momenta. The three
groups into which Levi-Civita splits his conditions according to this degree
\cite[p.~386]{LeviCivita} are therefore, up to sign,
\[
\mathrm{(I)}\;\; b_{ij}[K]=0,\qquad
\mathrm{(II)}\;\; b_{1,ij}=0,\qquad
\mathrm{(III)}\;\; b_{2,ij}=-a^{(ij)}\,U_{q_i}U_{q_j}=0\qquad(i\ne j),
\]
where $b_{1,ij}$ and $b_{2,ij}$ are the first- and second-order coefficients
for the probes $(q_i,p_i;q_j,p_j)$. In the notation of
Section~\ref{sec:linearised}, his second group is $\Lop_K(U)=0$, the
linearised operator of the kinetic energy applied to the potential.

The first-order coefficient is the subject of the rest of this section.
Section~\ref{sec:linearised} defines the linearised operator and computes it
for the polar-nodal normal forms, Section~\ref{sec:defect-functional} turns
the off-diagonal values into a measure of the departure from separation,
Section~\ref{sec:homoclinic} relates the linearised operator to the splitting
of separatrices, and Section~\ref{sec:levelset} treats Hamiltonians that
separate only on their level sets, with Kerr geodesics as the example.

\subsection{The linearised operator}
\label{sec:linearised}

On a prescribed product structure, the first-order coefficients for the
off-diagonal probes define the \emph{linearised separability operator}, which
assigns to a perturbation its first-order off-diagonal defects,
\begin{equation}
 \Lop_{H_0}(V)_{ij}=\frac{\partial}{\partial\varepsilon}b_{ij}[H_0+\varepsilon V]\Big|_{\varepsilon=0},
 \qquad i<j.
 \label{eq:linearised}
\end{equation}
It is linear in $V$ and of differential order at most two. When $H_0$ is
separable, its kernel consists of the perturbations that preserve separation
to first order, whereas exact preservation along the line
$H_0+\varepsilon V$ requires, in addition, that the quadratic and cubic
coefficients in \eqref{eq:cubic} vanish for the off-diagonal probes. Two
families of kernel directions are always present, and both originate in
symmetries of the problem, namely reparametrisations of the energy and
canonical flows that preserve the factors.

\begin{proposition}[Symmetry directions in the kernel]
\label{prop:kernel}
Suppose that $b_{ij}[H_0]=0$ for all $i\ne j$. Then the following hold.
\begin{enumerate}[label=(\alph*)]
 \item \emph{Energy reparametrisations.} For every smooth function $f$ of one
variable and every $\varepsilon\in\mathbb R$, the off-diagonal values of
$H_0+\varepsilon f(H_0)$ vanish, $b_{ij}[H_0+\varepsilon f(H_0)]=0$ for all
$i\ne j$, and in particular $\Lop_{H_0}(f(H_0))=0$.
 \item \emph{Factor-preserving canonical flows.}
 $\Lop_{H_0}(\{H_0,G\})=0$ for every $G=\sum_iG_i(q_i,p_i)$.
\end{enumerate}
\end{proposition}

\begin{proof}
\emph{(a)} The perturbed Hamiltonian is $\chi_\varepsilon(H_0)$ with
$\chi_\varepsilon(t)=t+\varepsilon f(t)$, so \eqref{eq:energychange} gives
\[
 b_{ij}[H_0+\varepsilon f(H_0)]=\chi_\varepsilon'(H_0)^3\,b_{ij}[H_0]=0
\]
for every $\varepsilon$, and differentiating at $\varepsilon=0$ gives the
first claim.

\emph{(b)} Each $G_i$ depends only on the $i$th factor, so the Hamiltonian
flow $\Phi_\varepsilon$ of $G$ is, locally, the product of the flows of the
$G_i$ on the factors $M_i$. It is therefore a symplectomorphism that maps
each factor to itself, that is, a canonical change of chart within the
factors. The covariance \eqref{eq:covariance}, applied with the probes
$q_i,p_i,q_j,p_j$, together with Proposition~\ref{prop:intrinsic}(c), gives
\[
 b_{ij}[\Phi_\varepsilon^*H_0]=\Phi_\varepsilon^*\,b_{ij}[H_0]=0
\]
for all small $\varepsilon$. Since $b_{ij}$ is a polynomial in the
derivatives of the Hamiltonian, its derivative along this curve at
$\varepsilon=0$ depends only on the first-order term
$\partial_\varepsilon\Phi_\varepsilon^*H_0|_{\varepsilon=0}=\Z_G(H_0)=\{H_0,G\}$,
and it equals $\Lop_{H_0}(\{H_0,G\})_{ij}$, which therefore vanishes.
\end{proof}

Beyond these symmetry directions, the linearised operator and the full
expansion \eqref{eq:cubic} can be computed in closed form for a family of
Hamiltonians that are not natural and that is central in satellite theory.
In the polar-nodal (Hill) variables $(r,\theta,\nu,R,\Theta,N)$, $r$ denotes
the radial distance, $\theta$ the argument of latitude, $\nu$ the right
ascension of the ascending node, $R$ the radial velocity, $\Theta$ the
modulus of the angular momentum and $N$ its polar component, and these
variables form a canonical chart on the open set where $r>0$ and
$\Theta>|N|$ \cite{Deprit1981}. The product structure considered here has
the three factors $(r,R)$, $(\theta,\Theta)$ and $(\nu,N)$, and the
Hamiltonians of the form
\begin{equation}
	H=\tfrac12R^2+K(r,\theta,\Theta,N),
	\label{eq:pn-H}
\end{equation}
with $K$ smooth, include the zonal problem of an artificial satellite
together with its normal forms and integrable intermediaries
\cite{CidLahulla1969,Deprit1981,AlfriendCoffey1984,Lara2021}. Since $K$ may depend on $\Theta$ in an arbitrary way, for
instance through the inclination, such Hamiltonians are in general not
natural, and the algebraic criteria of \cite{Waksjo,Benenti2016} do not
apply to them. On the angular factor, the bracket
\[
\{A,B\}_{\theta\Theta}
=\frac{\partial A}{\partial\theta}\frac{\partial B}{\partial\Theta}
-\frac{\partial A}{\partial\Theta}\frac{\partial B}{\partial\theta}
\]
denotes the Poisson bracket in the pair $(\theta,\Theta)$, in which $r$ and
$N$ enter as parameters. The integrable reference Hamiltonian is
\begin{equation}
	H_0=\tfrac12R^2-\frac{\mu}{r}+\frac{G}{2r^2},\qquad
	G=\Theta^2+2f(\theta,\Theta,N),
	\label{eq:pn-H0}
\end{equation}
where $\mu>0$ is the gravitational parameter and $f$ is a smooth function,
in the applications typically a trigonometric polynomial in $\theta$ with
coefficients depending on $\Theta$ and $N$. The term $r^{-2}f$ thus modifies
the centrifugal term of the Kepler problem, and $G$, which reduces to the
squared angular momentum when $f=0$, will be called the modified angular
integral.

\begin{proposition}[Polar-nodal normal forms]
	\label{prop:polar-nodal}
	Let $H$ be of the form \eqref{eq:pn-H} on a domain of the polar-nodal chart.
	Then the following hold.
	\begin{enumerate}[label=(\alph*)]
		\item \emph{Defect.} The off-diagonal values involving the factor $(\nu,N)$
		vanish identically, and the radial--angular value is
		\begin{equation}
			b_{r\theta}[H]=R\,\{K,\partial_rK\}_{\theta\Theta}.
			\label{eq:pn-defect}
		\end{equation}
		\item \emph{Integrable family.} The reference Hamiltonian \eqref{eq:pn-H0}
		satisfies $b_{r\theta}[H_0]\equiv0$, and $H_0$, $G$ and $N$ are first
		integrals of $H_0$ in involution.
		\item \emph{Exact perturbation defect.} For a perturbation
		$V=\rho(r)\,g(\theta,\Theta,N)$, with $\rho$ and $g$ smooth, and for every
		$\varepsilon$,
		\begin{equation}
			b_{r\theta}[H_0+\varepsilon V]=\varepsilon\,\Lop_{H_0}(V)_{r\theta},\qquad
			\Lop_{H_0}(V)_{r\theta}
			=\frac{(r^2\rho)'}{2r^4}\,R\,\{G,g\}_{\theta\Theta},
			\label{eq:pn-linear}
		\end{equation}
		so that the expansion \eqref{eq:cubic} terminates at first order.
		\item \emph{Kernel.} Let $U$ be the set of points with $r\in I$, $R\in J$
		and $(\theta,\Theta,N)\in W$, where $I$ and $J$ are open intervals and $W$
		is an open set. Then $\Lop_{H_0}(V)_{r\theta}$ vanishes on $U$ if and only
		if $r^2\rho$ is constant on $I$ or $\{G,g\}_{\theta\Theta}$ vanishes on $W$.
	\end{enumerate}
\end{proposition}

\begin{proof}
	\emph{(a)} For the pair $(q_i,p_i)=(r,R)$ and $(q_j,p_j)=(\theta,\Theta)$,
	the Hamiltonian \eqref{eq:pn-H} gives $H_R=R$ and $H_{R\theta}=H_{R\Theta}=0$,
	so only two of the four terms in \eqref{eq:LC} survive, and
	\[
	b_{r\theta}[H]=-H_RH_\Theta H_{r\theta}+H_RH_\theta H_{r\Theta}
	=R\bigl(K_\theta\,\partial_rK_\Theta-K_\Theta\,\partial_rK_\theta\bigr),
	\]
	which is \eqref{eq:pn-defect}. For the pairs that involve $(\nu,N)$, the
	Hamiltonian does not depend on $\nu$, so that $H_\nu$ and every second
	derivative containing $\nu$ vanish, and each of the four terms in
	\eqref{eq:LC} contains one such factor.
	
	\emph{(b)} The functions
	\[
	K_0=-\frac{\mu}{r}+\frac{G}{2r^2},\qquad
	\partial_rK_0=\frac{\mu}{r^2}-\frac{G}{r^3}
	\]
	depend on $(\theta,\Theta)$ only through $G$, so their angular bracket
	vanishes, and (a) gives $b_{r\theta}[H_0]\equiv0$. The function $G$ depends
	neither on $(r,R)$ nor on $\nu$, so that $\{G,H_0\}$ reduces to the angular
	bracket of $G$ with a function of $G$, which vanishes. Finally, $N$ is a
	first integral because $\nu$ is cyclic, and $\{G,N\}=0$ because $G$ does not
	depend on $\nu$.
	
	\emph{(c)} With $K=K_0+\varepsilon\rho g$, the radial derivative is
	$\partial_rK=\partial_rK_0+\varepsilon\rho'g$. Bilinearity of the bracket
	together with $\{g,g\}_{\theta\Theta}=0$ gives
	\[
	\{K,\partial_rK\}_{\theta\Theta}
	=\varepsilon\rho'\{K_0,g\}_{\theta\Theta}+\varepsilon\rho\{g,\partial_rK_0\}_{\theta\Theta}
	=\varepsilon\Bigl(\frac{\rho'}{2r^2}+\frac{\rho}{r^3}\Bigr)\{G,g\}_{\theta\Theta}
	=\varepsilon\,\frac{(r^2\rho)'}{2r^4}\,\{G,g\}_{\theta\Theta},
	\]
	and (a) yields the exact defect. Since it is linear in $\varepsilon$,
	comparison with \eqref{eq:cubic} identifies its coefficient with
	$\Lop_{H_0}(V)_{r\theta}$, as defined in \eqref{eq:linearised}, and shows
	that the quadratic and cubic coefficients vanish.
	
	\emph{(d)} ($\Leftarrow$) If $r^2\rho$ is constant on $I$, the radial factor
	in \eqref{eq:pn-linear} vanishes on $U$, and if $\{G,g\}_{\theta\Theta}$
	vanishes on $W$, so does the angular factor.
	
	($\Rightarrow$) Suppose that $(r^2\rho)'(r_0)\ne0$ for some $r_0\in I$. The
	open interval $J$ contains some $R_0\ne0$, and the vanishing of
	\eqref{eq:pn-linear} at the points with $r=r_0$, $R=R_0$ and
	$(\theta,\Theta,N)\in W$ then forces $\{G,g\}_{\theta\Theta}=0$ on $W$.
\end{proof}

The two alternatives in (d) correspond to two familiar ways of keeping a
normal form integrable. A perturbation with $r^2\rho=\kappa$ constant is
absorbed into the reference family by the replacement of $f$ with
$f+\varepsilon\kappa g$, and a perturbation in involution with $G$, for
instance $g=\varphi(G,N)$ for a smooth function $\varphi$, leaves $G$ a first
integral of the perturbed problem. Every other perturbation of the form
$\rho g$ destroys the separation on $U$, and its defect then factors into the
radial momentum, a radial profile that measures the departure of $\rho$ from
the inverse-square law and the angular bracket $\{G,g\}_{\theta\Theta}$. The
same bracket governs the evolution of the modified angular integral, since
along the flow of $H_0+\varepsilon V$
\[
\dot G=\{G,H_0+\varepsilon V\}=\varepsilon\,\rho\,\{G,g\}_{\theta\Theta},
\]
and the integral of this rate along homoclinic orbits of the unperturbed
problem, when such orbits exist, is the Melnikov function of the
perturbation \cite{HolmesMarsden1982}. The defect and the Melnikov integrand
therefore share the angular factor, and Section~\ref{sec:homoclinic} shows
that this is an instance of an exact identity between the linearised
operator and the Melnikov integrand, valid for every Hamiltonian with a
cyclic angle (Remark~\ref{rem:pn}). The consequences of this identity for
the normal forms of satellite theory are developed in a companion paper.

\subsection{The defect functional}
\label{sec:defect-functional}

The off-diagonal values vanish exactly where the Hamiltonian separates, so a
natural measure of how far a Hamiltonian is from separating is a weighted
$L^2$ norm of these values over a region of phase space. This measure, which is introduced here, is related to the least-squares fits
of a separability residual along an orbit used in the St\"ackel-fudge
literature \cite[Section~A.5.1]{VasilievDoc}, and Section~\ref{sec:stackel}
shows how such fits relate to the functional defined below.

\begin{definition}[Defect functional]
	\label{def:defect}
	Let $U$ be a relatively compact open domain, let $w_{ij}>0$, $i<j$, be smooth
	weights on its closure, and let $\dd\mu=\om^n/n!$ be the symplectic volume.
	The \emph{defect functional} of $H$ on $U$ with weights $w=(w_{ij})$ is defined as
	\begin{equation}
		\mathcal E_{U,w}(H)=\int_U\sum_{i<j}w_{ij}\,b_{ij}[H]^2\,\dd\mu.
		\label{eq:defect}
	\end{equation}
\end{definition}

\begin{proposition}[Properties of the defect]
\label{prop:defect}
The defect functional \eqref{eq:defect} has the following properties.
\begin{enumerate}[label=(\alph*)]
 \item \emph{Chart independence.} If the domain $U$ and the weights $w_{ij}$
are fixed as a subset of $M$ and functions on it, independently of the
canonical charts on the factors, then $\mathcal E_{U,w}(H)$ takes the same
value for every choice of these charts.
 \item \emph{Zero set.} $\mathcal E_{U,w}(H)=0$ if and only if
 $b_{ij}[H]\equiv0$ on $U$ for all $i<j$.
 \item \emph{Energy reparametrisation.} For every smooth function $\chi$
 with $\chi'\ne0$ on $H(U)$,
 \[
  \mathcal E_{U,w}(\chi(H))=\mathcal E_{U,\chi'(H)^6w}(H),
 \]
 so the zero set of the defect is invariant under reparametrisation, whereas
 its value is not.
 \item \emph{Perturbation expansion.} If $b_{ij}[H_0]=0$ on $U$ for all
 $i<j$, then $\mathcal E_{U,w}(H_0+\varepsilon V)$ is a polynomial of degree
 at most six in $\varepsilon$, with
 \begin{equation}
  \mathcal E_{U,w}(H_0+\varepsilon V)
  =\varepsilon^2\int_U\sum_{i<j}w_{ij}\bigl[\Lop_{H_0}(V)_{ij}\bigr]^2\dd\mu+O(\varepsilon^3),
  \label{eq:defect-leading}
 \end{equation}
 and $\mathcal E_{U,w}(H_0+\varepsilon V)=O(\varepsilon^4)$ when
 $\Lop_{H_0}(V)=0$.
\end{enumerate}
\end{proposition}

\begin{proof}
\emph{(a)} The off-diagonal values are invariant by
Proposition~\ref{prop:intrinsic}(c), the volume $\dd\mu$ is defined by $\om$
alone, and $U$ and $w$ do not depend on the charts.

\emph{(b)} If all $b_{ij}$ vanish on $U$, the integrand vanishes and so does
$\mathcal E_{U,w}(H)$. Conversely, the integrand is continuous and
nonnegative; if it were positive at a point of $U$, it would be positive on
a neighbourhood of that point, which has positive measure, and the integral
would be positive. Hence $\mathcal E_{U,w}(H)=0$ forces the integrand to
vanish on $U$, and since every $w_{ij}$ is positive, every $b_{ij}$ vanishes
there.

\emph{(c)} By \eqref{eq:energychange},
$b_{ij}[\chi(H)]^2=\chi'(H)^6\,b_{ij}[H]^2$, which gives the identity. Since
$\chi'(H)^6w_{ij}$ is again a smooth positive weight, (b) shows that the two
defects vanish together.

\emph{(d)} For the off-diagonal probes, \eqref{eq:cubic} with
$b_{ij}[H_0]=0$ reads
\[
 b_{ij}[H_0+\varepsilon V]=\varepsilon\,\Lop_{H_0}(V)_{ij}
 +\varepsilon^2\,b_{2,ij}+\varepsilon^3\,b_{ij}[V],
\]
where $b_{2,ij}$ denotes the quadratic coefficient. Its square is a
polynomial of degree six in $\varepsilon$ whose lowest-order term is
$\varepsilon^2\Lop_{H_0}(V)_{ij}^2$, and whose terms of order two and three
vanish when $\Lop_{H_0}(V)_{ij}=0$. Integrating against
$w_{ij}\,\dd\mu$ gives all three statements.
\end{proof}

For a prescribed family of product structures, the defect can be minimised
over the family parameters. The domain and the weights must then refer to
the same physical points for every member of the family, since a box fixed
in the trial coordinates would itself move along the family and bias the
objective. The defect measures how far $H$ is from separating in the chosen
coordinates, which is a different quantity from the errors in trajectories
or actions of an approximate separable model, and the relation between the
two is taken up in Sections~\ref{sec:stackel} and~\ref{sec:limitations}.

\subsection{Defects and homoclinic splitting}
\label{sec:homoclinic}

Chaos appears first near the separatrices of an integrable system, and the
Melnikov method measures the splitting of a separatrix under a perturbation
to first order \cite{HolmesMarsden1982,GuckenheimerHolmes}. Holmes and
Marsden \cite{HolmesMarsden1982} developed the method for perturbations of
Hamiltonians in which the actions conjugate to cyclic angles enter a system
with a homoclinic orbit, and for two degrees of freedom, the case considered
in this subsection, it is shown below
that the linearised operator is, up to an explicit factor, an angle
derivative of the Melnikov integrand. The identity
\eqref{eq:rotor-identity}, Theorem~\ref{thm:melnikov} and
Corollary~\ref{cor:dichotomy} are introduced here, whereas the Melnikov
function and its relation to the splitting distance are standard
\cite{HolmesMarsden1982}. The class considered is
\begin{equation}
 H=H_0(q,p,I)+\varepsilon V(q,p,\theta,I),\qquad
 \nu=\frac{\partial H_0}{\partial I},
 \label{eq:rotor}
\end{equation}
on the product of the pairs $(q,p)$ and $(\theta,I)$, where $\theta$ is an angle
of period $2\pi$ that is cyclic for $H_0$, $I$ is its conjugate action, and
$\nu(q,p,I)$ will be called the \emph{rotor frequency}. The class contains
every integrable system with two degrees of freedom written in action--angle
variables for one of them. The sector field of $H_0$ on the first factor is
\[
 \Y_1=\frac{\partial H_0}{\partial p}\,\partial_q-\frac{\partial H_0}{\partial q}\,\partial_p,
 \qquad\text{so that}\qquad \Y_1F=-\{H_0,F\}_{qp}
\]
for every function $F$, where $\{\cdot,\cdot\}_{qp}$ denotes the bracket in the
pair $(q,p)$, in which $\theta$ and $I$ enter as parameters.

\begin{proposition}[Defect of Hamiltonians with a cyclic angle]
\label{prop:rotor}
For the Hamiltonian \eqref{eq:rotor} the following hold.
\begin{enumerate}[label=(\alph*)]
 \item $b_{q\theta}[H_0]\equiv0$.
 \item For every smooth $V$,
 \begin{equation}
  \Lop_{H_0}(V)_{q\theta}=-\nu\,\Y_1\Bigl(\frac{\partial V}{\partial\theta}\Bigr)
  +\frac{\partial V}{\partial\theta}\,\Y_1(\nu).
  \label{eq:rotor-L}
 \end{equation}
 \item Where $\nu\ne0$, the linearised operator is an angle derivative,
 \begin{equation}
  \nu^{-2}\,\Lop_{H_0}(V)_{q\theta}
  =\frac{\partial}{\partial\theta}\,\bigl\{H_0,\nu^{-1}V\bigr\}_{qp}.
  \label{eq:rotor-identity}
 \end{equation}
\end{enumerate}
\end{proposition}

\begin{proof}
With $a=q$, $b=p$, $c=\theta$ and $d=I$, the derivations in
\eqref{eq:bordered} are $D_q=\partial_p$, $D_p=-\partial_q$, $D_\theta=\partial_I$
and $D_I=-\partial_\theta$. After a change of sign of one row and one column,
and with each pair written in the order (coordinate, momentum),
\eqref{eq:bordered} reads
\[
 b_{q\theta}[H]=\det\begin{pmatrix}
  0 & H_\theta & H_I\\
  H_q & H_{q\theta} & H_{qI}\\
  H_p & H_{p\theta} & H_{pI}
 \end{pmatrix}.
\]

\emph{(a)} For $H_0$ the second column vanishes.

\emph{(b)} For $H_0+\varepsilon V$ the second column equals
$\varepsilon\,\partial_\theta(V,V_q,V_p)^{\mathsf T}$, so the first-order
coefficient is obtained by evaluating the other two columns on $H_0$. These
columns do not depend on $\theta$, and the multilinearity of the determinant
gives
\begin{equation}
 \Lop_{H_0}(V)_{q\theta}=\frac{\partial}{\partial\theta}\det\begin{pmatrix}
  0 & V & \nu\\
  H_{0,q} & V_q & \nu_q\\
  H_{0,p} & V_p & \nu_p
 \end{pmatrix}
 =\frac{\partial}{\partial\theta}\Bigl(\nu\,\{H_0,V\}_{qp}-V\,\{H_0,\nu\}_{qp}\Bigr),
 \label{eq:rotor-det}
\end{equation}
where $\nu_q=H_{0,qI}$ and $\nu_p=H_{0,pI}$, and the second equality is the
expansion along the first row. Since $\nu$ does not depend on $\theta$ and
$\Y_1F=-\{H_0,F\}_{qp}$, this is \eqref{eq:rotor-L}.

\emph{(c)} The expression in brackets in \eqref{eq:rotor-det} equals
$\nu^2\{H_0,\nu^{-1}V\}_{qp}$, and $\nu$ does not depend on $\theta$.
\end{proof}

In this form, the Melnikov integrand of \eqref{eq:M} below appears directly
inside the operator. The cyclic angle confines the dependence on $\theta$ to
one column of the bordered determinant, the angle derivative leaves the
determinant by multilinearity, and the determinant that remains, with $V$
in place of $V_\theta$, is $\nu^2$ times the integrand $\{H_0,\nu^{-1}V\}_{qp}$.

When $H_0=H_0(I_1,I_2)$ is written in action--angle variables, the choice
$(q,p)=(\varphi_1,I_1)$ and $(\theta,I)=(\varphi_2,I_2)$ gives
$\Y_1=\omega_1\partial_{\varphi_1}$ and $\Y_1(\nu)=0$, where the frequencies
are $\omega_i=\partial_{I_i}H_0$, and Proposition~\ref{prop:rotor} reduces to a
statement about Fourier coefficients.

\begin{corollary}[Action--angle symbol]
\label{cor:aa}
Let $H_0=H_0(I_1,I_2)$ and let $V=\sum_{k\in\mathbb Z^2}V_k(I)\,e^{\mathrm i
k\cdot\varphi}$ be the Fourier series of a smooth perturbation. Then
\begin{equation}
 \Lop_{H_0}(V)_{12}=-\omega_1\omega_2\,\frac{\partial^2V}{\partial\varphi_1\partial\varphi_2}
 =\sum_{k\in\mathbb Z^2}\omega_1\omega_2\,k_1k_2\,V_k(I)\,e^{\mathrm ik\cdot\varphi},
 \label{eq:aa}
\end{equation}
and on every open set of actions where $\omega_1\omega_2\ne0$,
\[
 \Lop_{H_0}(V)_{12}\equiv0\iff V_k=0\ \text{whenever}\ k_1k_2\ne0
 \iff V=V_1(\varphi_1,I)+V_2(\varphi_2,I).
\]
\end{corollary}

\begin{proof}
The first equality in \eqref{eq:aa} is \eqref{eq:rotor-L} with
$\Y_1=\omega_1\partial_{\varphi_1}$ and $\Y_1(\nu)=0$, and termwise
differentiation gives the second. For the first equivalence, ($\Rightarrow$)
follows from the independence of the exponentials together with
$\omega_1\omega_2\ne0$, and ($\Leftarrow$) holds because every remaining
harmonic has $k_1k_2=0$. The second equivalence groups the harmonics with
$k_2=0$ and those with $k_1=0$.
\end{proof}

Corollary~\ref{cor:aa} connects the linearised operator with Poincar\'e's
theorem on the non-existence of analytic first integrals
\cite{Poincare1892,Kozlov1983}. The Poincar\'e set of a perturbation consists
of the harmonics $V_k\ne0$ whose resonance $k\cdot\omega=0$ is met in the
domain, and a harmonic with $k_1k_2=0$ can only resonate where a frequency
vanishes, so on a domain where $\omega_1\omega_2\ne0$ the Poincar\'e set is
contained in the spectral support of $\Lop_{H_0}(V)_{12}$. On a resonant torus
the corresponding coefficient in \eqref{eq:aa} equals $-(k_1\omega_1)^2V_k$,
since $\omega_2k_2=-\omega_1k_1$ there.

The Melnikov setting requires the following hypotheses on
\eqref{eq:rotor}.
\begin{enumerate}[label=(H\arabic*)]
 \item At an action $I=I_h$ the system $H_0(\cdot,\cdot,I_h)$ with one degree
 of freedom has a hyperbolic equilibrium $x_h$, with Lyapunov exponent
 $\lambda>0$, and a homoclinic orbit $\gamma(t)$, $t\in\mathbb R$.
 \item The rotor frequency does not vanish on the closure of $\gamma$, and
 $\nu_h=\nu(x_h,I_h)$ denotes its value at the equilibrium.
 \item $V$ is smooth and $2\pi$-periodic in $\theta$.
\end{enumerate}
In the full system, $x_h\times\{I=I_h\}$ is then a hyperbolic periodic orbit
on the energy surface $\{H_0=E\}$ with $E=H_0(x_h,I_h)$, along which $\theta$
rotates with frequency $\nu_h$, and the solutions on its homoclinic manifold
are $(\gamma(t),\theta_0+\Phi(t),I_h)$, labelled by the phase $\theta_0$, with
\[
 \Phi(t)=\int_0^t\nu\bigl(\gamma(s),I_h\bigr)\,\dd s .
\]
By (H2) and the implicit function theorem, the equation $H_0(q,p,\mathcal
I_0)=E$ defines near $\gamma$ a smooth function $\mathcal I_0(q,p)$, which
will be called the \emph{reduced action}, and $-\mathcal I_0$ is the
unperturbed Hamiltonian of the reduced system in which $\theta$ plays the
role of time \cite{HolmesMarsden1982}. Along the homoclinic manifold, three
functions of the phase are considered, namely the Melnikov function $M$ in
the form derived below, the integrated defect $D_\gamma$, which is introduced
here, and the Melnikov potential $\Psi$,
\begin{align}
 M(\theta_0)&=\int_{-\infty}^{\infty}\bigl\{H_0,\nu^{-1}V\bigr\}_{qp}
 \bigl(\gamma(t),\theta_0+\Phi(t),I_h\bigr)\,\dd t,\label{eq:M}\\
 D_\gamma(\theta_0)&=\int_{-\infty}^{\infty}\nu^{-2}\,\Lop_{H_0}(V)_{q\theta}
 \bigl(\gamma(t),\theta_0+\Phi(t),I_h\bigr)\,\dd t,\label{eq:Dg}\\
 \Psi(\theta_0)&=\int_{-\infty}^{\infty}V\bigl(\gamma(t),\theta_0+\Phi(t),I_h\bigr)\,\dd t.
 \label{eq:Psi}
\end{align}

\begin{theorem}[Melnikov--defect identity]
\label{thm:melnikov}
Under (H1)--(H3) the following hold.
\begin{enumerate}[label=(\alph*)]
 \item The integrals \eqref{eq:M} and \eqref{eq:Dg} converge absolutely, and
 $M$ is smooth, $2\pi$-periodic and of zero mean.
 \item $M'(\theta_0)=D_\gamma(\theta_0)$ for every $\theta_0$.
 \item If $V(x_h,\theta,I_h)=0$ for every $\theta$, then \eqref{eq:Psi}
 converges, $M=\Psi'$ and $D_\gamma=\Psi''$.
 \item $M$ is the Melnikov function of the reduced system
 \cite{HolmesMarsden1982}. Along a solution of $H$ on $\{H=E\}$ that follows
 $\gamma$ with phase $\theta_0$ during one passage, the reduced action changes
 by
 \begin{equation}
  \Delta\mathcal I_0=-\varepsilon\,M(\theta_0)+O(\varepsilon^2),
  \label{eq:jump}
 \end{equation}
 and in the section $\theta=\theta_0$ the stable and unstable manifolds of the
 perturbed periodic orbit are separated near $\gamma(0)$ by the signed
 distance $\varepsilon M(\theta_0)\,|\nabla\mathcal I_0(\gamma(0))|^{-1}+O(\varepsilon^2)$.
\end{enumerate}
\end{theorem}

\begin{proof}
\emph{(a)} The integrand of \eqref{eq:M} is $-\Y_1(\nu^{-1}V)$, and $\Y_1$
vanishes at $x_h$. Since $\gamma(t)$ tends to $x_h$ at the exponential rate
$\lambda$ as $t\to\pm\infty$ and $\nu^{-1}V$ is smooth near the closure of
$\gamma$, the integrand and all its derivatives in $\theta_0$ are bounded by a
constant multiple of $e^{-\lambda|t|}$, uniformly in $\theta_0$, which gives
absolute convergence, smoothness and periodicity, and allows differentiation
under the integral sign. The integrand of \eqref{eq:Dg} is the
$\theta_0$-derivative of that of \eqref{eq:M} by \eqref{eq:rotor-identity}, so
the same bound applies. For the mean, let $\bar V(q,p,I)$ denote the average
of $V$ over $\theta$. Since $\nu^{-1}\bar V$ does not depend on $\theta$, its
derivative along $\gamma$ is $\Y_1(\nu^{-1}\bar V)$, and averaging \eqref{eq:M}
over $\theta_0$ gives
\[
 \frac{1}{2\pi}\int_0^{2\pi}M\,\dd\theta_0
 =-\int_{-\infty}^{\infty}\Y_1\bigl(\nu^{-1}\bar V\bigr)\bigl(\gamma(t)\bigr)\,\dd t
 =-\Bigl[\nu^{-1}\bar V\bigl(\gamma(t)\bigr)\Bigr]_{t=-\infty}^{t=\infty}=0,
\]
because both limits equal the value at $x_h$.

\emph{(b)} This is the differentiation under the integral sign justified in
(a), combined with \eqref{eq:rotor-identity}.

\emph{(c)} Along the curve $t\mapsto(\gamma(t),\theta_0+\Phi(t))$,
\[
 \frac{\dd}{\dd t}\bigl(\nu^{-1}V\bigr)
 =\Y_1\bigl(\nu^{-1}V\bigr)+\nu\,\frac{\partial}{\partial\theta}\bigl(\nu^{-1}V\bigr)
 =-\bigl\{H_0,\nu^{-1}V\bigr\}_{qp}+\frac{\partial V}{\partial\theta}.
\]
If $V$ vanishes at $x_h$ for every $\theta$, then $V$ and $\partial_\theta V$ are of
order $|\gamma(t)-x_h|$, hence exponentially decaying along $\gamma$, so
\eqref{eq:Psi} converges and the boundary terms of $\nu^{-1}V$ vanish.
Integration gives $M=\int\partial_\theta V\,\dd t=\Psi'$, and (b) gives
$D_\gamma=\Psi''$.

\emph{(d)} On $\{H=E\}$ the implicit function theorem gives
$I=\mathcal I(q,p,\theta;\varepsilon)$ with
$\mathcal I=\mathcal I_0-\varepsilon\,\nu^{-1}V+O(\varepsilon^2)$, and $-\mathcal I$
is the reduced Hamiltonian with time $\theta$ \cite{HolmesMarsden1982}.
Differentiating $H_0(q,p,\mathcal I_0)=E$ gives
$\partial_q\mathcal I_0=-\nu^{-1}\partial_qH_0$ and
$\partial_p\mathcal I_0=-\nu^{-1}\partial_pH_0$, hence
$\{\mathcal I_0,F\}_{qp}=-\nu^{-1}\{H_0,F\}_{qp}$ for every $F$, where the
dependence of $F$ on $q$ and $p$ through $I=\mathcal I_0(q,p)$ does not
contribute because $\{H_0,\mathcal I_0\}_{qp}=0$. Along a solution of the
reduced system,
\[
 \frac{\dd\mathcal I_0}{\dd\theta}=\{\mathcal I_0,-\mathcal I\}_{qp}
 =\varepsilon\,\{\mathcal I_0,\nu^{-1}V\}_{qp}+O(\varepsilon^2)
 =-\varepsilon\,\nu^{-1}\bigl\{H_0,\nu^{-1}V\bigr\}_{qp}+O(\varepsilon^2),
\]
and the substitution $\dd\theta=\nu\,\dd t$, with the right-hand side
evaluated on the unperturbed solution, gives \eqref{eq:jump}. The same
computation shows that $M(\theta_0)$ equals the integral
$\int\{-\mathcal I_0,\nu^{-1}V\}_{qp}\,\dd\theta$ of the bracket of the
unperturbed reduced Hamiltonian and its first-order perturbation along the
homoclinic orbit, which is the Melnikov function of
\cite{HolmesMarsden1982}, and the distance formula is the standard Melnikov
estimate for this function.
\end{proof}

\begin{corollary}[First-order dichotomy]
\label{cor:dichotomy}
Under (H1)--(H3) the following hold.
\begin{enumerate}[label=(\alph*)]
 \item If $\Lop_{H_0}(V)_{q\theta}$ vanishes on $\gamma\times S^1$, then $M\equiv0$,
 and the splitting of the separatrix is of order $\varepsilon^2$.
 \item If $D_\gamma\not\equiv0$, then $M$ changes sign at least twice in every
 period. For all sufficiently small $\varepsilon\ne0$ the stable and unstable
 manifolds of the perturbed periodic orbit then cross topologically in every
 section $\theta=\mathrm{const}$, and the section map has positive topological
 entropy.
 \item If, in addition, $D_\gamma(\theta_*)\ne0$ at a zero $\theta_*$ of $M$, the
 crossing at the corresponding phase is transverse, and the section map has
 a hyperbolic invariant set on which it is conjugate to a shift.
\end{enumerate}
\end{corollary}

\begin{proof}
\emph{(a)} By Theorem~\ref{thm:melnikov}(b), $M'=D_\gamma=0$, so $M$ is
constant, and Theorem~\ref{thm:melnikov}(a) shows that the constant is zero.

\emph{(b)} If $D_\gamma\not\equiv0$, then $M$ is not constant, and a
continuous periodic function of zero mean that does not vanish identically
takes both signs, hence changes sign at least twice per period. Once
$\varepsilon$ is small enough for the first-order term in
Theorem~\ref{thm:melnikov}(d) to dominate, the signed distance between the
manifolds changes sign between two phases at which $M$ has opposite signs,
so the manifolds cross topologically. A topological crossing of the stable
and unstable manifolds of a hyperbolic periodic point of a surface
diffeomorphism implies positive topological entropy \cite{BurnsWeiss1995}.

\emph{(c)} A simple zero of the Melnikov function gives a transverse
homoclinic point for small $\varepsilon\ne0$ \cite{HolmesMarsden1982}, and the
Smale--Birkhoff theorem gives the hyperbolic set \cite{GuckenheimerHolmes}.
\end{proof}

Pointwise values of the defect carry no information on stability, as the
two quadratic Hamiltonians compared in Section~\ref{sec:limitations} show. Its integral $D_\gamma$ along
a homoclinic orbit carries such information at first order: by
Corollary~\ref{cor:dichotomy}, if $D_\gamma\not\equiv0$ the stable and unstable
manifolds cross and the section map has positive topological entropy, and if
the linearised defect vanishes along the orbit the splitting is of order
$\varepsilon^2$. The
integrated defect also controls the width of the resulting chaotic layer.
With $h=\mathcal I_0-I_h$, the time spent near $x_h$ by an orbit at distance $h$
from the separatrix grows like $-\lambda^{-1}\ln|h|$, and the separatrix map
$h_{n+1}=h_n-\varepsilon M(\theta_n)$,
$\theta_{n+1}=\theta_n-\nu_h\lambda^{-1}\ln|h_{n+1}|+\mathrm{const}$
\cite{Chirikov1979} stretches phases by more than a factor of order one for
$|h|$ below
\begin{equation}
 w=\varepsilon\,\frac{|\nu_h|}{\lambda}\,\max_{\theta_0}|D_\gamma(\theta_0)| .
 \label{eq:width}
\end{equation}
This estimate of the half-width of the layer is heuristic, like the Chirikov
criterion \cite{Chirikov1979} on which it rests.

\begin{remark}[Polar-nodal normal forms]
\label{rem:pn}
For the family \eqref{eq:pn-H0} at fixed $N$, the angular system with
Hamiltonian $G$ can be written in action--angle variables $(\psi,J)$, and
$H_0$ then belongs to the class \eqref{eq:rotor} with $(q,p)=(r,R)$,
$(\theta,I)=(\psi,J)$ and $\nu=\frac{G'(J)}{2r^2}$. For $V=\rho(r)\,g$, formula
\eqref{eq:rotor-identity} gives
\[
 \Lop_{H_0}(V)=-\frac{G'(J)\,(r^2\rho)'\,R}{2r^4}\,\frac{\partial g}{\partial\psi}
 =\frac{(r^2\rho)'}{2r^4}\,R\,\{G,g\}_{\theta\Theta},
\]
because $\{G,g\}_{\psi J}=-G'(J)\,\partial_\psi g$ and both the bracket and the
off-diagonal value are unchanged by the canonical change
$(\theta,\Theta)\mapsto(\psi,J)$ within the angular factor
(Proposition~\ref{prop:intrinsic}(c)). This is \eqref{eq:pn-linear}, so
Proposition~\ref{prop:polar-nodal}(c) is the special case of
Proposition~\ref{prop:rotor} in which the rotor is the angular motion, and
the factor that distinguishes the defect from the Melnikov integrand is the
weight $\nu^{-2}$ together with one angle derivative. The homoclinic orbits
relevant for normal forms belong to the angular system, and treating them
with Theorem~\ref{thm:melnikov} requires the radial motion to be written in
action--angle variables, which is done in a companion paper.
\end{remark}

\begin{figure}[ht]
\centering
\includegraphics[width=\textwidth]{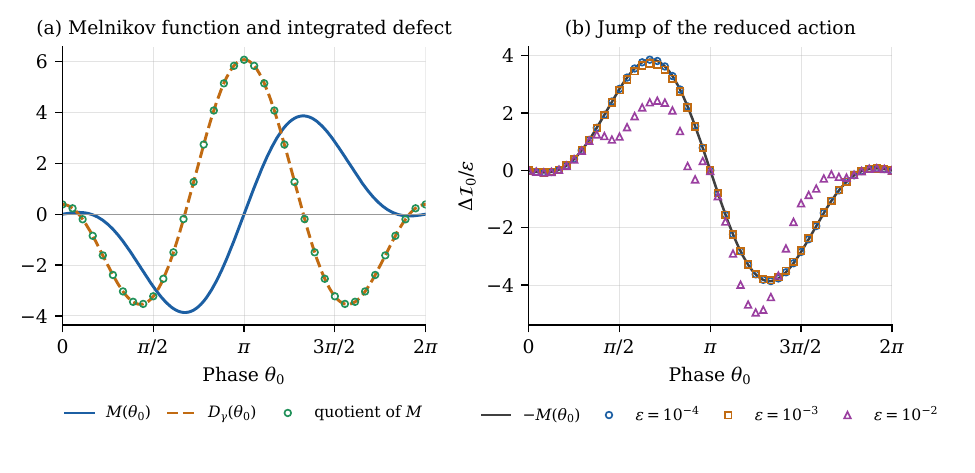}
\caption{Pendulum coupled to a rotor (Example~\ref{ex:rotor}). (a) Melnikov
function $M$ of \eqref{eq:M}, its numerical derivative and the integrated
defect $D_\gamma$ of \eqref{eq:Dg}, as functions of the phase $\theta_0$ of the
homoclinic solution; the last two curves coincide, as
Theorem~\ref{thm:melnikov}(b) states. (b) Change of the reduced action over
one passage along the separatrix, divided by $\varepsilon$, for three values of
$\varepsilon$, compared with the first-order prediction $-M$ of
\eqref{eq:jump}, which is accurate at $\varepsilon=10^{-4}$ and $10^{-3}$ and
fails at $\varepsilon=10^{-2}$.}
\label{fig:melnikov}
\end{figure}

\begin{example}[Pendulum coupled to a rotor]
\label{ex:rotor}
The Hamiltonian
\[
 H_0=\omega_0I+\tfrac12p^2+(1+kI)\cos q,\qquad
 V=(1-\cos q)\bigl(\cos\theta+\tfrac12\sin q\,\sin2\theta\bigr),
\]
with $\omega_0=1$ and $k=0.3$, has the rotor frequency $\nu=\omega_0+k\cos q$,
which lies in $[0.7,1.3]$, and for fixed $I$ it is a pendulum with hyperbolic
equilibrium $q=0$ and amplitude $A=1+kI$. At $I_h=0.5$ the homoclinic orbit is
\[
 q(t)=4\arctan e^{\sqrt A\,t},\qquad
 p(t)=2\sqrt A\,\operatorname{sech}\bigl(\sqrt A\,t\bigr),
\]
and $V$ vanishes to second order at $q=0$, so the hyperbolic periodic orbit is
not displaced by the perturbation. The identity \eqref{eq:rotor-identity}
holds symbolically with zero residual. On a grid of $145$ phases, the mean of
$M$ vanishes to rounding error, and the three-point difference quotient of $M$
agrees with $D_\gamma$ to $5.0\times10^{-3}$, against $\max|D_\gamma|=6.08$, the
difference being the error of the difference quotient. Integrations of $H$ at
$49$ phases, from the unperturbed separatrix at $t=-14\,A^{-\frac12}$ to
$t=7\,A^{-\frac12}$ on the perturbed energy surface, give jumps of $\mathcal
I_0$ whose deviation from \eqref{eq:jump}, in units of $\varepsilon\max|M|$ with
$\max|M|=3.87$, has the median $3.5\times10^{-5}$ and the maximum
$3.4\times10^{-4}$ at $\varepsilon=10^{-4}$, and the median $1.8\times10^{-3}$ and the
maximum $3.5\times10^{-2}$ at $\varepsilon=10^{-3}$, whereas at $\varepsilon=10^{-2}$
the median deviation of $0.17$ shows that the first-order regime has been
left (Figure~\ref{fig:melnikov}).
\end{example}

\subsection{Level sets and Kerr geodesics}
\label{sec:levelset}

By Proposition~\ref{prop:algebra}(f), a reparametrisation $\chi(H)$ of the
energy multiplies the five-input operator, and in particular every
off-diagonal value $b_{ij}$, by $\chi'(H)^3$. Multiplication of the Hamiltonian
by a function of the phase-space point, restricted to one energy surface,
changes the time parametrisation of the flow on that surface, as in the
Jacobi--Maupertuis principle, in the Mino time of Kerr geodesics
\cite{Mino2003} and in the choice of lapse for a Hamiltonian constraint. The
following lemma, which is introduced here, gives the behaviour of the
five-input operator under this operation.

\begin{lemma}[Level-set covariance]
\label{lem:levelset}
Let $K$ and $\sigma$ be smooth functions on an open subset of a symplectic
product. At every point of the level set $\{K=0\}$ and for all $i\ne j$,
\begin{equation}
 b_{ij}[\sigma K]=\sigma^3\,b_{ij}[K].
 \label{eq:levelset}
\end{equation}
\end{lemma}

\begin{proof}
Put $F=\sigma K$. On $\{K=0\}$ the identity $\dd F=\sigma\,\dd K+K\,\dd\sigma$
reduces to $\dd F=\sigma\,\dd K$, so the sector fields satisfy
$\Y_i^F=\sigma\Y_i^K$. The second derivatives are
\[
 \partial_a\partial_bF=\sigma\,\partial_a\partial_bK+\partial_a\sigma\,\partial_bK
 +\partial_aK\,\partial_b\sigma+K\,\partial_a\partial_b\sigma,
\]
and the last term vanishes on $\{K=0\}$. In the mixed block, with $a$ running
over the variables of the $i$th factor and $b$ over those of the $j$th, the
second term contracts on $(\Y_i^K,\Y_j^K)$ to $(\Y_i^K\sigma)\,\dd_jK(\Y_j^K)$,
where $\dd_jK$ is the part of $\dd K$ in the $j$th factor and
\[
 \dd_jK(\Y_j^K)=K_{q_j}K_{p_j}-K_{p_j}K_{q_j}=0,
\]
and the third term vanishes in the same way with $i$ and $j$ exchanged.
Hence the mixed Hessian of $F$ on $(\Y_i^F,\Y_j^F)$ is $\sigma^3$ times the
mixed Hessian of $K$ on $(\Y_i^K,\Y_j^K)$, and
Proposition~\ref{prop:intrinsic}(a) gives \eqref{eq:levelset}.
\end{proof}

For a constant $\sigma=t$, \eqref{eq:levelset} is the homogeneity of
Proposition~\ref{prop:algebra}(e), which holds at every point, and for a
nonconstant $\sigma$ the restriction to $\{K=0\}$ is needed because the terms
of the first and second derivatives of $\sigma K$ that contain $K$ as a factor
vanish only there.

\begin{corollary}[Separation on energy surfaces]
\label{cor:levelset}
Let $H$ be smooth, and suppose that for every constant $E$ in an interval
there is a smooth positive function $\sigma_E$ such that
$b_{ij}[\sigma_E(H-E)]=0$ on $\{H=E\}$. Then $b_{ij}[H]=0$ at every point at
which $H$ takes a value in that interval.
\end{corollary}

\begin{proof}
At a point $z$ put $E=H(z)$. Since $b_{ij}$ depends only on derivatives of its
argument, $b_{ij}[H](z)=b_{ij}[H-E](z)$, and Lemma~\ref{lem:levelset} with
$K=H-E$ gives $b_{ij}[H-E](z)=\sigma_E(z)^{-3}\,b_{ij}[\sigma_E(H-E)](z)=0$.
\end{proof}

Separation on an energy surface therefore does not depend on the time
parametrisation, and by \eqref{eq:levelset} a change of parametrisation
multiplies the integrand of the defect functional \eqref{eq:defect} on the
surface by $\sigma^6$, in analogy with Proposition~\ref{prop:defect}(c).

Geodesic motion in the Kerr metric \cite{Kerr1963} of mass $M_\bullet$ and spin
parameter $a$ is generated by $\frac12g^{\alpha\beta}p_\alpha p_\beta$, and in
Boyer--Lindquist coordinates \cite{BoyerLindquist1967} the conserved momenta $p_t=-E$ and $p_\phi=L$ reduce it to the factors
$(r,p_r)$ and $(\theta,p_\theta)$, with the Hamiltonian
\begin{gather}
 H_K=\frac{A_r+A_\theta}{2\Sigma},\qquad \Sigma=r^2+a^2\cos^2\theta,\qquad
 \Delta=r^2-2M_\bullet r+a^2,\label{eq:kerr-H}\\
 A_r=\Delta p_r^2-\frac{\bigl((r^2+a^2)E-aL\bigr)^2}{\Delta},\qquad
 A_\theta=p_\theta^2+\frac{(L-aE\sin^2\theta)^2}{\sin^2\theta}.\notag
\end{gather}
On the mass shell $H_K=-\frac12m^2$ of a particle of mass $m$, the functions
\[
 K_r=A_r+m^2r^2,\qquad K_\theta=A_\theta+m^2a^2\cos^2\theta
\]
satisfy $2\Sigma\,(H_K+\frac12m^2)=K_r+K_\theta$, which is the separation of
Carter \cite{Carter1968}. The sector fields $\Y_r$ and $\Y_\theta$ of the
sum-separable function $K_r+K_\theta$ are the Hamiltonian fields of $K_r$ and
$K_\theta$ in their own factors, and they commute.

\begin{proposition}[Kerr geodesics]
\label{prop:kerr}
On the domain $\Delta>0$, $0<\theta<\pi$ the following hold.
\begin{enumerate}[label=(\alph*)]
 \item $b_{r\theta}[H_K]\equiv0$.
 \item For every smooth deformation $h(r,\theta,p_r,p_\theta)$,
 \begin{equation}
  \Lop_{H_K}(h)_{r\theta}=-\frac{1}{(2\Sigma)^3}\,\Y_r\Y_\theta\bigl(2\Sigma h\bigr),
  \label{eq:kerr-L}
 \end{equation}
 where $m^2$ in $K_r$ and $K_\theta$ is set to $-2H_K$ at the point of
 evaluation.
 \item If $\Sigma h=f(r,p_r)+g(\theta,p_\theta)$, then $\Lop_{H_K}(h)_{r\theta}=0$.
\end{enumerate}
\end{proposition}

\begin{proof}
\emph{(a)} The function $K_r+K_\theta$ is sum-separable, so its off-diagonal
value vanishes identically for every $m^2$, and
Corollary~\ref{cor:levelset} with $E=-\frac12m^2$ and $\sigma_E=2\Sigma$ gives
the claim.

\emph{(b)} Fix a point $z$, put $H_\varepsilon=H_K+\varepsilon h$ and
$m_\varepsilon^2=-2H_\varepsilon(z)$, so that $z$ lies on the zero set of
$H_\varepsilon+\frac12m_\varepsilon^2$. Lemma~\ref{lem:levelset} with
$\sigma=2\Sigma$ gives
$b_{r\theta}[H_\varepsilon](z)=(2\Sigma)^{-3}\,b_{r\theta}[K_r+K_\theta+2\varepsilon\Sigma h](z)$,
where $K_r$ and $K_\theta$ carry $m_\varepsilon^2$. For a sum-separable function the
mixed Hessian vanishes and the sector fields act on separate variables, so
the first-order coefficient of the off-diagonal value of
$K_r+K_\theta+\varepsilon G$ is $-\Y_r\Y_\theta G$. The dependence of
$m_\varepsilon^2$ on $\varepsilon$ changes the sector fields only in terms that are
already of order $\varepsilon$ and therefore contributes at order
$\varepsilon^2$, and differentiation at $\varepsilon=0$ gives \eqref{eq:kerr-L}.

\emph{(c)} $\Y_\theta f=0$, and $\Y_rg=0$.
\end{proof}

The direct evaluation of \eqref{eq:LC} for \eqref{eq:kerr-H} confirms (a),
and \eqref{eq:kerr-L} was checked for
$h=\frac12(\alpha p_r^2+\beta p_\theta^2+\gamma E^2)$ with the three non-separable
functions
\[
\alpha=1+\frac{r\sin\theta}{10},\qquad
\beta=\exp\Bigl(-\frac{r\cos\theta}{5}\Bigr),\qquad
\gamma=\frac{r^2\theta}{1+r\theta}
\]
at $20$ random points, where the two sides of \eqref{eq:kerr-L} differ by at most
$1.7\times10^{-49}$ in 50-digit arithmetic and by a relative residual of at most
$5.6\times10^{-16}$ in double precision. Item (c) recovers, at first order, the condition used to construct deformed
black-hole metrics that keep a Carter-like constant
\cite{VigelandYunesStein2011,Johannsen2013,PapadopoulosKokkotas2018}, namely
the additive separability of the $\Sigma$-weighted deformed Hamiltonian. For
any other deformation, \eqref{eq:kerr-L} gives, pointwise and in closed form,
the first-order obstruction to separation in the Boyer--Lindquist product
structure. A nonzero value excludes a Carter-like constant obtained by
separation in these coordinates, but it does not exclude integrability
through a different splitting or through a deformed integral of another
form, since even the integrable oscillator of Example~\ref{ex:oscillator}
has a nonzero off-diagonal value in its original pairs. A dynamical
conclusion follows only from the integrated defect, through
Theorem~\ref{thm:melnikov} and Corollary~\ref{cor:dichotomy}. In Mino time
the mass-shell Hamiltonian $K_r+K_\theta$ is sum-separable, so after
action--angle variables are introduced for the polar motion it belongs to
the class \eqref{eq:rotor}, and its radial separatrices, the zoom--whirl
orbits homoclinic to the unstable spherical orbits, are known in closed form
\cite{LevinPerezGiz2009}. Lemma~\ref{lem:levelset} transfers the defect
\eqref{eq:kerr-L} to Mino time, so Theorem~\ref{thm:melnikov} applies to
deformed Kerr geodesics near these separatrices.

\section{Exact coordinate identification}
\label{sec:euler}

This section investigates the defect as a coordinate-selection tool in the
simplest setting, a family of product structures that contains a separating
member, and the results show that minimising the defect recovers that member
exactly. A pair of coupled linear oscillators already displays the mechanism,
and the two-centre problem of Euler, which separates in elliptic coordinates
\cite{WaalkensDullinRichter}, provides a nonlinear case in which the
separating focal distance is known exactly and the minimisation can be
carried out in closed form (Theorem~\ref{thm:euler}).

\subsection{Coupled linear oscillators}

In the linear case, the family of product structures is generated by
rotations of the configuration plane, and the defect functional is an
explicit function of the rotation angle.

\begin{example}[Rotated oscillator]
	\label{ex:oscillator}
	Consider the Hamiltonian
	\[
	H=\tfrac12(p_1^2+p_2^2+q_1^2+q_2^2)+\epsilon\,q_1q_2,\qquad 0<|\epsilon|<1,
	\]
	together with the canonical rotations
	\[
	\mathbf Q=\mathbf R_\varphi\mathbf q,\qquad
	\mathbf P=\mathbf R_\varphi\mathbf p,\qquad
	\mathbf R_\varphi=\begin{pmatrix}\cos\varphi&\sin\varphi\\-\sin\varphi&\cos\varphi\end{pmatrix}.
	\]
	In the rotated pairs the coupling coefficient becomes $\epsilon\cos2\varphi$,
	so that \eqref{eq:LC} gives
	$b^{(\varphi)}_{12}=-\epsilon\cos2\varphi\,P_1P_2$. On a ball $U$ centred at
	the origin and with unit weight, the rotation invariance of $U$ and of
	$\dd\mu$ yields
	\[
	\mathcal E(\varphi)=C_U\,\epsilon^2\cos^22\varphi,\qquad
	C_U=\int_UP_1^2P_2^2\,\dd\mu>0,
	\]
	where $C_U$ does not depend on $\varphi$. The minima $\varphi=\pi/4+m\pi/2$ are therefore exactly the normal-mode
	pairs. Because the rotations mix the original factors, the minimisation acts
	on the product structure itself and selects the splitting into normal modes
	among all rotated splittings.
\end{example}

\subsection{The two-centre problem}

The two-centre problem of Euler is the main object of this section, and it
also serves as the nonlinear model for the galactic application of
Section~\ref{sec:stackel}, where the focal distance of prolate coordinates
has to be chosen in the same way. A particle of unit mass attracted by two
centres at $(\pm c,0)$, $c>0$, has the Hamiltonian \cite{WaalkensDullinRichter}
\begin{equation}
	H=\tfrac12(p_x^2+p_y^2)+V,\qquad
	V=-\frac{\mu_1}{r_+}-\frac{\mu_2}{r_-},\qquad
	r_\pm=\sqrt{(x\mp c)^2+y^2},
	\label{eq:euler-H}
\end{equation}
with $\mu_1,\mu_2>0$, and it separates in the elliptic coordinates whose
foci coincide with the two centres. To test the defect, the focal distance is
treated as unknown: for a trial focal distance $a>0$, the elliptic
coordinates
\[
x=a\cosh\xi\cos\eta,\qquad y=a\sinh\xi\sin\eta
\]
and their cotangent lift give
\begin{equation}
	H=\frac{p_\xi^2+p_\eta^2}{2D_a}+V,\qquad
	D_a=a^2(\sinh^2\xi+\sin^2\eta),\qquad
	b_a:=\B_{\xi,p_\xi;\eta,p_\eta}[H].
\end{equation}
The physical Hamiltonian is the same for every $a$, and the parameter only
selects the product structure.

\begin{proposition}[Exact two-centre defect]
	\label{prop:euler-defect}
	Let $V_{xy}=\partial^2V/\partial x\,\partial y$ denote the mixed Cartesian
	derivative of the two-centre potential \eqref{eq:euler-H}. Then the
	off-diagonal value of the elliptic structure with focal distance $a$ is
	\begin{equation}
		b_a=-\frac{p_\xi p_\eta}{D_a^3}\,\frac{\partial^2(D_aV)}{\partial\xi\,\partial\eta}
		=(a^2-c^2)\,\frac{V_{xy}\,Q_a}{W_a},
		\label{eq:euler-defect}
	\end{equation}
	where $W_a$ and $Q_a$ are the following functions of the Cartesian
	phase-space point:
	\begin{equation}
		W_a=D_a^2=\bigl((x-a)^2+y^2\bigr)\bigl((x+a)^2+y^2\bigr),\qquad
		Q_a=p_\xi p_\eta=(x^2-y^2-a^2)\,p_xp_y+xy\,(p_y^2-p_x^2).
		\label{eq:euler-cartesian}
	\end{equation}
\end{proposition}

\begin{proof}
	Since $\partial_{\xi\eta}D_a=0$, the purely kinetic terms cancel in
	\eqref{eq:LC}, which gives the first equality. For the second, let
	$r_e=\sqrt{(x-e)^2+y^2}$ be the distance to the point $(e,0)$. The coordinate
	fields are
	\[
	\partial_\xi=a\sinh\xi\cos\eta\,\partial_x+a\cosh\xi\sin\eta\,\partial_y,\qquad
	\partial_\eta=-a\cosh\xi\sin\eta\,\partial_x+a\sinh\xi\cos\eta\,\partial_y,
	\]
	and the chain rule yields
	\[
	\frac{\partial^2}{\partial\xi\,\partial\eta}\Bigl(\frac{D_a}{r_e}\Bigr)
	=-3D_a(a^2-e^2)\,\frac{y\,(x-e)}{r_e^5}
	=-D_a(a^2-e^2)\,\frac{\partial^2}{\partial x\,\partial y}\Bigl(\frac{1}{r_e}\Bigr).
	\]
	The two terms of $V$ correspond to $e=c$ and $e=-c$, for both of which
	$a^2-e^2=a^2-c^2$, so that
	\[
	\frac{\partial^2(D_aV)}{\partial\xi\,\partial\eta}=-D_a(a^2-c^2)\,V_{xy},
	\]
	and the Cartesian expressions \eqref{eq:euler-cartesian} follow from the
	coordinate map and its cotangent lift.
\end{proof}

The factor $a^2-c^2$ shows at once that $b_c\equiv0$, in agreement with the
classical separation of the problem: at $a=c$ the product
\[
D_cV=-c\bigl[(\mu_1+\mu_2)\cosh\xi+(\mu_1-\mu_2)\cos\eta\bigr]
\]
is the sum of a function of $\xi$ and a function of $\eta$. With the local
focal estimate $\delta_{\mathrm{loc}}$ introduced in
Section~\ref{sec:stackel}, \eqref{eq:euler-defect} states that the
two-centre potential has $\delta_{\mathrm{loc}}=c^2$
at every point where $V_{xy}\ne0$. The Cartesian form \eqref{eq:euler-defect}
is smooth away from the centres and the trial foci $(\pm a,0)$, and it agrees
with the canonical evaluation in every local elliptic chart: the two branches
related by $(\xi,\eta)\mapsto(-\xi,-\eta)$ reverse both momenta and therefore
leave $b_a$ unchanged. No global elliptic chart is thus needed, and the only
requirement on the integration domain in the next theorem is that its closure
be compact and avoid the centres and the trial foci.

\begin{theorem}[Nonlinear coordinate identification]
	\label{thm:euler}
	Let $A=[a_-,a_+]$ with $0<a_-<c<a_+$, let
	$\Gamma_A=\{(\pm a,0):a\in A\}$, and let $\Omega$ be a nonempty open domain
	whose closure is compact and disjoint from $\Gamma_A$. Let $w>0$ be smooth
	near $\overline\Omega$, let $P_0>0$, and put
	\[
	\mathcal E_{\mathrm E}(a)
	=\int_\Omega\int_{p_x^2+p_y^2<P_0^2}w\,b_a^2\,\dd p_x\,\dd p_y\,\dd x\,\dd y .
	\]
	Then the following hold.
	\begin{enumerate}[label=(\alph*)]
		\item \emph{Factorisation.} For every $a\in A$,
		\begin{equation}
			\mathcal E_{\mathrm E}(a)=\frac{\pi P_0^6}{24}\,(a^2-c^2)^2\,I(a),\qquad
			I(a)=\int_\Omega\frac{w\,V_{xy}^2}{W_a}\,\dd x\,\dd y>0.
			\label{eq:euler-objective}
		\end{equation}
		\item \emph{Unique minimiser.} The separating value $a=c$ is the unique
		global minimiser of $\mathcal E_{\mathrm E}$ on $A$.
		\item \emph{Nondegeneracy.} $\mathcal E_{\mathrm E}$ is smooth on $A$, and
		\[
		\mathcal E_{\mathrm E}''(c)=\frac{\pi P_0^6}{3}\,c^2\,I(c)>0.
		\]
	\end{enumerate}
\end{theorem}

\begin{proof}
	\emph{(a)} Write $Q_a=\mathcal A\,p_xp_y+\mathcal B\,(p_y^2-p_x^2)$ with
	$\mathcal A=x^2-y^2-a^2$ and $\mathcal B=xy$, so that
	$\mathcal A^2+4\mathcal B^2=W_a$. Over the momentum disk,
	\[
	\int p_x^2p_y^2=\frac{\pi P_0^6}{24},\qquad
	\int(p_y^2-p_x^2)^2=\frac{\pi P_0^6}{6},
	\]
	while the cross moment vanishes by reflection symmetry. Hence
	$\int Q_a^2=\pi P_0^6W_a/24$, and substituting \eqref{eq:euler-defect} gives
	the factorisation. Direct differentiation gives
	\[
	V_{xy}=-3y\Bigl(\frac{\mu_1(x-c)}{r_+^5}+\frac{\mu_2(x+c)}{r_-^5}\Bigr),
	\]
	which is negative for $x>c$, $y>0$. Since $V_{xy}$ is real-analytic on the
	connected set $\mathbb R^2\setminus\{(\pm c,0)\}$, it cannot vanish
	identically on any nonempty open set, and together with $w>0$ and $W_a>0$
	this gives $I(a)>0$.
	
	\emph{(b)} By (a), $\mathcal E_{\mathrm E}\ge0$ on $A$, with equality exactly
	when $a^2=c^2$, that is, at $a=c$.
	
	\emph{(c)} The positive distance between $\overline\Omega$ and $\Gamma_A$
	bounds $W_a$ away from zero uniformly in $a\in A$, so $I$ is smooth on $A$.
	Differentiating \eqref{eq:euler-objective} twice, all terms that still
	contain the factor $a^2-c^2$ vanish at $a=c$, and the remaining term is
	$\tfrac{\pi P_0^6}{24}\cdot8c^2I(c)$.
\end{proof}

\begin{figure}[t]
\centering
\includegraphics[width=\textwidth]{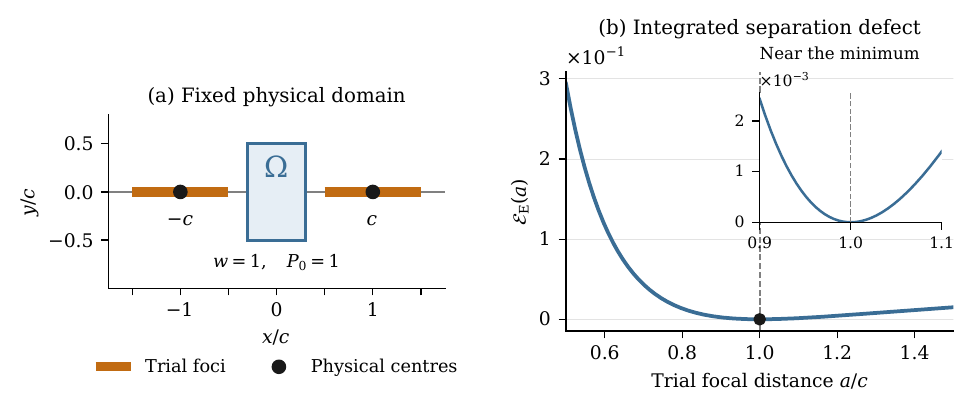}
\caption{Euler coordinate identification on the fixed axis-crossing domain
$\Omega=(-0.3,0.3)\times(-0.5,0.5)$ with $c=1$, $\mu_1=1$, $\mu_2=0.7$,
$P_0=1$ and $w=1$. (a) The domain avoids the whole focal set of the trial
interval $[0.5,1.5]$. (b) The exact objective \eqref{eq:euler-objective},
with an inset that resolves the nondegenerate minimum at $a=c$.}
\label{fig:euler}
\end{figure}

The domain may cross either axis, provided that it avoids the whole trial
focal set, and Figure~\ref{fig:euler} shows such an axis-crossing example.
The identification also survives energy reparametrisation, since by
\eqref{eq:energychange} replacing $H$ by $\chi(H)$ with $\chi'\ne0$
multiplies the integrand by $\chi'(H)^6$; the momentum moments then no longer
factor out, but the prefactor $(a^2-c^2)^2$ remains, so that $a=c$ is still
the unique minimiser. The same computation finally covers the spatial
two-centre problem with the centres on the $z$-axis. After axial reduction,
the prolate coordinates $R=a\sinh\xi\sin\eta$, $z=a\cosh\xi\cos\eta$ give
the same defect with $(x,y)$ replaced by $(z,R)$, $R>0$, because the
centrifugal term contributes
\[
D_a\,\frac{L_z^2}{2R^2}
=\frac{L_z^2}{2}\Bigl(\frac{1}{\sinh^2\xi}+\frac{1}{\sin^2\eta}\Bigr)
\]
to $D_aH$, and this sum has no mixed derivative.

\section{St\"ackel residuals and focal-distance selection}
\label{sec:stackel}

Unlike the two-centre problem, a general axisymmetric galactic potential has
no exactly separating member in the family of prolate product structures. At
fixed axial angular momentum $L_z$, the meridional Hamiltonian and the
prolate coordinates of focal distance $\Delta$ are
\begin{equation}
	H=\tfrac12(p_R^2+p_z^2)+\frac{L_z^2}{2R^2}+\Phi(R,z),\qquad
	R=\Delta\sinh u\sin v,\quad z=\Delta\cosh u\cos v,
	\label{eq:fudge-H}
\end{equation}
with $u>0$, $0<v<\pi$, $D=\Delta^2(\sinh^2u+\sin^2v)$ and cotangent-lifted
momenta $p_u,p_v$; in this section $u$ and $v$ denote these prolate
coordinates and are unrelated to the response rates of
Section~\ref{sec:operator}. The potential is said to be of
\emph{St\"ackel form} in these coordinates when
\[
\Phi(R,z)=\frac{\mathcal U(u)+\mathcal V(v)}{\sinh^2u+\sin^2v}
\]
for two functions $\mathcal U$ and $\mathcal V$ of one variable, that is,
when $D\Phi$ is the sum of a function of $u$ and a function of $v$
\cite{Stackel,DeZeeuw1985}. Since the centrifugal term contributes such a
sum to $DH$ as well, the Hamilton--Jacobi equation then separates in
$(u,v)$, and the radial and vertical actions reduce to one-dimensional
quadratures. The St\"ackel fudge \cite{Binney2012,SandersBinney2016}
evaluates these quadratures at a single phase-space point of a general
axisymmetric potential, as if the potential were of St\"ackel form for a
suitably chosen $\Delta$, while $L_z$ is conserved exactly. The accuracy of
the resulting estimates depends on the focal distance, and the choice of
$\Delta$ is the subject of this section. Proposition~\ref{prop:fudge} below
shows that the separability operator measures exactly the departure from
St\"ackel form, through the residual $C_\Delta$
\begin{equation}
	g=\Phi_{Rz},\qquad
	A=3z\Phi_R-3R\Phi_z+Rz(\Phi_{RR}-\Phi_{zz})+(z^2-R^2)g,\qquad
	C_\Delta=A-g\Delta^2,
	\label{eq:residual}
\end{equation}
which is built from the first and second derivatives of the potential.

\begin{proposition}[Operator--residual factorisation]
	\label{prop:fudge}
	Let $b_\Delta=\B_{u,p_u;v,p_v}[H]$ denote the off-diagonal value of the
	prolate structure with focal distance $\Delta$ for the Hamiltonian
	\eqref{eq:fudge-H}. Then the following hold.
	\begin{enumerate}[label=(\alph*)]
		\item \emph{Factorisation.} The off-diagonal value is proportional to
		the residual \eqref{eq:residual},
		\begin{equation}
			b_\Delta=-\frac{p_up_v}{D^3}\,\frac{\partial^2(D\Phi)}{\partial u\,\partial v}
			=-\frac{p_up_v}{D^2}\,C_\Delta.
			\label{eq:fudge-B}
		\end{equation}
		\item \emph{Weighted square.} The square of the off-diagonal value is the
		square of the residual with a nonnegative weight,
		\begin{equation}
			b_\Delta^2=w_\Delta\,C_\Delta^2,\qquad w_\Delta=\frac{p_u^2p_v^2}{D^4}.
			\label{eq:fudge-weight}
		\end{equation}
	\end{enumerate}
\end{proposition}

\begin{proof}
	\emph{(a)} Since $\partial_{uv}D=0$, the purely kinetic terms cancel in
	\eqref{eq:LC}, as in the proof of Proposition~\ref{prop:euler-defect}. The
	centrifugal term contributes
	\[
	D\,\frac{L_z^2}{2R^2}
	=\frac{L_z^2}{2}\Bigl(\frac{1}{\sinh^2u}+\frac{1}{\sin^2v}\Bigr)
	\]
	to $DH$, which is the sum of a function of $u$ and a function of $v$ and has
	no mixed derivative, and this gives the first equality. For the second equality in \eqref{eq:fudge-B}, which expresses the mixed derivative through the residual, the
	coordinate fields are
	\[
	\partial_u=\Delta\cosh u\sin v\,\partial_R+\Delta\sinh u\cos v\,\partial_z,\qquad
	\partial_v=\Delta\sinh u\cos v\,\partial_R-\Delta\cosh u\sin v\,\partial_z,
	\]
	and the chain rule yields
	\[
	\frac{\partial^2(D\Phi)}{\partial u\,\partial v}=D\,C_\Delta .
	\]
	
	\emph{(b)} Squaring \eqref{eq:fudge-B} gives \eqref{eq:fudge-weight}.
\end{proof}

Proposition~\ref{prop:fudge} gives the five-input operator, evaluated on the
prolate probes $(u,p_u;v,p_v)$, a direct meaning for the St\"ackel fudge.
Where $g\ne0$, the local focal estimate $\delta_{\mathrm{loc}}=A/g$ turns the
residual into $C_\Delta=g(\delta_{\mathrm{loc}}-\Delta^2)$, so that
\[
b_\Delta=-\frac{p_up_v\,g}{D^2}\,\bigl(\delta_{\mathrm{loc}}-\Delta^2\bigr).
\]
At a point with $p_up_vg\ne0$ the operator therefore vanishes exactly when
$\Delta^2=\delta_{\mathrm{loc}}$, while zeros caused by $p_u=0$ or $p_v=0$
impose no condition on $\Delta$. The quantity $\delta_{\mathrm{loc}}$ is
precisely the local focal prescription of Sanders \cite{Sanders2012}, with
the sign correction given in \cite[Eq.~(15) and footnote~1]{SandersBinney2016},
so this established prescription is recovered as the pointwise zero of the
separability operator.

A single focal distance must, however, serve a whole orbit, and the operator
offers three natural ways of combining the local information along it, which
differ in how the momentum weight $w_\Delta$ of \eqref{eq:fudge-weight} is
treated. Let $\mathbf x_m$, $m=1,\ldots,N$, be \emph{training points}, that is, phase-space points sampled along the integrated orbit as described in Section~\ref{sec:numerics}, at which the selectors \eqref{eq:sel-med}--\eqref{eq:sel-C} below evaluate the potential derivatives and the momenta that determine their focal distance. They are the same physical points
for every trial focal distance, and they are kept separate from the later
\emph{validation points} at which the action spreads are evaluated
(Section~\ref{sec:numerics}). Write $A_m$, $g_m$ and $\delta_{\mathrm{loc},m}$ for the values at
$\mathbf x_m$ of the residual at zero focal distance $A=C_0$, the mixed
derivative $g=\Phi_{Rz}$ of \eqref{eq:residual} and the local focal estimate
$\delta_{\mathrm{loc}}=A/g$, respectively, and let
$I=[\delta_-,\delta_+]$ be a search interval for $\delta=\Delta^2$, with
$\proj_I$ the projection onto it, which sends a number to the nearest point
of $I$. The three selectors built on the orbit itself are
\begin{align}
	\delta_{\mathrm{med}}&=\proj_I\bigl(\operatorname{median}_m\delta_{\mathrm{loc},m}\bigr),
	\label{eq:sel-med}\\
	\delta_B&\in\argmin_{\delta\in I}F_B(\delta),\qquad
	F_B(\delta)=\frac1N\sum_{m=1}^Nb_{\sqrt\delta}(\mathbf x_m)^2,
	\label{eq:sel-B}\\
	\delta_C&=\argmin_{\delta\in I}F_C(\delta)
	=\proj_I\Bigl(\frac{\sum_mg_mA_m}{\sum_mg_m^2}\Bigr),\qquad
	F_C(\delta)=\frac1N\sum_{m=1}^N(A_m-g_m\delta)^2.
	\label{eq:sel-C}
\end{align}
In terms of the operator, $\delta_{\mathrm{med}}$ is the median of its
pointwise zeros, taken over the samples at which $\delta_{\mathrm{loc},m}$
is defined or over a continuous model-specific extension, $F_B$ is the
sampled analogue of the defect functional of Definition~\ref{def:defect},
with the training points in place of the symplectic volume, and $F_C$ is the
sampled square of the operator after division by its momentum weight, by
\eqref{eq:fudge-weight}. Since both the residual and the weight depend on
$\Delta$, $F_B$ is minimised numerically. The quantities $A_m$ and $g_m$ do
not depend on $\Delta$, so $F_C$ is a quadratic polynomial in $\delta$,
strictly convex when $\sum_mg_m^2>0$, and its unconstrained minimiser
\[
\frac{\sum_mg_mA_m}{\sum_mg_m^2}
=\frac{\sum_mg_m^2\,\delta_{\mathrm{loc},m}}{\sum_mg_m^2}
\]
is the $g^2$-weighted mean of the local estimates over the samples with
$g_m\ne0$. Minimising the deviation of $C_\Delta$ from zero over the region
occupied by an orbit was proposed by Vasiliev
\cite[Section~A.5.1, Eq.~(41b)]{VasilievDoc}, and \eqref{eq:sel-C} is a
closed-form implementation of this idea on a fixed sample. The selectors differ on two degenerate sets as follows. On the set $g=\Phi_{Rz}=0$,
which contains the equatorial plane of every potential symmetric about it,
the local estimate is undefined, whereas $C_\Delta$ and $b_\Delta$ remain
smooth. At the turning points of either prolate coordinate, where
$p_up_v=0$ and the St\"ackel quadratures have their endpoints, the
hypothesis $H_{p_i}\ne0$ of the momentum-graph criterion of
Proposition~\ref{prop:LC} fails, and the operator objective $F_B$ suppresses
these regions through the smooth weight $w_\Delta$, whereas Sanders'
prescription and Vasiliev's criterion depend only on the position and give
them their full sample share. Written as a defect functional, $F_C$ would
therefore carry the weight $w_\Delta^{-1}$, which is singular at the turning
points and lies outside the class of smooth weights of
Section~\ref{sec:defects}.

A fourth action-free selector, due to Binney \cite[Appendix~I]{Binney2014},
uses the shell orbit $\gamma_{\mathrm{sh}}$ of the same energy $E$ and axial
angular momentum $L_z$, that is, the orbit with $J_R=0$, which leaves the
equatorial plane at a radius $R_s$ with $p_R=0$ and returns to it with
$p_R=0$. In a St\"ackel potential such an orbit lies on a coordinate surface
$u=\mathrm{const}$, which in the meridional plane is the ellipse
\[
 \frac{R^2}{R_s^2}+\frac{z^2}{R_s^2+\Delta^2}=1,
\]
and fitting this equation to points $(R_k,z_k)$ of one passage of
$\gamma_{\mathrm{sh}}$ above the plane by linear least squares in the single
unknown $(R_s^2+\Delta^2)^{-1}$ gives
\begin{equation}
 \delta_{\mathrm{sh}}=\proj_I\biggl(\frac{\sum_kz_k^4}{\sum_kz_k^2\bigl(1-R_k^2R_s^{-2}\bigr)}-R_s^2\biggr).
 \label{eq:sel-sh}
\end{equation}
In terms of the operator, $p_u=D\dot u$ vanishes identically along
$\gamma_{\mathrm{sh}}$ exactly when the orbit lies on a surface
$u=\mathrm{const}$, and then $b_\Delta=0$ along $\gamma_{\mathrm{sh}}$ by
\eqref{eq:fudge-B}, whatever the residual. The shell-orbit fit thus makes
the operator vanish through its momentum factor on a neighbouring orbit of
the same energy, at the cost of one additional orbit integration for every
$(E,L_z)$, whereas $\delta_B$ makes it small on the orbit itself. For the
100 orbits of Section~\ref{sec:numerics}, the minimiser of $F_B$ on the
points of $\gamma_{\mathrm{sh}}$ agrees with $\Delta_{\mathrm{sh}}$ to within
$0.92\%$ of $R_0$.

Through Proposition~\ref{prop:fudge}, the four selectors are thus statements
about the same operator, namely its pointwise zero, its sampled square, its
sampled square with the momentum weight removed and its zero on the shell
orbit, and none of them uses action estimates. Since $F_B$ suppresses the
turning regions and $F_C$ does not, the two fits can select different focal
distances, and the consequences for the actions are examined in
Section~\ref{sec:weighting}.

\subsection{Numerical design}
\label{sec:numerics}

The four selectors are compared on two standard models, the flattened
logarithmic potential $\Phi_{\log}$ \cite{Binney1981}, \cite[Section~2.3.2]{BinneyTremaine2008} and the Miyamoto--Nagai potential
$\Phi_{\mathrm{MN}}$ \cite{MiyamotoNagai},
\begin{equation}
 \Phi_{\log}=\frac{v_0^2}{2}\log\Bigl(R^2+\frac{z^2}{q^2}+r_c^2\Bigr),\qquad
 \Phi_{\mathrm{MN}}=-\frac{\mu_{\mathrm{MN}}}{\sqrt{R^2+\bigl(a_d+\sqrt{z^2+b_d^2}\bigr)^2}},
 \label{eq:potentials}
\end{equation}
with $(v_0^2,q,r_c)=(1.01,0.9,0.1)$, $(a_d,b_d)=(1,0.2)$ and
$\mu_{\mathrm{MN}}=[1+(a_d+b_d)^2]^{3/2}$, chosen so that both have unit
circular speed at $R=1$. With $v_c(R)=\sqrt{R\Phi_R(R,0)}$ and
$\zeta=\sqrt{z^2+b_d^2}$, direct differentiation gives the continuous local
estimates
\begin{equation}
 \delta_{\mathrm{loc,log}}=(1-q^2)\Bigl(R^2+\frac{z^2}{q^2}+2r_c^2\Bigr),\qquad
 \delta_{\mathrm{loc,MN}}=a_d^2+\frac{a_db_d^2}{\zeta}+\frac{a_db_d^2\bigl(R^2+(a_d+\zeta)^2\bigr)}{3(a_d+\zeta)\zeta^2},
 \label{eq:dloc-models}
\end{equation}
both of which are positive and remain defined at $z=0$, where $A/g$ itself
is indeterminate.

The local estimate $\delta_{\mathrm{loc,MN}}$ in \eqref{eq:dloc-models} is
not constant, so by Proposition~\ref{prop:fudge} the defect $b_\Delta$
cannot vanish on any open set, whatever the focal distance, and no prolate
structure separates the Miyamoto--Nagai model. Its two limits are exactly
separable: for $a_d\to0$ it reduces to the spherical Plummer model, with
$\delta_{\mathrm{loc}}=0$, and for $b_d\to0$ to the Kuzmin disc
\cite{Kuzmin1956,Toomre1963}, with $\delta_{\mathrm{loc}}=a_d^2$ off the
plane. The second limit is not uniform, since for small $b_d$ and heights
$z$ comparable to $b_d$
\[
\delta_{\mathrm{loc,MN}}-a_d^2\simeq\frac{(R^2+a_d^2)\,b_d^2}{3(z^2+b_d^2)},
\]
so the departure from St\"ackel form is concentrated in a layer of thickness
comparable to $b_d$ about the plane, with an amplitude that does not
decrease with $b_d$ and grows with $R$. Along orbits with large vertical excursions at large radii, the local
estimate therefore varies by a factor between about two and four and a
half, so that a single focal distance can match it only on part of such an
orbit, and these orbits produce the largest action spreads in
Section~\ref{sec:results}.

The orbits are launched from the equatorial plane, $z_0=0$, at the radius
$R_0$, with axial angular momentum $L_z=0.8R_0v_c(R_0)$ and with radial and
vertical momenta $p_{R,0}=f_Rv_c(R_0)$ and $p_{z,0}=f_zv_c(R_0)$, given as
fractions $f_R$ and $f_z$ of the circular speed. For each potential, the
$2\times5\times5=50$ launch states are those of the grid
\[
R_0\in\{1,2\},\qquad
f_R\in\{0,0.15,0.30,0.45,0.60\},\qquad
f_z\in\{0.05,0.15,0.30,0.45,0.60\},
\]
which ranges from weakly to strongly excited motion in both directions at
two length scales. The grid was fixed before the selectors were compared
and does not sample a galactic distribution function. Each orbit is integrated with the eighth-order Runge--Kutta method DOP853 \cite{HairerNorsettWanner1993} of SciPy \cite{SciPy2020}, at relative and absolute
tolerances $3\times10^{-12}$ and $10^{-13}$, for 128 launch circular periods
$T_c=2\pi R_0/v_c(R_0)$, the first 64 supplying 4097 training points and the
last 64 supplying 512 validation points, and all selectors use the interval
$0.03R_0\le\Delta\le3R_0$, with $F_B$ minimised by a scan at 1001 points and
local refinement. The actions of all four selectors are computed by the same
St\"ackel quadrature, that of galpy~1.12.0 at Gauss--Legendre order 64
\cite{Bovy2015,galpy1120}, whose relative discretisation error, bounded by comparison
with order 96, is at most $1.47\times10^{-6}$, so that on each orbit the
four estimates differ only through the focal distance.

The true actions are constant along an orbit, so the accuracy of each
selector can be measured, without reference actions, by how much its action
estimates vary along the orbit, a criterion used to assess action
estimators in \cite{Binney2012,SandersBinney2016}. For a selector
$s\in\{\mathrm{med},B,C,\mathrm{sh}\}$ and an action $\ell\in\{R,z\}$, let
$\widehat J^{(s)}_{\ell,m}$, $m=1,\ldots,N_v=512$, be the St\"ackel-fudge
estimates at the validation points and $\bar J^{(s)}_\ell$ their mean. The
\emph{relative action spread} $s^{(s)}_\ell$ is the standard deviation of
these estimates divided by their mean, and the \emph{joint spread}
$s^{(s)}_{\mathrm{joint}}$ is the root mean square of the radial and
vertical spreads,
\begin{equation}
 s^{(s)}_\ell=\frac{1}{\bar J^{(s)}_\ell}
 \Bigl[\frac1{N_v}\sum_{m=1}^{N_v}\bigl(\widehat J^{(s)}_{\ell,m}-\bar J^{(s)}_\ell\bigr)^2\Bigr]^{1/2},
 \qquad
 s^{(s)}_{\mathrm{joint}}=\Bigl[\tfrac12\bigl((s^{(s)}_R)^2+(s^{(s)}_z)^2\bigr)\Bigr]^{1/2}.
 \label{eq:spread}
\end{equation}
Since the mean square deviation from a reference action is the sum of the
mean square deviation from the mean and the squared offset of the mean, the
spread misses only the bias of the estimates, which reference actions
obtained independently from a toy-torus generating-function fit
\cite[Section~3.2]{SandersBinney2016} show to be small on this grid
(Section~\ref{sec:results} and Appendix~\ref{app:reference}).

\subsection{Results}
\label{sec:results}

\begin{figure}[h]
\centering
\includegraphics[width=0.92\textwidth]{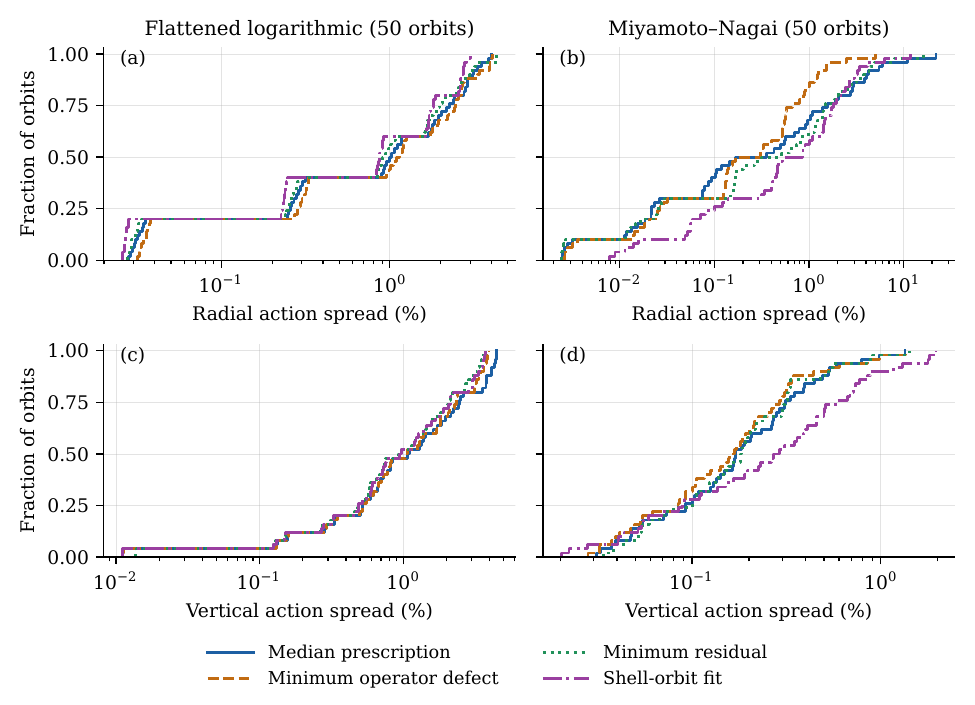}
\caption{Empirical cumulative distributions of the relative action spreads
	\eqref{eq:spread} over all 50 orbits of each potential, radial in the top
	row and vertical in the bottom row, for the logarithmic (left) and
	Miyamoto--Nagai (right) potentials. A curve gives the fraction of orbits
	whose spread is at or below the value on the horizontal axis, so that a
	curve further to the left corresponds to more accurate actions. The legend
	entries \emph{median prescription}, \emph{minimum operator defect},
	\emph{minimum residual} and \emph{shell-orbit fit} denote the selectors
	$\Delta_{\mathrm{med}}$, $\Delta_B$, $\Delta_C$ and $\Delta_{\mathrm{sh}}$,
	all evaluated on the same validation points with the same St\"ackel
	quadrature.}
\label{fig:cdf}
\end{figure}

All 100 integrations, shell-orbit searches and St\"ackel calculations
completed, with a largest energy drift, scaled by
$\max\{|H(0)|,v_c(R_0)^2\}$, of $1.18\times10^{-10}$, and none of the 400
focal choices lies on a boundary of the search interval. None of the four
selectors computes actions, so all belong to the class of inexpensive
prescriptions used to tabulate the focal distance as a function of the
integrals of motion \cite{Binney2014,VasilievAgama}, and the median
prescription, the simplest established choice in this class, serves as the
baseline in the paired comparison. Tables~\ref{tab:errors}
and~\ref{tab:paired} and Figure~\ref{fig:cdf} summarise the outcome on all
50 orbits of each potential.

\begin{table}[h]
\centering
\caption{Relative action spreads \eqref{eq:spread} of the four selectors
	\eqref{eq:sel-med}--\eqref{eq:sel-C} and \eqref{eq:sel-sh}, in percent,
	over all 50 orbits of each potential, in the form \emph{median (90th
	percentile)}. Every column is aggregated separately, and the joint spread
	is formed on each orbit before the median and percentile are taken.
	Smaller values indicate more accurate actions.}
\label{tab:errors}
\begin{tabular}{llccc}
\toprule
Potential & Selector & $100s_R$ & $100s_z$ & $100s_{\mathrm{joint}}$\\
\midrule
Logarithmic & $\Delta_{\mathrm{med}}$ & 1.010 (3.093) & 1.086 (4.105) & 1.714 (3.158)\\
 & $\Delta_B$ & 1.118 (3.312) & 1.057 (3.503) & 1.669 (2.908)\\
 & $\Delta_C$ & 0.922 (2.982) & 0.945 (3.161) & 1.460 (2.596)\\
 & $\Delta_{\mathrm{sh}}$ & 0.865 (2.696) & 0.964 (3.350) & 1.454 (2.560)\\
\midrule
Miyamoto--Nagai & $\Delta_{\mathrm{med}}$ & 0.260 (4.098) & 0.170 (0.527) & 0.407 (2.898)\\
 & $\Delta_B$ & 0.236 (1.223) & 0.163 (0.448) & 0.312 (0.907)\\
 & $\Delta_C$ & 0.416 (3.854) & 0.182 (0.530) & 0.523 (2.727)\\
 & $\Delta_{\mathrm{sh}}$ & 0.693 (3.070) & 0.280 (0.907) & 0.811 (2.172)\\
\bottomrule
\end{tabular}
\end{table}

\begin{table}[h]
\centering
\caption{Paired comparison of the selectors $B$, $C$ and $\mathrm{sh}$ with
	the median prescription, orbit by orbit, over all 50 orbits of each
	potential. On each orbit, the \emph{defect ratio} is the mean of
	$b_\Delta^2$ over the validation points at the focal distance of the
	selector divided by the same mean at the focal distance of the median
	prescription, and the \emph{spread ratio} is the corresponding quotient of
	the joint spreads. The table lists the median of each ratio, with the 10th
	and 90th percentiles of the spread ratio in brackets, and a value below 1
	favours the selector. The last two columns count the orbits on which the
	selector reduces the joint spread and both $s_R$ and $s_z$.}
\label{tab:paired}
\begin{tabular}{llcccc}
\toprule
Potential & Selector & Defect ratio & Spread ratio [10th, 90th] & Joint & Both\\
\midrule
Logarithmic & $B$ & 0.160 & 0.992 [0.903, 1.077] & 26/50 & 1/50\\
 & $C$ & 0.210 & 0.891 [0.793, 0.945] & 48/50 & 48/50\\
 & $\mathrm{sh}$ & 0.654 & 0.875 [0.806, 0.939] & 50/50 & 50/50\\
\midrule
Miyamoto--Nagai & $B$ & 0.625 & 0.941 [0.298, 1.086] & 27/50 & 26/50\\
 & $C$ & 0.862 & 1.076 [0.871, 1.455] & 22/50 & 19/50\\
 & $\mathrm{sh}$ & 2.584 & 1.956 [0.686, 2.810] & 12/50 & 12/50\\
\bottomrule
\end{tabular}
\end{table}

\paragraph{Logarithmic potential.} The unweighted residual fit improves on
the baseline almost uniformly. It reduces both action spreads on 48 of the
50 orbits, the exceptions being the two launches with $(f_R,f_z)=(0,0.6)$,
and the median paired ratio of the joint spread is 0.891. The
operator-weighted selector lowers the validation mean of $b_\Delta^2$ on 49
of the 50 orbits, with median ratio 0.160, but it improves the joint spread
on only 26 orbits and both components on one, so in this potential a much
smaller operator defect does not translate into more accurate actions.

\paragraph{Miyamoto--Nagai potential.} Here the ranking is reversed, and the
advantage of the operator-weighted selector lies in the worst cases. Its
median joint spread is close to that of the baseline, with median ratio
0.941 and improvements on 27 orbits, but it reduces the 90th-percentile
joint spread by a factor of about three, from $2.898\%$ to $0.907\%$. The
nine largest baseline spreads occur on strongly vertically excited orbits,
with $f_z\ge0.45$, and are dominated by their radial component, and on these
nine orbits the operator-weighted selector reduces the joint spread by
factors between 2.9 and 4.4. The unweighted residual fit improves the joint
spread on only 22 orbits, with median ratio 1.076.

The operator accounts for both the location of these failures and the
direction of the correction. On the nine orbits, between $59\%$ and $68\%$
of the training points lie above half of the maximal height of the orbit,
where the median local estimate is between 1.15 and 1.37, close to the value $a_d^2=1$ of the Kuzmin limit (Section~\ref{sec:numerics}), whereas closer to the plane it is between 1.62 and
2.01. The median prescription, which treats all samples alike, is therefore
pulled towards the values far from the plane, where the orbit lingers near
its vertical turning points, and it gives $\Delta_{\mathrm{med}}^2$ between
1.21 and 1.44. The operator objective suppresses these turning regions
through the weight $w_\Delta$ and selects the larger values $\Delta_B^2$
between 1.48 and 1.64, which on each of the nine orbits give the more
accurate actions.

The five-input operator also measures the part of the departure from
St\"ackel form that no focal distance removes. By \eqref{eq:fudge-weight},
the minimum $F_B(\delta_B)$ in \eqref{eq:sel-B} is the mean square of the
operator values $b_{\Delta_B}(\mathbf x_m)$ over the training points, and
the \emph{relative residual defect}, introduced here, normalises this
minimum,
\begin{equation}
	\rho_B=\frac{F_B(\delta_B)}{\dfrac1N\displaystyle\sum_{m=1}^Nw_{\Delta_B}(\mathbf x_m)\,A_m^2}.
	\label{eq:rhoB}
\end{equation}
By \eqref{eq:fudge-B}, each summand of the denominator is the square of the
operator value obtained when the momentum factor is kept at $\Delta_B$ and
the residual $C_{\Delta_B}$ is replaced by $C_0=A$. This is the residual of
spherical coordinates, the limit $\Delta\to0$ of the prolate family, since
in spherical coordinates $(r,\theta)$ it reads
$A=r^{-1}\partial_r\partial_\theta(r^2\Phi)$ and vanishes exactly for the
potentials $f(r)+h(\theta)/r^2$. The two means in \eqref{eq:rhoB} therefore
share the momentum factor of the operator and differ only in the residual,
and $\rho_B$ is the fraction of the spherical mean square that the selected
focal distance leaves in place. Its value is nonnegative, and $\rho_B=0$
exactly when the operator vanishes at every training point for one focal
distance, which is the vanishing of the off-diagonal value required by
Proposition~\ref{prop:LC} on an open domain, imposed here only at the
training points. With the momentum factor held at $\Delta_B$, the operator
values depend linearly on $\delta=\Delta^2$ through $C_\Delta=A-g\delta$.
Writing $w_m=w_{\Delta_B}(\mathbf x_m)$, this gives
\[
\rho_B=1-\frac{\bigl(\sum_mw_mg_mA_m\bigr)^2}{\sum_mw_mg_m^2\,\sum_mw_mA_m^2}
\]
whenever $\delta_B$ also minimises the sum with the weights held fixed, and
in that case $\rho_B\le1$, with equality when no focal distance improves on
spherical coordinates. Since $\delta_B$ also shifts the weight, this
expression is approximate for the selected focal distance. All values found
here are below $0.09$, so on every orbit the selected focal distance removes
more than $91\%$ of the spherical mean square of the operator. As a ratio of
two mean squares of operator values, $\rho_B$ is dimensionless, can be
compared across orbits and models, and requires neither actions nor
quadratures. Its rank correlation with the joint spread of the baseline is 0.902 in the
Miyamoto--Nagai model and 0.627 in the logarithmic model, and with the joint
spread of the operator-weighted selector 0.931 and 0.610, respectively. In
the Miyamoto--Nagai model its nine largest values include seven of the nine
worst orbits of the baseline, and in the logarithmic model they belong
exactly to the nine worst orbits, although they are only slightly separated
from the others. The
five-input operator therefore plays two roles in the St\"ackel fudge, since
the minimiser of $F_B$ selects the focal distance and its normalised minimum
ranks the orbits a priori, at no extra cost, by the part of their departure
from St\"ackel form that no focal distance removes.

\paragraph{Shell-orbit fit.} In the logarithmic model Binney's selector is
the most accurate of the four, since it improves on the baseline on all 50
orbits and in both components, with median paired ratio 0.875, and it is
practically equivalent to the unweighted fit. In the Miyamoto--Nagai model
it improves on the baseline on only 12 orbits, among them the nine worst,
and it is less accurate than the operator-weighted selector on every orbit.
The local estimate \eqref{eq:dloc-models} explains this pattern. On the 30
orbits with $f_z\le0.30$, the shell orbit rises between 1.2 and 14 times
higher above the plane than the orbit itself, into the region where
$\delta_{\mathrm{loc,MN}}$ approaches $a_d^2$, so that $\Delta_{\mathrm{sh}}<\Delta_{\mathrm{med}}$ and the baseline is more
accurate on all of them. On these 30 orbits, the shell-orbit focal distance
multiplies the validation mean of $b_\Delta^2$ on the orbit itself by a
median factor of 4.55, against 2.58 over all 50 orbits 
(Table~\ref{tab:paired}).

\paragraph{Launch radius.} The gain of the operator-weighted selector in the
Miyamoto--Nagai model comes mainly from $R_0=2$, with median paired ratios
$1.000$ at $R_0=1$ and $0.775$ at $R_0=2$. In the logarithmic model each
orbit at $R_0=2$ is the image, under a scaling of positions and times by the
factor two, of the orbit with the same $(f_R,f_z)$ at $R_0=1$ in the model
with half the core radius, so the logarithmic grid probes 25 orbit shapes at
two nearly equivalent scales, whose joint spreads differ by at most
$6.2\%$.

\paragraph{Consistency with the reference actions.} The reference actions
of Appendix~\ref{app:reference} pass the acceptance checks on 98 of the 100
orbits. The two Miyamoto--Nagai launches $(R_0,f_R,f_z)=(1,0,0.6)$ and
$(1,0.15,0.6)$ fail the convergence checks, but no inference about chaotic
motion can be drawn from this. On the accepted orbits, the relative RMS errors
of the first three selectors with respect to the reference actions
reproduce the medians and 90th percentiles of the spreads to within $2.2\%$
in the logarithmic model and $5.3\%$ in the Miyamoto--Nagai model, with the
same improvement counts, so the bias is a minor part of the error and the
conclusions do not depend on the reference actions.

\paragraph{Defect reduction and action accuracy.} Each selector is computed
from its own action-free measure of the departure from St\"ackel form,
namely the local estimates of Sanders \cite{Sanders2012} in
\eqref{eq:sel-med}, the operator objective $F_B$ of \eqref{eq:sel-B},
introduced in this paper, the residual criterion $F_C$ of Vasiliev
\cite{VasilievDoc} in \eqref{eq:sel-C} and the shell-orbit ellipse fit of
Binney \cite{Binney2014} in \eqref{eq:sel-sh}, and by
Proposition~\ref{prop:fudge} all four are expressed through the same
operator $b_\Delta$. The accuracy, however, is measured by the action spread
\eqref{eq:spread}, and the relation between the two depends on where the
residual is large relative to the turning points that bound the quadratures
(Section~\ref{sec:weighting}), so none of the four selectors can guarantee
a smaller spread on every orbit. Since the grid was fixed in
advance, the orbits without improvement are reported together with the
successes. Conceptually, the operator places the local prescription of
Sanders, Vasiliev's residual criterion, Binney's shell-orbit fit and a new
operator-weighted selector in a single framework, as the pointwise zero,
the unweighted fit, the zero on the shell orbit and the weighted fit on the
orbit itself of one and the same quantity. Practically, at the same
action-free cost, the unweighted fit and the shell-orbit fit give the better
typical accuracy on this grid in a mildly flattened model such as the
logarithmic one, whereas the operator-weighted fit gives the smaller tail of
large errors, as for the strongly vertically excited orbits of the
Miyamoto--Nagai model, and no weighting dominates in both models.

\subsection{Momentum weighting and action sensitivity}
\label{sec:weighting}

The operator objective $F_B$ of \eqref{eq:sel-B} and Vasiliev's residual
criterion $F_C$ of \eqref{eq:sel-C} are sampled squares of the same
residual $C_\Delta$ of \eqref{eq:residual} over the same training points,
and by \eqref{eq:fudge-weight} they differ only by the momentum weight
$w_\Delta$,
\[
F_B(\delta)=\frac1N\sum_{m=1}^Nw_{\sqrt\delta}(\mathbf x_m)\,C_{\sqrt\delta}(\mathbf x_m)^2,
\qquad
F_C(\delta)=\frac1N\sum_{m=1}^NC_{\sqrt\delta}(\mathbf x_m)^2 .
\]
This weight is produced by the five-input operator itself. The derivatives
of the Hamiltonian \eqref{eq:fudge-H} with respect to the prolate momenta
are the coordinate velocities $\dot u=H_{p_u}=D^{-1}p_u$ and
$\dot v=H_{p_v}=D^{-1}p_v$, so that \eqref{eq:fudge-B} and
\eqref{eq:fudge-weight} take the form
\begin{equation}
	b_\Delta=-\dot u\,\dot v\,C_\Delta,\qquad w_\Delta=\dot u^2\dot v^2,
	\label{eq:velocity-factor}
\end{equation}
and the operator multiplies the residual by the product of the two
coordinate velocities. The weight therefore vanishes exactly at the turning
points of either prolate coordinate, where the St\"ackel quadratures have
their endpoints. At these points the hypothesis $H_{p_i}\ne0$ of
Proposition~\ref{prop:LC} fails, so $F_B$ gives no weight to the points at
which the momentum-graph form of the separability criterion does not apply,
whereas $F_C$ counts them like any other sample. The sensitivity of the
actions to the parameters of the quadratures near these endpoints is
investigated below, together with the share of the residual that each
objective assigns to the turning regions, and the results show that $F_B$
fits the residual away from the turning points, whereas $F_C$ gives the
turning regions at least their share of the samples.

Each St\"ackel action is a quadrature of the form \eqref{eq:action-variation}
below, with $x$ one of the prolate coordinates, $p_x$ its conjugate
momentum, $P(x,s)=p_x^2$ the square of this momentum, written through the
separated equation as a function of $x$, and $s$ a parameter of the
quadrature, such as the focal distance. If $P(x,s)$ is of class $C^2$ and $a(s)<b(s)$ are simple zeros of
$P(\cdot,s)$ with $P>0$ between them, the substitution $x=a+(b-a)\,t$ makes
the $s$-derivative of the integrand bounded, and differentiation under the
integral sign, followed by an integration by parts whose boundary terms
vanish with $\sqrt P$, shows that
\begin{equation}
	J(s)=\frac1\pi\int_{a(s)}^{b(s)}\sqrt{P(x,s)}\,\dd x
	\qquad\text{satisfies}\qquad
	J'(s)=\frac1{2\pi}\int_{a(s)}^{b(s)}\frac{\partial_sP(x,s)}{\sqrt{P(x,s)}}\,\dd x .
	\label{eq:action-variation}
\end{equation}
Since $\sqrt{P(x,s)}=|p_x|$, the kernel $P^{-1/2}=|p_x|^{-1}$ of
\eqref{eq:action-variation} is the inverse of the separated momentum, and
it is largest near the turning points, where the operator weight
\eqref{eq:velocity-factor} vanishes. The kernel is integrable,
however, and a neighbourhood of width $\ell$ of a simple zero contributes
only $O(\sqrt\ell)$ to $J'(s)$ when $\partial_sP$ is bounded, so
\eqref{eq:action-variation} alone does not imply that the endpoints dominate
the response of the actions to the focal distance.

The share of the residual carried by the turning regions is measured on the
near-turning set
\begin{equation}
	T_\eta=\Bigl\{m:\min\Bigl(\frac{|p_{u,m}|}{\sigma_u},\frac{|p_{v,m}|}{\sigma_v}\Bigr)\le\eta\Bigr\}
	\label{eq:near-turning}
\end{equation}
of the validation points, where $\sigma_u$ and $\sigma_v$ are the RMS values
of $p_u$ and $p_v$ over the training points and $\eta>0$ is a threshold. At
the focal distance of the operator-weighted selector and $\eta=0.1$, the
set $T_\eta$ makes up a median $8.59\%$ and $9.38\%$ of the samples in the
logarithmic and Miyamoto--Nagai models, respectively, but it carries only
$0.044\%$ and $0.033\%$ of the weight $w_\Delta$ and $0.161\%$ and $0.118\%$
of the squared operator $b_\Delta^2$, which is the summand of $F_B$, against
$15.16\%$ and $15.08\%$ of the unweighted squared residual $C_\Delta^2$,
which is the summand of $F_C$. The same ordering holds for $\eta=0.05$ and
$\eta=0.2$, for the other selectors and on the training samples.

The squared operator thus concentrates the fit on the part of the orbit
where both coordinate velocities are bounded away from zero and
Proposition~\ref{prop:LC} applies, and \eqref{eq:velocity-factor} identifies
this concentration as a mechanism through which $F_B$ and $F_C$ can rank
orbits differently. Its effect has opposite signs in the two models. The
median ratio of the joint spread of the unweighted fit to that of the
operator-weighted selector is $0.901$ in the logarithmic model, where the
unweighted fit is more accurate on 48 of the 50 orbits, and $1.492$ in the
Miyamoto--Nagai model, where the operator-weighted selector is more accurate
on 36 of them. These include the nine worst orbits of
Section~\ref{sec:results}, on which the suppression of the vertical turning
regions raises the selected focal distance above the median prescription
and the joint spread of the unweighted fit is between 2.4 and 3.6 times
larger. The scaled focal displacement $|\Delta_B-\Delta_C|/R_0$ alone does not
predict which objective gives the smaller joint spread on a given orbit,
since its rank correlation with the same ratio of joint spreads, that of the
unweighted fit to the operator-weighted selector, is $-0.562$ in the
logarithmic model and $0.822$ in the Miyamoto--Nagai model, so that a larger
displacement favours the unweighted fit in the first model and the
operator-weighted selector in the second. A quantitative
criterion would require control of how the residual $C_\Delta$ perturbs each
separated momentum, near and away from the turning points, and this is left
for future work.

\section{Limitations}
\label{sec:limitations}

The results of Sections~\ref{sec:operator}--\ref{sec:stackel} rest on
assumptions that restrict their scope, and the resulting limitations are
summarised in the following list.
\begin{itemize}
	\item \emph{Scope of the theory.} All statements are local and refer to a
	prescribed product structure. Separation requires the off-diagonal values
	to vanish on an open set, and pointwise zeros, the zero set of a single
	scalar evaluation and the full kernel of Proposition~\ref{prop:covariance}
	are different objects. The diagonal entries are not invariant and carry no
	separability information.
	\item \emph{Generality.} The theory applies to arbitrary smooth
	Hamiltonians, and the defect is computed in closed form for the non-natural
	polar-nodal normal forms of satellite theory
	(Proposition~\ref{prop:polar-nodal}) and for Kerr geodesics
	(Proposition~\ref{prop:kerr}). The coordinate-selection applications,
	however, concern natural Hamiltonians, for which algebraic criteria
	\cite{Waksjo,Benenti2016} also find exactly separating coordinates, and
	selection problems for Hamiltonians with momentum-linear or non-quadratic
	kinetic terms remain to be developed.
	\item \emph{Defects and accuracy.} The defect measures the departure from
	separation, which is a different quantity from trajectory and action
	errors, and its pointwise values carry no stability information. The
	oscillator of Example~\ref{ex:oscillator} has $b_{12}=-\epsilon p_1p_2\ne0$
	in the original pairs but is linearly stable, whereas
	$\tfrac12(p_1^2-q_1^2)+\tfrac12(p_2^2+q_2^2)$ has zero defect and an
	exponentially growing tangent direction. Only the integral of the defect
	along a homoclinic orbit carries dynamical information
	(Corollary~\ref{cor:dichotomy}), it does so at first order in the
	perturbation, and the layer width \eqref{eq:width} derived from it is a
	heuristic estimate.
	\item \emph{Numerical comparison.} The 100 launches form a deterministic
	grid in two potentials, so the medians, percentiles and counts describe
	that grid and carry no population-level or causal interpretation, and the
	logarithmic grid contains only 25 distinct orbit shapes, at two nearly
	equivalent scales. The action spread \eqref{eq:spread} is blind to a bias
	common to all points of an orbit, and the reference actions that show this
	bias to be small on 98 of the 100 orbits pass stated convergence checks but
	do not come with certified error bounds. The shell-orbit fit is implemented
	as a linear least-squares ellipse fit to one passage of the shell orbit, and
	other implementations of Binney's prescription may give slightly different
	focal distances. The ranking of the orbits by the relative residual defect
	\eqref{eq:rhoB} is an empirical finding on this grid, and $\rho_B=0$ is only
	a necessary condition for St\"ackel form at the training points.
	\item \emph{Mechanism.} The velocity factor \eqref{eq:velocity-factor} and
	the near-turning shares of Section~\ref{sec:weighting} show how differently
	the two objectives weight an orbit, but they neither predict which selector
	is more accurate on a given orbit nor provide a bound on the action error.
\end{itemize}

These limitations also indicate the natural next steps. A quantitative
relation between the residual $C_\Delta$ and the error of each separated
momentum would turn the weighting comparison into a prediction, selectors
that interpolate between $F_B$ and $F_C$ could combine the strengths of
both, and the linearised operator $\Lop_{H_0}$ could be used to classify the
perturbations that preserve separation to first order.

\section{Conclusions}
\label{sec:conclusions}

This paper has introduced the five-input operator $\B_{a,b;c,d}[H]$, a
scalar built from Poisson brackets with the Hamiltonian, and has used it as
one calculus for Hamiltonian separability, for the obstructions that
perturbations create and for the selection of approximately separating
coordinates. On a prescribed symplectic product structure, its off-diagonal
values reproduce the Levi-Civita expressions
\cite{LeviCivita,LlibreRamirez2020} (Proposition~\ref{prop:canonical}). They
are shown here to equal a mixed Hessian and a symplectic pairing of sector
fields (Proposition~\ref{prop:intrinsic}), which explains their symmetry and
their invariance under canonical changes within the factors, and to express
the Frobenius obstruction of the Lagrangian distribution spanned by the
sector fields by one invariant scalar per pair of degrees of freedom
(Proposition~\ref{prop:frobenius}). Their vanishing yields a local
separation theorem that requires only nonvanishing sector fields and
therefore also holds at the turning points of the separated motions, where
the momentum-graph criterion of Proposition~\ref{prop:LC} fails
(Theorem~\ref{thm:sector}). The same scalar takes the form of a bordered
Hessian \eqref{eq:bordered}, and evaluated on the angles of an integrable
Hamiltonian it reproduces the isoenergetic nondegeneracy determinant of
Arnold (Remark~\ref{rem:arnold}).

Because the operator is cubic in the Hamiltonian, perturbations of
separable systems have exact, terminating defect expansions with an
explicit linearised operator, whose kernel contains the symmetry directions
of the unperturbed system (Proposition~\ref{prop:kernel}). For Hamiltonians
with a cyclic angle, the linearised operator is an angle derivative of the
Melnikov integrand, so that the defect integrated along a homoclinic orbit
decides the first-order splitting of the separatrix
(Theorem~\ref{thm:melnikov} and Corollary~\ref{cor:dichotomy}). The defect
is obtained in closed form for the polar-nodal normal forms of satellite
theory and for Kerr geodesics and their deformations
(Propositions~\ref{prop:polar-nodal} and~\ref{prop:kerr}).

Used as a coordinate-selection tool, the defect identifies the separating
focal distance of the Euler problem exactly (Theorem~\ref{thm:euler}). For
axisymmetric galactic potentials, the operator on prolate coordinates
factorises into the product of the two coordinate velocities and the
St\"ackel residual (Proposition~\ref{prop:fudge} and
\eqref{eq:velocity-factor}). This factorisation places the local
prescription of Sanders \cite{Sanders2012}, the residual criterion of
Vasiliev \cite{VasilievDoc} and the shell-orbit fit of Binney
\cite{Binney2014} in one formulation, as the pointwise zero, the unweighted
fit and the zero on the shell orbit of the same quantity, next to the
operator-weighted selector introduced here. The velocity factor identifies
a mechanism behind the different accuracy of the weighted and unweighted
fits in the two test models, and the relative residual defect
\eqref{eq:rhoB}, the normalised minimum of the operator objective, singles
out most of the least accurate orbits before any action is computed. None
of the four selectors is the most accurate on every orbit, and the
contribution of the operator to the St\"ackel fudge lies in this common
formulation and in the a priori diagnostic.

The identity of Theorem~\ref{thm:melnikov} applies to every separatrix of a
Hamiltonian with a cyclic angle, and natural applications are the
separatrices of the hierarchical three-body problem and the zoom--whirl
separatrices of deformed Kerr geodesics \cite{LevinPerezGiz2009}. Within the
local and first-order limits listed in Section~\ref{sec:limitations}, the
five-input operator thus offers one computable quantity that connects the
geometric criteria of separation, the first-order obstructions created by
perturbations and the practical selection of approximately separating
coordinates.

\appendix
\section{Reference-action protocol}
\label{app:reference}

This appendix describes how the reference actions used for the consistency
check of Section~\ref{sec:results} are computed, which checks they must pass, how the
protocol was validated on an exactly separable model, and the fallback used
for four launches.

\paragraph{Generating-function fit.} Let $J^{\mathrm{toy}}_\ell$ and
$\theta^{\mathrm{toy}}_\ell$, $\ell\in\{R,z\}$, be the exact actions and
angles of an isochrone toy potential, evaluated at the integrated physical
samples. The two action equations
\begin{equation}
 J^{\mathrm{toy}}_\ell(t)=J^{\mathrm{ref}}_\ell
 +\sum_{\mathbf n\in\mathcal N_K}n_\ell A_{\mathbf n}
 \cos\bigl(n_R\theta^{\mathrm{toy}}_R(t)+n_z\theta^{\mathrm{toy}}_z(t)\bigr)
 \label{eq:gf-fit}
\end{equation}
are fitted jointly, where $\mathcal N_K$ contains one representative of each
pair $\pm\mathbf n$ with $0<n_R^2+n_z^2\le K^2$ and $n_z$ even. The
coefficients $A_{\mathbf n}$ are shared by the two equations, which enforces
the generating-function structure, and the azimuthal action is the known
$L_z$. The isochrone length is $0.8R_0$, and its amplitude is matched to the
launch circular speed and doubled, if necessary, until all sampled toy
energies are negative.

\paragraph{Acceptance.} The selected reference uses $K=16$. It is compared
componentwise with the fit at $K=12$, with the fits to the two temporal
halves at $K=16$ and with a fit using the toy length $0.9R_0$, and it is
accepted if all of the following conditions hold:
\begin{enumerate}[label=(\roman*)]
 \item the reference actions are positive;
 \item the fit matrix has full rank and a condition number below $10^8$;
 \item each retained Fourier phase covers at least $2\pi$, with successive
 increments below $\pi$;
 \item the largest relative change in either action over the comparison fits
 is at most $2\times10^{-4}$.
\end{enumerate}
Insufficient phase resolution triggers only a refinement of the reference
cadence, capped at 1024 samples per $T_c$, and leaves the training and
validation samples unchanged. In addition, the energy drift, scaled by
$\max\{|H(0)|,v_c(R_0)^2\}$, must stay below $10^{-8}$, and the 64-point
St\"ackel actions must agree with the 96-point quadrature to $10^{-4}$ at 32
validation points per orbit and selector.

\paragraph{Controls.} As a control, the protocol is applied to the exactly
separable flattened Kuzmin--Kutuzov potential \cite{DejongheDeZeeuw,galpyKK}
with axis ratio $5$ of the coordinate surfaces and focal distance $\Delta=0.5$,
normalised to $v_c(1)=1$, along the orbit launched from
$(R_0,f_R,f_z)=(1,0.3,0.3)$, on which both toy actions vary. The order-16
reference differs from adaptive separated quadratures by $1.1\times10^{-8}$
radially and by $5.3\times10^{-9}$ vertically, and the order-24 reference by
$7.6\times10^{-12}$ and $5.1\times10^{-11}$, below the variation of the
quadrature actions along the orbit, which is about $10^{-10}$.

\paragraph{Fallback for circular toy launches.} Four launches have $f_R=0$
and $f_z=0.6$, so that their launch speed equals $v_c(R_0)$, the spherical
toy orbit is circular at launch and its radial angle is undefined, and all
four failed the phase checks at the cadence cap. A fallback that changes
only the reference chart was therefore introduced at the revision stage,
after a first amplitude multiplier of 1.2 had failed: the matched amplitude
is multiplied by 1.02, and the toy lengths $0.9R_0$ and $R_0$ are used for
the selected and comparison fits at 1024 samples per $T_c$, with the
acceptance thresholds, physical samples, selectors and St\"ackel actions
unchanged. The fallback recovers both logarithmic cases and the
Miyamoto--Nagai case $(2,0,0.6)$, whereas for $(1,0,0.6)$ the reference
discrepancy still exceeds the threshold. All attempts are retained in the
case inventory of the archive, and since the reference actions serve only
for the consistency check, the fallback affects none of the spreads in
Tables~\ref{tab:errors} and~\ref{tab:paired}.


\section*{CRediT authorship contribution statement}

\textbf{Vladimir Martinusi:} Conceptualization, Methodology, Formal analysis,
Software, Validation, Investigation, Visualization, Writing -- original
draft, Writing -- review \& editing.

\section*{Declaration of competing interest}

The author declares no known competing financial interests or personal
relationships that could have appeared to influence the work reported in
this paper.

\section*{Funding}

This research did not receive any specific grant from funding agencies in
the public, commercial, or not-for-profit sectors.

\section*{Data availability}

The source archive \cite{CodeArchive} contains the integration, shell-orbit
and reference-fitting scripts, the launch-state specification, every
attempted case with its acceptance diagnostics, the training and validation
arrays, the symbolic identity checks carried out with SymPy
\cite{SymPy2017} (including those of Sections~\ref{sec:homoclinic}
and~\ref{sec:levelset}), the separatrix integrations of
Example~\ref{ex:rotor}, the exact-model controls, the action spreads, paired
statistics, relative residual defects and turning-point diagnostics of both
models, the software versions together with input and output hashes, and
the scripts that generate all tables and figures. The archive is deposited
in Zenodo under the MIT licence and is openly available at
\url{https://doi.org/10.5281/zenodo.24048471}.

\section*{Declaration of generative AI and AI-assisted technologies in the
manuscript preparation process}

During the preparation of this work the author used Claude (Anthropic) in
order to check derivations symbolically, write and run the reproduction
code, and edit the language of the manuscript. After using this
tool/service, the author reviewed and edited the content as needed and
takes full responsibility for the content of the published article.

\bibliographystyle{unsrtnat}
\bibliography{references}

\end{document}